\documentclass[11pt]{amsart}
\usepackage{fullpage}
\usepackage{amsfonts,amssymb,amsmath,amsthm}
\usepackage{hyperref}
\hypersetup{pdftitle={Injectivity and stability for the sub-Riemannian X-ray transform on H-type groups}, pdfauthor={Steven Flynn}, pdfkeywords={X-ray transform, sub-Riemannian geometry, H-type groups, group Fourier transform, conjugate points, limited-data tomography}}
\usepackage{color}
\usepackage[normalem]{ulem} %
\usepackage{comment}
\usepackage{tikz-cd}
\usepackage[backend=bibtex,date=year,maxbibnames=99,sortcites=true]{biblatex}
\DeclareFieldFormat{doi}{} %
\DeclareFieldFormat{isbn}{} %
\DeclareFieldFormat{issn}{} %
\DeclareFieldFormat{url}{} %
\DeclareFieldFormat{urldate}{} %

\DeclareMathOperator{\sgn}{sgn}
\theoremstyle{plain}
\newtheorem{theorem}{Theorem}[section]
\newtheorem*{theorem*}{Theorem}
\newtheorem{thmA}{Theorem}

 \newtheorem{corollary}[theorem]{Corollary}
 
 \newtheorem{lemma}[theorem]{Lemma}
 \newtheorem{proposition}[theorem]{Proposition}
 \theoremstyle{definition}
 \newtheorem{definition}[theorem]{Definition}
 
 \theoremstyle{remark}
 \newtheorem{remark}[theorem]{Remark}
 \numberwithin{equation}{section}
  \newtheorem{example}[theorem]{Example}

\def\fz{\mathfrak{z}}
\def\fg{\mathfrak{g}}
\def\fv{\mathfrak{v}}

\def\cA{\mathcal{A}}

\def\cE{\mathcal{E}}
\def\cK{\mathcal{K}}
\def\cN{\mathcal{N}}
\def\cJ{\mathcal{J}}
\def\cU{\mathcal{U}}
\def\cG{\mathcal{G}}
\def\cD{\mathcal{D}}
\def\0s{\setminus\{0\}}

\def\cH{\mathcal{H}}
\def\cF{\mathcal{F}}
\def\cL{\mathcal{L}}
\def\ddt{\frac{d}{dt}\big|_{t=0}}

\def\CP{\mathbb{CP}}

\def\mz{\setminus\{0\}}

\def\eps{\varepsilon}
\def\R{{\mathbb R}}%
\def\C{{\mathbb C}}%
\def\N{{\mathbb N}}%
\def\Z{{\mathbb Z}}%
\mathchardef\mhyphen="2D %

\def \heis {\mathbb{H}}

\numberwithin{equation}{section}

\begin{document}
	\title[Injectivity and stability for the sub-Riemannian X-ray transform on H-type groups]{Injectivity and stability for the sub-Riemannian X-ray transform on H-type groups}

	\author[S. Flynn]{Steven Flynn}
	\address[S. Flynn]{Department of Mathematics, University College London, Gower Street, London WC1E 6BT, United Kingdom}
	\email{steven.flynn@ucl.ac.uk}

\begin{abstract}
	We prove injectivity and a stability estimate for the X-ray transform along sub-Riemannian geodesics of $H$-type groups, a setting in which almost every geodesic carries conjugate points. The geodesics are helices labeled by a charge $\lambda$ valued in the dual of the center. We prove a Fourier slice theorem, which identifies the X-ray transform over geodesics with charge $\lambda$ as the restriction of the group Fourier transform to their \textit{visible frequencies}, followed by an explicit multiplier. We characterize the sets of charges whose geodesics determine an integrable function, as exactly those whose visible frequencies are dense. In particular, when the center of the group has dimension at least two, charges of a single magnitude suffice provided their directions fill a great circle. Laguerre asymptotics on the spectrum of the normal operator at fixed charge yields a stability estimate at every nonzero visible frequency with the loss of exactly half a derivative in the sub-Laplacian, except on the first Heisenberg group. Finally, we prove that data limited to a band or a cap of charges determines all sufficiently high central frequencies of an unknown function.
\end{abstract}

\maketitle

\makeatletter
\renewcommand\l@subsection{\@tocline{2}{0pt}{3pc}{5pc}{}}
\makeatother

\setcounter{tocdepth}{1}
\tableofcontents

\section{Introduction}

Let $G$ be an $H$-type group (defined in Section \ref{H-type}), the model case of step-two sub-Riemannian geometry, with its left-invariant sub-Riemannian metric, and let $\cG$ be the set of maximal oriented sub-Riemannian geodesics in $G$. The \textit{X-ray transform} of a function $f$ on $G$ is
\begin{align*}
	I(f)(\gamma):=\int_\gamma f\, ds, \quad \gamma\in \cG.
\end{align*}
This paper asks which subfamilies of $\cG$ determine $f$, and how stably.

The Lie algebra of $G$ is $\fg=\fv\oplus\fz$, with $\fz$ the center (Section \ref{H-type}). For the first Heisenberg group $\heis_1=\R^3$, the subspace $\fv$ is the horizontal plane and $\fz$ the vertical axis. Every geodesic has a \textit{charge} $\lambda\in \fz^*$, the central component of its momentum, a constant of motion (Section \ref{sRGeometry}). A geodesic of charge $\lambda\neq0$ is a helix \eqref{param} over a circle of radius $|\lambda|^{-1}$, advancing in the direction $\hat\lambda=\lambda/|\lambda|$ along its central axis by $\pi|\lambda|^{-2}$ per turn. For $\lambda\in \fz^*\mz$ let
\begin{align*}
	\Gamma_\lambda^*:=\{\mu\in \fz^*: \mu\cdot\hat\lambda\in 2|\lambda|^2\Z\},
\end{align*}
the \textit{visible frequencies} of charge $\lambda$, the $\mu$ for which $e^{i\mu\cdot u}$ is unchanged by the central advance. They form parallel hyperplanes orthogonal to $\lambda$ with spacing $2|\lambda|^2$. On $\heis_1$ they are the points $2\lambda^2\Z$.
\begin{thmA}\label{t:intro}
	Let $G$ be an $H$-type group and $Z\subset \fz^*\mz$. The following are equivalent.
	\begin{itemize}
		\item[(a)] If $f\in L^1(G)$ and $I(f)(\gamma)=0$ for almost every geodesic $\gamma$ with charge $\lambda\in Z$, then $f=0$.
		\item[(b)] The set of visible frequencies $\bigcup_{\lambda\in Z}\Gamma_\lambda^*$ is dense in $\fz^*$.
	\end{itemize}
	In particular $I$ is injective on $L^1(G)$. If $\dim\fz=1$ then (a) holds if and only if $Z$ accumulates at zero, that is, if and only if the radii of the helices are unbounded. If $\dim\fz\geq2$ then (a) holds for $Z=\{\rho\theta: \theta\in C\}$, for every $\rho>0$ and every great circle $C$ of the unit sphere of $\fz^*$.
\end{thmA}

 The geodesics of charge $\lambda\neq0$ are the left translates of the model helices $\gamma_{\nu, \lambda}$ of \eqref{param}, with \textit{horizontal direction} $\nu\in S(\fv^*)$ the initial horizontal component of momentum. The transform then splits into the family of transforms $I_{\nu, \lambda}f$, the integrals of $f$ over these translates (Section \ref{s:geodesicsXray}). Theorem \ref{t:intro} follows from Theorem \ref{MainResult} and Corollary \ref{c:injSets}. In fact, one horizontal direction per charge suffices (Theorem \ref{t:singleDirection}).

The transforms of charge $\lambda$ recover the group Fourier transform $\cF(f)(\mu)$ of $f$ only at frequency $\mu\in \Gamma_\lambda^*$ (Theorem \ref{FourierSliceTheorem}). %

For $G\neq\heis_1$, at every nonzero frequency $\mu\in\Gamma_\lambda^*$, the data $I_{\nu, \lambda}f$ for all directions $\nu$ at charge $\lambda$ determine $\cF(f)(\mu)$ with the loss of exactly half a derivative in the sub-Laplacian, as for the Euclidean X-ray transform \cite{strichartzLpEstimatesRadon1981} (Theorem \ref{t:stability}). The stability estimate holds despite the fact that every geodesic of nonzero charge carries conjugate points.

Finally, a band of charges ($\dim\fz=1$) or a spherical cap ($\dim\fz\geq2$) determines $\cF(f)(\mu)$ for $|\mu|$ above an explicit threshold (Theorem \ref{PW2}), and for $f$ of compact central support the rest by analytic continuation (Corollary \ref{c:centralSupport}).

Geodesic X-ray transforms are studied extensively in Riemannian geometry \cite{sharafutdinovIntegralGeometryTensor1994, ilmavirtaIntegralGeometryManifolds2018, paternainGeometricInverseProblems2023}, for noninvasive imaging \cite{nattererMathematicsComputerizedTomography1986} and as linearizations of boundary and spectral rigidity \cite{stefanovRigidityMetricsSame1998, guilleminInverseSpectralResults1980}, and on groups and symmetric spaces with a left-invariant Riemannian metric \cite{helgasonIntegralGeometryRadon2011, ilmavirtaRadonTransformsCompact2016, ilmavirtaTorusComputedTomography2020a, kleinFunkTransformCompact2009}. Tomographic data are often restricted to a subfamily of curves \cite{gelfandDifferentialFormsIntegral1969, greenleafNonlocalInversionFormulas1989}. For curves tangent to a distribution $\cD\subset TM$, an integrable $\cD$ gives a family of leafwise problems, as in single-axis tilt electron tomography \cite{quintoLocalTomographyElectron2008} or parallel-beam computed tomography \cite{nattererMathematicsComputerizedTomography1986}. More generally, the same reduction applies when $\cD$ and its iterated brackets span a distribution of constant rank less than $\dim M$. If $\cD$ is bracket generating there is no such reduction, and if $\cD$ carries a metric the horizontal curves of stationary energy are the \textit{sub-Riemannian geodesics}.

Since $[\fv, \fv]=\fz$, the subspace $\fv\subset\fg=T_eG$ determines a left-invariant \textit{horizontal} distribution $\cD_q:=dL_q(\fv)$, $q\in G$, which is bracket generating. By Chow's theorem, any two points of $G$ can be joined by a horizontal curve \cite{montgomeryTourSubriemannianGeometries2006}. A sub-Riemannian geodesic is a horizontal curve that locally minimizes length, measured by the inner product on $\fv$, among horizontal curves with the same endpoints. Every geodesic is normal, that is, the projection of a trajectory of the Hamiltonian flow of Section \ref{HamFormalism}.

The simplest case is the $\lambda=0$ geodesics of $\heis_1$, a three-parameter family of lines in $\R^3$, along which the transform is not injective \cite{strichartzLpHarmonicAnalysis1991}; see \cite{rubinRadonTransformHeisenberg2012, peyerimhoffXrayTransform2step2016} for Radon and X-ray transforms on two-step nilpotent groups, with a left-invariant Riemannian metric. To the author's knowledge, X-ray tomography along sub-Riemannian geodesics is unexplored aside from \cite{flynnInjectivityHeisenbergXray2021, flynnSingularValueDecomposition2023}; the first proves the Fourier slice theorem on $\heis_1$ and deduces injectivity on $L^1(\heis_1)$ from the geodesics of all charges in a neighborhood of zero. New here are the $H$-type groups of every dimension, on which the direction $\nu$ is a parameter and, for $\dim\fz\geq2$, the charge and the central frequency need not be parallel. Also new are the characterization of the injectivity sets and the stability estimate.

Under the projection $G\to G/\exp(\fz)\cong\fv$, the geodesics are the horizontal lifts of the trajectories of a charged particle in $\fv$ moving in the constant generalized magnetic field $\omega$ of \eqref{product}. The charge $\lambda$ determines the Lorentz force $\omega_\lambda=\lambda\circ\omega$ (Section \ref{sRGeometry}). The X-ray transform along magnetic geodesics is studied in \cite{dairbekovBoundaryRigidityMagnetic2007}, and for a constant field in $\R^3$ in \cite[Section 10.3]{stefanovGeodesicXrayTransform2012}. The Fourier transform in the central variable turns $I$ into a family of attenuated magnetic ray transforms on $\fv$ \cite{holmanAttenuatedGeodesicXray2018, paternainGeometricInverseProblems2023} (Remarks \ref{r:twistedSlice} and \ref{r:parallelTransport}).

Conjugate points are an important obstruction to well-posedness of geometric inverse problems. Their existence inhibits stable inversion on Riemannian manifolds \cite{holmanAttenuatedGeodesicXray2018}, with complete loss of stability in two dimensions \cite{monardGeodesicRayTransform2015}. On $H$-type groups conjugate points are ubiquitous. On a geodesic $\gamma$ of charge $\lambda\neq0$, translates of a point $\gamma(s)$ by the lattice $\Gamma_\lambda$ are all conjugate \eqref{e:helicalSym}. The lattice $\Gamma_\lambda$ plays the role of the antipodal map on the sphere, where every great circle through $x$ passes through $-x$ (see \cite[Section 10.2]{stefanovGeodesicXrayTransform2012}).
The periodicity of the conjugate points is what makes the spectrum of the normal operator computable.

Our methods rest on the unitary representation theory of $H$-type groups. In the Riemannian setting stability is read off from the fact that the normal operator $I^*I$ is an elliptic pseudodifferential operator of order $-1$ \cite{stefanovBoundaryRigidityStability2005a}, a statement which may fail in the presence of conjugate points \cite{stefanovGeodesicXrayTransform2012}. On a graded nilpotent group the left-invariant pseudodifferential operators of the calculus of \cite{fischerQuantizationNilpotentLie2016} are Fourier multipliers. Our Fourier Slice Theorem \ref{FourierSliceTheorem} fits exactly in this framework, in the same way one treats the classical Radon transform with the classical Fourier Transform \cite{helgasonIntegralGeometryRadon2011}. However, the full normal operator $I^*I$ is undefined with respect to Lebesgue measure on the set of charges. Instead, we compute the eigenvalues of the normal operator for the transform along geodesics of fixed charge $\lambda$, the average of $I_{\nu, \lambda}^*I_{\nu, \lambda}$ over the horizontal directions $\nu$. Its symbol $\cN_\lambda$, which is operator-valued, is a joint spectral multiplier of the sub-Laplacian $\Delta_{\rm{sr}}=\sum_j X_j^2$ and the central derivatives. A pointwise-in-$\mu$ analysis of Laguerre asymptotics shows that, for $G\neq\heis_1$, the symbol agrees to leading order in the spectral parameter with a multiple of the inverse square root of the sub-Laplacian (Theorem \ref{t:stability}, Lemma \ref{l:eigenAsymp}), the sub-Riemannian counterpart of the Euclidean identity $I^*I=c_n(-\Delta)^{-1/2}$.

Sections \ref{H-type}--\ref{s:geodesicsXray} define $H$-type groups and the transforms $I_{\nu, \lambda}$, and Section \ref{MainResults} states the main results. Sections \ref{sRGeometry}--\ref{GFT} develop the geodesic flow, the quotient $G_\lambda=G/\Gamma_\lambda$, the representation theory of $G$ and $G_\lambda$ and the Fourier slice theorem, Section \ref{s:proofs} computes the spectral decomposition of the multipliers and proves the remaining main theorems, and Section \ref{s:concluding} discusses the horizontal lines and the outlook; Appendix \ref{momentumFunctions} collects the momentum functions.

\subsection{$H$-type groups}\label{H-type}
$H$-type groups, first introduced by Kaplan \cite{kaplanFundamentalSolutionsClass1980, kaplanRiemannianNilmanifoldsAttached1981}, are the simplest class of nilpotent groups beyond the classical Heisenberg groups. The Iwasawa $N$-groups among them (Example \ref{ex:htype}) are the boundary at infinity, minus a point, of the rank one symmetric spaces of noncompact type aside from the real hyperbolic spaces \cite{follandHardySpacesHomogeneous1982b, cowlingHTypeGroupsLwasawa}, where the sub-Riemannian structure is the natural one induced on the boundary, and their nilmanifolds are sufficiently diverse to construct the first examples of closed isospectral Riemannian manifolds that are not locally isometric \cite{gordonIsospectralClosedRiemannian1993}.

\begin{definition}\label{d:htype}
	Let $\fg$ be a real Lie algebra with an inner product $\langle\cdot, \cdot\rangle$, let $\fz$ be its center and $\fv$ the orthogonal complement of $\fz$, so that $\fg=\fv\oplus\fz$, and suppose that $\fg$ is 2-step, $0\neq[\fg, \fg]\subset\fz$. For $\mu\in \fz^*$ let $\omega_\mu$ be the skew bilinear form on $\fv$ and $J_\mu: \fv\to\fv$ the skew-adjoint map with
	\begin{align*}
		\omega_\mu(V, W):=\mu([V, W])=\langle J_\mu V, W\rangle, \qquad V, W\in \fv.
	\end{align*}
	The pair $(\fg, \langle\cdot, \cdot\rangle)$ is of \textit{$H$-type} if
	\begin{align*}
		J_\mu^2=-|\mu|^2I \qquad \text{for all } \mu\in \fz^*,
	\end{align*}
	and the \textit{$H$-type group} $G$ is the connected, simply connected Lie group with Lie algebra $\fg$.
\end{definition}
The inner product identifies $\fg^*$ with $\fg$, so that $|\mu|$ and $\mu\cdot\lambda$ denote the induced norm and inner product on $\fz^*$, and $J_\lambda$ acts on vectors $V\in \fv$ and on covectors $\nu\in \fv^*$ alike. The $H$-type condition says that $|\mu|^{-1}J_\mu$ is an orthogonal complex structure on $\fv$ for every $\mu\neq0$. In particular $\omega_\mu$ is non-degenerate for $\mu\neq0$, $[\fv, \fv]=\fz$, and $\dim\fv$ is even. We always set $2n:=\dim\fv$ and $m:=\dim\fz$.

\begin{remark}
	Polarizing the $H$-type condition gives
	\begin{equation}\label{e:htypePolar}
		J_\mu J_\lambda+J_\lambda J_\mu=-2(\mu \cdot \lambda) I, \qquad \mu, \lambda\in \fz^*,
	\end{equation}
	and conversely. Thus $\mu\mapsto J_\mu$ extends to a representation on $\fv$ of the Clifford algebra of $(\fz^*, |\cdot|^2)$, and $H$-type algebras with center $\fz$ correspond to Clifford modules over $\fz^*$ \cite{kaplanFundamentalSolutionsClass1980}. Every $m$ occurs. For $m\geq2$ and orthonormal $\mu, \lambda\in \fz^*$, the anticommuting complex structures $J_\mu$ and $J_\lambda$ generate a quaternionic structure on $\fv$, so $n$ is even.
\end{remark}

\begin{example}\label{ex:htype}
	In each row $\fv$ is a complex, quaternionic or octonionic vector space with its standard inner product, $\fz$ is a space of imaginary units, and $J_\mu$ is left multiplication by $\mu$, or by $i\mu$ when $m=1$.
	\begin{center}
	\begin{tabular}{cllll}
		$m$ & $\fv$ & $\fz$ & $J_\mu V$ & $G$\\
		\hline
		$1$ & $\C^n$ & $\R$ & $i\mu V$ & Heisenberg group $\heis_n$\\
		$2$ & quaternionic, $n$ even & ${\rm Span}\{i, j\}$ & $(\mu_1 i+\mu_2 j)V$ & the case of Figure \ref{fig:visible}\\
		$3$ & quaternionic, $n$ even & imaginary quaternions & $\mu V$ & quaternionic Heisenberg group\\
		$7$ & octonions, $n=4$ & imaginary octonions & $\mu V$ & octonionic Heisenberg group
	\end{tabular}
	\end{center}
	The rows $m=1, 3, 7$ are the Iwasawa $N$-groups of the complex, quaternionic and octonionic hyperbolic spaces \cite{cowlingHTypeGroupsLwasawa}.
\end{example}

\begin{remark}
	The $H$-type condition depends on the inner product, which we regard as part of the structure. Throughout, $\langle\cdot, \cdot\rangle$ is a fixed inner product on $\fg$ for which $(\fg, \langle\cdot, \cdot\rangle)$ is $H$-type, and the sub-Riemannian metric is its restriction to $\fv$ (Section \ref{HamFormalism}). Horizontal covectors $\nu\in \fv^*$ labelling geodesics are always taken of unit length, $|\nu|=1$, so that geodesics are parameterized by arclength.
\end{remark}

Since $\fg$ is nilpotent, $\exp: \fg\to G$ is a diffeomorphism, and we identify $G$ with $\fg$ through it. The group law is then the Dynkin product of the 2-step algebra $\fg$,
\begin{align}
	U\cdot V:=U+V+\tfrac{1}{2}[U, V], \qquad \exp(U)\exp(V)=\exp(U\cdot V). \label{Product}
\end{align}
Fix orthonormal bases $X_1, \ldots, X_{2n}$ of $\fv$ and $U_1, \ldots, U_m$ of $\fz$, and give $G$ the coordinates
\begin{align*}
	(x, u)=(x_1, \ldots, x_{2n}, u_1, \ldots, u_m)\longmapsto \exp\Big(\sum_{j=1}^{2n} x_jX_j+\sum_{k=1}^m u_kU_k\Big).
\end{align*}
In these coordinates \eqref{Product} reads
\begin{align}
	(x, u)(x', u')=(x+x', u+u'+\tfrac{1}{2}\omega( x, x')), \label{product}
\end{align}
where $\omega(x, x'):=(\omega_1(x, x'), \ldots, \omega_m(x, x'))\in \R^m$ and $\omega_k(x, x'):=\langle U_k, [\sum_i x_iX_i, \sum_j x'_jX_j]\rangle$, so that $\mu\cdot\omega(x, x')=\omega_\mu(x, x')$ for $\mu\in \fz^*$. For $\heis_1$, with $X_1=1$ and $X_2=i$ in $\fv=\C$, $\omega(x, x')=x_1x_2'-x_2x_1'$.

\subsection{Geodesics and the X-ray transform}\label{s:geodesicsXray}
 A unit-speed geodesic through the origin is determined by its initial momentum $(\nu, \lambda)\in S(\fv^*)\times\fz^*$. The horizontal momentum, initially $\nu$, evolves by the rotation $s\mapsto e^{sJ_\lambda}\nu$. Here $S(\fv^*)$ and $S(\fz^*)$ denote the unit spheres in $\fv^*$ and $\fz^*$. Let $\cG_\lambda\subset\cG$ be the geodesics of charge $\lambda$. Left translations are isometries and preserve the charge, so $G$ acts on each $\cG_\lambda$.

When $\lambda=0$ the geodesics are the horizontal lines and the X-ray transform along them is the horizontal line transform of the Heisenberg group \cite{strichartzLpHarmonicAnalysis1991, xiaoKplaneTransformHeisenberg2020}. We study $\lambda\neq 0$.

When $\lambda\neq0$, set $R:=|\lambda|^{-1}$, $\hat\lambda:=\lambda/|\lambda|$, and
\begin{align}\label{param}
	\gamma_{\nu, \lambda}(s):=\Big(e^{sJ_\lambda}J_\lambda^{-1}\nu,\ \tfrac{s}{2}R\hat\lambda\Big)\in G, \qquad \nu\in S(\fv^*),\ s\in \R.
\end{align}
This is the geodesic with momentum $(\nu, \lambda)$, translated so that its horizontal projection is the circle of radius $R$ centered at the origin, traversed at unit speed (Section \ref{spaceGeodesics}). The curve does not pass through the origin, which is its \textit{guiding center} (Appendix \ref{momentumFunctions}). Its geometry is summarized as follows.
\begin{center}
	\begin{tabular}{ll}
		radius of the horizontal circle & $R=|\lambda|^{-1}$\\
		period & $2\pi R=2\pi|\lambda|^{-1}$\\
		central advance per period & $\pi R^2\hat\lambda=\pi|\lambda|^{-2}\hat\lambda$\\
		first conjugate point & at arclength $2\pi R$ (Remark \ref{conjRmk})
	\end{tabular}
\end{center}
 The helical symmetry of the geodesic reads
\begin{align}\label{e:helicalSym}
	\gamma_{\nu, \lambda}(s+2\pi R)=(0, \pi R^2\hat\lambda)\,\gamma_{\nu, \lambda}(s).
\end{align}
Thus the geodesic returns to itself, up to translation, along the \textit{central lattice}
\begin{align*}
	\Gamma_\lambda:=\{(0, j\pi R^2\hat\lambda)\in G: j\in \Z\},
\end{align*}
a discrete subgroup of the center.

Every geodesic of nonzero charge is a left translate of some $\gamma_{\nu, \lambda}$. For $\lambda\neq0$ and $\nu\in S(\fv^*)$ let $\cG_{\nu, \lambda}:=\{(x, u)\gamma_{\nu, \lambda}: (x, u)\in G\}$ be the $G$-orbit of $\gamma_{\nu, \lambda}$. By Lemma \ref{l:fibers}, $\cG_{\nu, \lambda}$ depends on $\nu$ only through the circle $\{e^{sJ_\lambda}\nu: s\in \R\}\subset S(\fv^*)$, the orbits $\cG_{\nu, \lambda}$ partition $\cG_\lambda$, and the stabilizer of $\gamma_{\nu, \lambda}$ is $\Gamma_\lambda$, so the orbit--stabilizer theorem identifies
\begin{align}\label{e:orbitStabMain}
	\cG_{\nu, \lambda}\cong G/\Gamma_\lambda=:G_\lambda, \qquad (x, u)\gamma_{\nu, \lambda}\longleftrightarrow(x, u)\Gamma_\lambda.
\end{align}
The quotient $G_\lambda$ is again a Lie group, with the center $\fz$ replaced by the cylinder $\fz/\pi R^2\hat\lambda\Z$. We give it the Haar measure induced from $G$ and $\cG_{\nu, \lambda}$ the corresponding measure (Section \ref{spaceGeodesics}). The identification \eqref{e:orbitStabMain} depends on $\nu$, but the group $\Gamma_\lambda$ does not.

\begin{definition}
	For $\nu\in S(\fv^*)$ and $\lambda\in \fz^*\mz$, the X-ray transform restricted to $\cG_{\nu, \lambda}$, written in the coordinates \eqref{e:orbitStabMain}, is
	\begin{align*}
		I_{\nu, \lambda}f(x, u):=\int_\R f\left((x, u)\gamma_{\nu, \lambda}(s)\right)ds, \qquad (x, u)\in G_\lambda,\ f\in C_c^\infty(G),
	\end{align*}
	which is well defined on $G_\lambda$ by \eqref{e:helicalSym}. We also write $If(x, u, \nu, \lambda):=I_{\nu, \lambda}f(x, u)$. By Proposition \ref{factorization}, $I_{\nu, \lambda}$ extends to a bounded map $L^1(G)\to L^1(G_\lambda)$.
\end{definition}

\section{Main Results}\label{MainResults}

\begin{definition}
	Let $f\in L^1(G)$. Call a subset $Z\subset \fz^*\mz$ a \textit{zero set} (\textit{of charges}) of $If$ if for every $\lambda\in Z$ there exists $\nu\in S(\fv^*)$ such that $I_{\nu, \lambda}f=0$. Call $Z$ an \textit{injectivity set} of the X-ray transform $I$ (for $L^1$ functions) if the only $f\in L^1(G)$ with zero set $Z$ is $f\equiv 0$.
\end{definition}
A zero set requires the vanishing of $I_{\nu, \lambda}f$ for only \textit{one} direction $\nu$ per charge $\lambda\in Z$. The injectivity results (Theorem \ref{MainResult}, Corollary \ref{c:injSets}) and Theorem \ref{t:singleDirection} require only a zero set, while Theorem \ref{t:stability} and Corollary \ref{chargeFrequency} use the data of all directions $\nu\in S(\fv^*)$. Theorem \ref{PW2} holds for zero sets, with better constants for all-direction data.

The frequencies of $f$ seen by the geodesics of charge $\lambda$ are the characters of $\fz$ trivial on $\Gamma_\lambda$,
\begin{align*}
	\Gamma_\lambda^*:=\{\mu\in \fz^*: \mu\cdot\pi R^2\hat\lambda\in 2\pi\Z\}=\{\mu\in \fz^*: \mu\cdot\hat\lambda\in 2|\lambda|^2\Z\}=\{\mu\in \fz^*: \mu\cdot\lambda\in 2|\lambda|^3\Z\},
\end{align*}
which we call the \textit{visible frequencies} of charge $\lambda$. The set $\Gamma_\lambda^*$ is a countable family of parallel hyperplanes in $\fz^*$ orthogonal to $\lambda$ with spacing $2|\lambda|^2$.

\begin{theorem}[Characterization of Injectivity Sets] \label{MainResult}
	A subset $Z\subset \fz^*\mz$ is an injectivity set of $I$ if and only if $\bigcup_{\lambda\in Z}\Gamma_\lambda^*$ is dense in $\fz^*$.
\end{theorem}

\begin{corollary}\label{c:injSets}
	The X-ray transform is injective on $L^1(G)$. That is, if $f\in L^1(G)$ and $I_{\nu, \lambda}f=0$ for almost every $(\nu, \lambda)\in S(\fv^*)\times\fz^*$, then $f=0$. Moreover, $Z\subset \fz^*\mz$ is an injectivity set of $I$ if either
	\begin{itemize}
		\item[(i)] $Z$ contains a sequence $\{\lambda_i\}$ in $\fz^*\mz$ accumulating at zero, or
		\item[(ii)] $m>1$ and the set of unit vectors $\hat\mu\in S(\fz^*)$ whose orthogonal hyperplane $\hat\mu^\perp$ meets $Z$ is dense in $S(\fz^*)$. This holds, in particular, whenever the set of directions $\{\lambda/|\lambda|: \lambda\in Z\}$ is dense in some great circle of $S(\fz^*)$.
	\end{itemize}
\end{corollary}

 When $m=1$, $\Gamma_\lambda^*$ is a discrete set of spacing $2|\lambda|^2$, and Theorem \ref{MainResult} shows that (i) is necessary as well as sufficient. When $m>1$, $\Gamma_\lambda^*$ contains the hyperplane $\lambda^\perp$, and (ii) shows that one can also recover $f$ from charges of a single magnitude.

  In Theorem \ref{FourierSliceTheorem} below, $\cF(f)(\pi):=\int_G f(g)\pi(g)^*\,dg$ where $\pi=\pi_\mu$, $\mu\in \fz^*\mz$ is the irreducible unitary representation of $G$ on Bargmann--Fock space $\cF(\C^n)$ with central character $e^{i\mu\cdot u}$ (Section \ref{BargmannFock}), $\pi_{(\eta, 0)}$, $\eta\in \fv^*$, are the characters of $G$ trivial on the center, and $\cF(f)(\mu):=\cF(f)(\pi_\mu)$. $\cF_\lambda$ is defined similarly on $G_\lambda$ (Section \ref{GFT}) \cite{bahouriFrequencySpaceHeisenberg2019, fischerQuantizationNilpotentLie2016}. By \eqref{e:orbitStabMain} the data $I_{\nu, \lambda}f$ live on the group $G_\lambda$.
\begin{theorem}[$H$-type Fourier Slice Theorem]\label{FourierSliceTheorem}
	Let $f\in L^1(G)$, $\lambda \in \fz^*\setminus\{0\}$ and $\nu\in S(\fv^*)$.
	\begin{itemize}
		\item[1.] For $\eta\in \fv^*$, let $\eta_{\nu, \lambda}$ be the orthogonal projection of $\eta$ onto ${\rm Span}\{\nu, J_{\hat\lambda}\nu\}$. Then
		\begin{align}\label{scalerSlice}
			\cF_\lambda(I_{\nu, \lambda}f)(\pi_{(\eta, 0)})=2\pi|\lambda|^{-1}J_0(|\eta_{\nu, \lambda}|\,|\lambda|^{-1})\cF(f)(\pi_{(\eta, 0)})
		\end{align}
		\item[2.] For $\mu\in \fz^*\mz$, with $\lambda\cdot \mu=2k|\lambda|^3$, for some $k\in \mathbb{Z}$,
		\begin{align}
			\mathcal{F}_\lambda\left(I_{\nu, \lambda} f\right)(\pi_\mu)=2\pi|\lambda|^{-1}\mathcal{J}_{\nu, \lambda}( \mu)\circ\mathcal{F}(f)(\pi_\mu). \label{slice}
		\end{align}
	\end{itemize}
\end{theorem}

Equation \eqref{scalerSlice} is an equality of complex valued functions of $\eta\in \fv^*$, while Equation \eqref{slice} is an equality of bounded operators acting on Bargmann-Fock space $\cF(\C^n)$, defined in Section \ref{BF}.

The operator-valued multiplier $\cJ_{\nu, \lambda}(\mu)$ is the strong operator integral of the representation $\pi_\mu$ over the closed loop $\gamma_{\nu, \lambda}$ in $G/\Gamma_\lambda$,
\begin{align*}
	\cJ_{\nu, \lambda}(\mu)=\frac{|\lambda|}{2\pi}\int_0^{2\pi|\lambda|^{-1}}\pi_\mu(\gamma_{\nu, \lambda}(s))ds,
\end{align*}
the analog of the Bessel multiplier $J_0(\rho|\zeta|)=\frac{1}{2\pi}\int_0^{2\pi}e^{i\rho\,\zeta\cdot e^{i\theta}}\,d\theta$ of the mean value transform $M^\rho f(z)=\frac{1}{2\pi}\int_0^{2\pi} f(z+\rho e^{i\theta})\,d\theta$ on $\R^2.$%

Part 2 says that $I_{\nu, \lambda}f$ determines $\cJ_{\nu, \lambda}(\mu)\circ\cF(f)(\mu)$ on the visible frequencies $\mu\in \Gamma_\lambda^*$.

\begin{remark}[Twisted convolution form]\label{r:twistedSlice}
	For $f\in L^1(G)$ and $\mu\in \fz^*$, and for $g\in L^1(G_\lambda)$ and $\mu\in \Gamma_\lambda^*$, let
	\begin{align}\label{e:centralFT}
		f^\mu(x):=\int_\fz f(x, u)e^{-i\mu\cdot u}du, \qquad g^\mu(x):=\int_{\fz/\Gamma_\lambda}g(x, u)e^{-i\mu\cdot u}du
	\end{align}
	be the central Fourier transforms. Since $\pi_\mu(x, u)=e^{i\mu\cdot u}\pi_\mu(x, 0)$, the group Fourier transform is the Weyl transform $W_\mu$ (Section \ref{s:Weyl}) of the central Fourier transform,
	\begin{align}\label{e:Weyl}
		\cF(f)(\mu)=W_\mu(f^\mu), \qquad \cF_\lambda(g)(\mu)=W_\mu(g^\mu),
	\end{align}
	and $W_\mu$ turns the $\mu$-twisted convolution $f\times_\mu g(x):=\int_\fv f(y)g(x-y)e^{\frac{i}{2}\omega_\mu(x, y)}dy$ into a composition of bounded operators. In these terms Part 2 says that, for $\mu\cdot\lambda=2k|\lambda|^3$ and $x(s):=J_\lambda^{-1}e^{sJ_\lambda}\nu$ the horizontal circle of $\gamma_{\nu, \lambda}$,
	\begin{align}\label{e:twistedSlice}
		(I_{\nu, \lambda}f)^\mu(x)=\int_0^{2\pi|\lambda|^{-1}}e^{ik|\lambda|s}e^{\frac{i}{2}\omega_\mu(x, x(s))}f^\mu(x+x(s))\, ds=f^\mu\times_\mu\sigma^\mu_{\nu, \lambda}(x),
	\end{align}
	where $\sigma^\mu_{\nu, \lambda}$ is the arc measure of the circle $s\mapsto-x(s)$ weighted by the phase $e^{ik|\lambda|s}$, and $W_\mu(\sigma^\mu_{\nu, \lambda})=2\pi|\lambda|^{-1}\cJ_{\nu, \lambda}(\mu)$ \eqref{multiplier}. This is the classical projection-slice statement, convolution with the measure on the curve, with twisted convolution in place of convolution. At $\mu=0$ the phases disappear, $\times_0$ is ordinary convolution, and \eqref{e:twistedSlice} convolves $f^0=\int f\, du$ with the arc measure on the circle of radius $|\lambda|^{-1}$ in ${\rm Span}\{\nu, J_{\hat\lambda}\nu\}$, whose Fourier transform is the Bessel multiplier in Part 1.
    The multiplier in Part 2 is operator valued because twisted convolution is not diagonalized by the Euclidean Fourier transform on $\fv$.
\end{remark}

We study the normal operator of the operator-valued symbol $\cJ_{\nu, \lambda}(\mu)$ averaged over the horizontal directions. For $\lambda\in \fz^*\mz$ and $\mu\in \Gamma_\lambda^*\mz$, set
\begin{align*}
	\cN_{\nu, \lambda}(\mu):=\cJ_{\nu, \lambda}(\mu)^*\circ\cJ_{\nu, \lambda}(\mu), \qquad \cN_\lambda(\mu):=\int_{S(\fv^*)}\cN_{\nu, \lambda}(\mu)\, d\sigma(\nu),
\end{align*}
where $d\sigma$ is the normalized surface measure on $S(\fv^*)$. By Theorem \ref{FourierSliceTheorem}, $\cN_{\nu, \lambda}(\mu)$ plays the role of the symbol of the normal operator $I_{\nu, \lambda}^*I_{\nu, \lambda}$ at the frequency $\pi_\mu$, and $\cN_\lambda(\mu)$ that of the normal operator averaged over all geodesics of charge $\lambda$. Let $\Delta_{\rm{sr}}=\sum_{j=1}^{2n}X_j^2$ be the sub-Laplacian. The operator $d\pi_\mu(-\Delta_{\rm{sr}})$ on $\cF(\C^n)$ has pure point spectrum $\{|\mu|(2l+n): l\in \N\}$. Write $\cE_l$ for the eigenspace with eigenvalue $|\mu|(2l+n)$, which has dimension $\binom{l+n-1}{l}$, and ${\rm Pr}_l$ for the orthogonal projection onto it (Section \ref{BargmannFock}).

For $n>1$ and every nonzero visible central frequency $\mu$, Part 2 of Theorem \ref{t:stability} compares the $L^2$ norm of the data with a Sobolev norm of order $-\tfrac12$ in the sub-Laplacian.

\begin{theorem}[Averaged Normal Operator and Stability]\label{t:stability}
	Let $\lambda\in \fz^*\mz$ and $\mu\in \Gamma_\lambda^*\mz$, and write $\mu\cdot\lambda=2k|\lambda|^3$ and $r:=|\lambda|^{-2}|\mu|$.
	\begin{itemize}
		\item[1.] $\cN_\lambda(\mu)=\sum_{l\in \N}c_l\, {\rm Pr}_l$ is diagonalized by the spectral decomposition of the sub-Laplacian $\Delta_{\rm{sr}}$, with eigenvalues $c_l=c_l(k, r)\geq0$, given by \eqref{e:avgEval} with $|\mu|$ replaced by $r$.
		 Moreover, if $n>1$ then,
		\begin{align}\label{e:mainAsymp}
			c_l(k, r)=\frac{\Gamma(n)}{\sqrt{\pi}\,\Gamma(n-\tfrac12)}\big(r(2l+n)\big)^{-1/2}(1+o(1)), \qquad l\to\infty,
		\end{align}
		 and $\cN_\lambda(\mu)$ is injective (Proposition \ref{p:Inj}). %
		\item[2.]
		Let $f\in L^1(G)$ such that $\cF(f)(\mu)$ is Hilbert--Schmidt. If $n>1$ there are constants $0<C_1\leq C_2$, depending on $n$, $k$, $r$ and $|\lambda|$, such that, with $\Lambda_\mu:=d\pi_\mu(-\Delta_{\rm{sr}})$,
		\begin{align}\label{e:mainStab}
			C_1\big\|\Lambda_\mu^{-1/4}\cF(f)(\mu)\big\|_{HS}^2\leq\int_{S(\fv^*)}\big\|\cF_\lambda(I_{\nu, \lambda}f)(\pi_\mu)\big\|_{HS}^2\, d\sigma(\nu)\leq C_2\big\|\Lambda_\mu^{-1/4}\cF(f)(\mu)\big\|_{HS}^2.
		\end{align}
	\end{itemize}
\end{theorem}

Part 1 says that, on each central frequency, the averaged normal operator is a function of the sub-Laplacian, and \eqref{e:mainAsymp} identifies that function to leading order for $G\neq\heis_1$, at each fixed $\mu\in\Gamma_\lambda^*\mz$, as $l\to\infty$, with $K_n:=\Gamma(n)/(\sqrt{\pi}\,\Gamma(n-\tfrac12))$:
\begin{align*}
	\cN_\lambda(\mu)=K_n|\lambda|\,\Lambda_\mu^{-1/2}\,(1+o(1)).
\end{align*}

\begin{remark}
 At each fixed $\mu$ there are finitely many pairs $(k, r)$ such that $\mu\cdot\lambda=2k|\lambda|^3$ and $r=|\mu||\lambda|^{-2}$ as $\lambda$ ranges over a band or a cap of charges (Theorem \ref{PW2}). Hence, for $n>1$ and a band or a cap, the constants $C_1$ and $C_2$ may be chosen uniformly in $\lambda$.
Uniformity in $\mu$ requires Laguerre asymptotics uniform in $(k, r)$ and is left to a subsequent article. For $n=1$ no estimate of the form \eqref{e:mainStab} holds (Remark \ref{r:nOne}).
\end{remark}

The vanishing of $I_{\nu, \lambda}f$ in a single direction $\nu$ already determines much of $\cF(f)$ on $\Gamma_\lambda^*$. For $k\in \N^*$ let
\begin{align*}
	\cA_k:=\{a\in \N: L_a^{(k)}(k)=0\}
\end{align*}
be the set of degrees $a$ for which $L_a^{(k)}$ vanishes at the point $k$. One has $\cA_k=\emptyset$ for every odd $k$ \cite{flynnInjectivityHeisenbergXray2021}, while $L_2^{(2)}(2)=0$, so $2\in \cA_2$. For this reason, we let $\Gamma_{\lambda, \rm odd}^*:=\{\mu\in \fz^*: \mu\cdot\lambda\in 2|\lambda|^3(2\Z+1)\}$ be the odd part of $\Gamma_\lambda^*$.
\begin{theorem}[Single-Direction Data]\label{t:singleDirection}
	Let $f\in L^1(G)$, $\lambda\in \fz^*\mz$ and $\nu\in S(\fv^*)$, and suppose $I_{\nu, \lambda}f=0$. Let $\mu\in \Gamma_\lambda^*\mz$ with $\mu\cdot\lambda=2k|\lambda|^3$.
	\begin{itemize}
		\item[(a)] If $k=0$ or $k$ is odd, then $\cF(f)(\mu)=0$.
		\item[(b)] If $k\neq0$ is even, then ${\rm ran}\,\cF(f)(\mu)\subset\ker\cJ_{\nu, \lambda}(\mu)$, and $\ker\cJ_{\nu, \lambda}(\mu)$ is the closed span of an explicit orthonormal family in $\cF(\C^n)$ indexed by $\{\alpha\in \N^n: \alpha_1\in \cA_{|k|}\}$, given in the proof. In particular $\cF(f)(\mu)=0$ whenever $\cA_{|k|}=\emptyset$.
	\end{itemize}
\end{theorem}
Thus a single direction at charge $\lambda$ determines $\cF(f)$ on $\lambda^\perp\cup\Gamma^*_{\lambda, \rm odd}$. The proof of Theorem \ref{MainResult} uses nothing more, since $\bigcup_{\lambda\in Z}(\lambda^\perp\cup\Gamma^*_{\lambda, \rm odd})$ is dense in $\fz^*$ whenever $\bigcup_{\lambda\in Z}\Gamma^*_\lambda$ is dense. Since $2\in \cA_2$, the operator $\cJ_{\nu, \lambda}(\mu)$ has a nontrivial kernel for $\mu\cdot\lambda=4|\lambda|^3$, and $I_{\nu, \lambda}f=0$ does not force $\cF(f)(\mu)=0$ there. Averaging over all directions removes this obstruction when $n>1$ but not when $n=1$, where $\cN_\lambda(\mu)$ has a nontrivial kernel for $\mu\cdot\lambda=4|\lambda|^3$ (Corollary \ref{invCriteria}), so that on $\heis_1$ even the data of all directions at a single charge miss part of $\cF(f)(\mu)$ on the even hyperplanes.
\begin{corollary}[All-Direction Data]\label{chargeFrequency}
	Let $f\in L^1(G)$ and $\lambda\in \fz^*\mz$, and suppose $I_{\nu, \lambda}f=0$ for all $\nu\in S(\fv^*)$. If $n>1$ then $\cF(f)(\mu)=0$ for all $\mu\in \Gamma_\lambda^*\mz$. If $n=1$ then $\cF(f)(\mu)=0$ for all $\mu\in \Gamma_{\lambda, \rm odd}^*$.
\end{corollary}

Define the Paley--Wiener space
\begin{align*}
	L^1_{\Omega}(G):=\{f\in L^1(G): f^\mu=0 \text{ for all } \mu\in \Omega^c\}.
\end{align*}
Equivalently, $f\in L^1(G)$ is in $L^1_\Omega(G)$ if and only if $\cF f(\mu)=0$ for all $\mu\in \Omega^c$, by the injectivity of the Weyl transform on $L^1(\fv)$ (Section \ref{s:Weyl}).
For $\lambda_0\in \fz^*\mz$ with $|\lambda_0|=1$ and $0<h\leq 1$ define the spherical cap of height $h$
$$C_h(\lambda_0):=\{\lambda\in S(\fz^*): \lambda\cdot\lambda_0\geq 1-h \}.$$
For compact $Z\subset S(\fz^*)$ let $\delta:=\min_{\hat\mu\in S(\fz^*)}\big(\max_{\lambda\in Z}\hat\mu\cdot\lambda-\min_{\lambda\in Z}\hat\mu\cdot\lambda\big)$ be its minimal width, so that $\delta>0$ if and only if $Z$ lies in no affine hyperplane, and $\delta=h$ for the cap $C_h(\lambda_0)$.
\begin{theorem}[Frequency-Support Theorems]\label{PW2}
	Let $f\in L^1(G)$ and let $Z\subset \fz^*\mz$ be a zero set of $If$. Then $f\in L^1_\Omega(G)$ with $\Omega=\{\mu\in \fz^*: |\mu|\leq c\}$ in the following cases.
	\begin{itemize}
		\item[1.] $m=1$ and $Z=\{\lambda\in \fz^*: |\lambda|\in [1, 1+\eps]\}$ for some $\eps>0$, with $c=2k_0$, where $k_0$ is the least odd integer $\geq 2/(2\eps+\eps^2)$.
		\item[2.] $m>1$ and $Z\subset S(\fz^*)$ is compact, connected, and not contained in any affine hyperplane of $\fz^*$, with $c=4/\delta$, where $\delta$ is the minimal width of $Z$. In particular, for the spherical cap $Z=C_h(\lambda_0)$, $0<h\leq 1$, $c=4/h$.
	\end{itemize}
	If moreover $I_{\nu, \lambda}f=0$ for all $\nu\in S(\fv^*)$ and all $\lambda\in Z$, one may take $c=2/\delta$ in part 2, with $2/h$ for the cap, and, when $n>1$, $c=2k_1$ in part 1, where $k_1$ is the least integer $\geq 1/(2\eps+\eps^2)$.
\end{theorem}
\begin{remark}[Holes]\label{r:holes}
	The ball $|\mu|\leq c$ in Theorem \ref{PW2} is not the exact set of undetermined frequencies. Let $Z$ be as in part 2 and suppose $I_{\nu, \lambda}f=0$ for all $\nu\in S(\fv^*)$ and $\lambda\in Z$. By Corollary \ref{chargeFrequency} the data determine $\cF f$ at every nonzero frequency in $\bigcup_{\lambda\in Z}\Gamma_\lambda^*$. The complement of this set is
	\begin{align*}
		\fz^*\setminus\bigcup_{\lambda\in Z}\Gamma_\lambda^*=\coprod_{k\geq 0}\big(H_k\sqcup -H_k\big), \qquad H_k:=\{\mu\in \fz^*: 2k<\mu\cdot\lambda<2k+2 \text{ for all } \lambda\in Z\}.
	\end{align*}
	 Each $H_k=\{\min_{\lambda\in Z}\mu\cdot\lambda>2k\}\cap\{\max_{\lambda\in Z}\mu\cdot\lambda<2k+2\}$ is open and convex, and the proof of Theorem \ref{PW2} shows $H_k\subset\{|\mu|<2/\delta\}$, so only finitely many are nonempty. Thus the data determine $\cF f(\mu)$ off a finite union of open convex holes; see Figure~\ref{fig:visible}.
\end{remark}

\begin{figure}[t]
\centering
\begin{tikzpicture}[scale=0.6, >=stealth]
  \begin{scope}
    \clip (-7,-7) rectangle (7,7);
    \foreach \s in {1,-1}{
      \foreach \k in {0,1,2}{
        \begin{scope}[scale=\s, even odd rule]
          \clip (-13,-13) rectangle (13,13)
            (135:{2*\k+2}) arc (135:45:{2*\k+2}) -- (45:20) -- (135:20) -- cycle;
          \fill[blue!18] (0,{2*\k*sqrt(2)}) -- ({-sqrt(2)},{(2*\k+1)*sqrt(2)})
            -- (0,{(2*\k+2)*sqrt(2)}) -- ({sqrt(2)},{(2*\k+1)*sqrt(2)}) -- cycle;
        \end{scope}
      }
    }
    \foreach \r in {2,4,6}{ \draw[gray!70, dashed, thin] (0,0) circle (\r); }
    \foreach \t in {45,52.5,...,135}{
      \foreach \k in {-3,...,3}{
        \draw[gray!55, line width=0.25pt]
          ({2*\k*cos(\t)},{2*\k*sin(\t)}) +({\t+90}:14) -- +({\t-90}:14);
      }
    }
    \foreach \k in {-3,...,3}{
      \draw[black, line width=0.8pt]
        ({2*\k*cos(110)},{2*\k*sin(110)}) +(200:14) -- +(20:14);
    }
  \end{scope}
  \draw[very thick] (45:1) arc (45:135:1);
  \draw[->, thick] (0,0) -- (110:1);
  \node at (-0.9,1.3) {$\lambda$};
  \node at (1.35,1.0) {$Z$};
  \node at (0.55,1.55) {$H_0$};
  \node at (0,3.5) {$H_1$};
  \draw[->] (1.8,6.6) -- (0.1,5.85);
  \node at (2.3,6.7) {$H_2$};
  \node at (0,-1.4) {$-H_0$};
  \node at (0,-3.5) {$-H_1$};
  \draw[->] (1.8,-6.6) -- (0.1,-5.85);
  \node at (2.3,-6.7) {$-H_2$};
  \node[fill=white, inner sep=1pt] at (200:2.4) {\scriptsize $|\mu|=2$};
  \node[fill=white, inner sep=1pt] at (200:4.3) {\scriptsize $4$};
  \node[fill=white, inner sep=1pt] at (200:6.3) {\scriptsize $6$};
  \node[anchor=west] at (7.1,4.67) {$\Gamma^*_{\lambda}$};
\end{tikzpicture}
\caption{The visible frequencies for $m=2$ and $Z=\{\lambda\in S(\fz^*): 45^\circ\leq\arg\lambda\leq135^\circ\}$, so that $\delta=1-\cos45^\circ$ and $2/\delta\approx6.8$. The thin lines are $\mu\cdot\lambda=2k$ for $\lambda\in Z$ and $|k|\leq3$, and the $k$-th family envelops the dashed circle $|\mu|=2|k|$. The bold lines are $\Gamma_{\lambda}^*$ for one $\lambda\in Z$. The shaded regions are the holes $\pm H_k$, $0\leq k\leq2$, whose union is the complement of $\bigcup_{\lambda\in Z}\Gamma_\lambda^*$. $H_k=\emptyset$ for $k\geq3$, and every hole lies in $|\mu|<6<2/\delta$.}
\label{fig:visible}
\end{figure}
\begin{remark}
	By Proposition \ref{homogeneity}, if the conclusion of Theorem \ref{PW2} holds for $Z$ with constant $c$, then it holds for $tZ=\{t\lambda: \lambda\in Z\}$ with constant $ct^2$, for every $t>0$. In particular part 1 covers every band $\{|\lambda|\in [t, t(1+\eps)]\}$.
\end{remark}
As $\lambda$ ranges over a band or a cap, the hyperplanes $\Gamma_\lambda^*$ sweep out every frequency $|\mu|\geq c$.
For functions with compact central support the ball $|\mu|\leq c$ is recovered by analytic continuation.
\begin{corollary}\label{c:centralSupport}
	Let $f\in L^1(G)$ be supported in $\exp(\fv\times K)$ for some compact $K\subset \fz$. If $Z$ is a zero set of $If$, with $Z$ as in Theorem \ref{PW2}, then $f=0$. In particular, an arbitrarily small spherical cap of charges is an injectivity set for functions with compact central support.
\end{corollary}

In Euclidean limited-angle tomography the data determine a function only microlocally, at the covectors conormal to some measured line, the \textit{visible} singularities of Quinto \cite{quintoSingularitiesXrayTransform1993}. The visible set is a cone in $T^*\R^n$, and the invisible part is recovered, if at all, by analytic continuation. On an $H$-type group the natural frequency space is not $T^*G$ but the unitary dual $\widehat G$, parameterized up to a Plancherel-null set by the central frequency \cite{bahouriFrequencySpaceHeisenberg2019, fischerQuantizationNilpotentLie2016}, and Corollary \ref{chargeFrequency} and Theorem \ref{PW2} are visible-frequency theorems in this sense. The visible set is a neighborhood of infinity in $\widehat G$ rather than a cone, the determination there is stable at every central frequency when $n>1$ (Theorem \ref{t:stability}), and the invisible part, a bounded set of central frequencies, is recovered only by analytic continuation, for $f$ compactly supported in the central variable (Corollary \ref{c:centralSupport}).

\section{Sub-Riemannian Geometry of $H$-type groups}\label{sRGeometry}
\subsection{Hamiltonian Formalism}\label{HamFormalism}
We define the sub-Riemannian metric on $G$ by declaring $X_1, \ldots X_{2n}$ orthonormal, and the length of $U_1, \ldots U_m$ to be infinite. This amounts to restricting the inner product $\langle , \rangle$ discussed in Section \ref{H-type} from all of $\fg=\fv\oplus \fz$ to the first stratum $\fv$. Then any finite length smooth path in $G$ must be tangent to the nonintegrable distribution $\cD_q=dL_q(\fv)$. We call such a path \textit{horizontal}. The length of a horizontal path equals the length of its projection to the plane by the map
\begin{align*}
	dx: (x, u)\mapsto x\in \R^{2n}.
\end{align*}
A minimizing $H$-type geodesic is the shortest horizontal path joining two points of $G$. That any two points in $G$ are connected by a horizontal path is guaranteed by Chow's Theorem and the fact that $\cD$ satisfies the H\"{o}rmander condition (i.e. $\cD$ is bracket generating.)

The fiber-quadratic Hamiltonian $H: T^*G\to \R$ given in canonical coordinates by $$H=\tfrac{1}{2}\left(P_{X_{1}}^2+\cdots+P_{X_{2n}}^2\right)=\tfrac{1}{2}|\xi+\tfrac{1}{2}J_\mu x|^2$$
generates the geodesics. That is, any solution for Hamilton's equations for $H$ projects, via the canonical projection $T^*G\to G$, to a sub-Riemannian geodesic, and conversely, all geodesics arise in this way, since $G$ has step two and so every geodesic is normal \cite{montgomeryTourSubriemannianGeometries2006}. To obtain geodesics parameterized by arclength, we only take solutions for which $H=1/2$. (Thus, we define the unit cotangent bundle $U^*G$ as the set of all $(q, p)\in T^*G$ for which $H(q, p)=1/2$).

Geodesics are best understood by their projection to $\R^{2n}$: they are ``generalized helices."
Indeed
\begin{align*}
	\dot{p}_u=-\frac{\partial H}{\partial u}=0
\end{align*}
so that $\lambda=p_u$ is a constant of motion.
 If we choose a basis for $\fv$ as in \eqref{symplecticBasis} that puts $J_\lambda$ into block diagonal form then $\fv$ is a direct sum of two-planes each invariant under $J_\lambda$. The geodesics projected to any one of these two-planes are either all circles or lines.
 If we interpret $\lambda$ as the ``generalized charge" of a particle, then Hamilton's equations for $H$, viewed as a Hamiltonian on $T^*\R^{2n}$, are just the Lorentz equations for a particle of charge $\lambda$ in a constant magnetic field.

\subsection{Geodesic Flow}

The Hamiltonian flow written in the coordinate $(x, u, \nu, \lambda)$ is exactly the left trivialization $T^*G\cong G\times \fg^*$. That is
\begin{align*}
	\varphi_s: G\times \fg^*\to G\times \fg^*
\end{align*}
The unique solution to Hamilton's equations for a curve starting at the origin with initial momentum $(\nu, \lambda)\in \fv^*\oplus\fz^*$ is
\begin{align*}
	s\mapsto \left(x(s), u(s), e^{sJ_\lambda}\nu, \lambda\right)\in G\times \fg^*
\end{align*}
with $x(s)$ and $u(s)$ given by
\begin{align}\label{geodesics}
	x(s)=
	\begin{cases}
		\frac{e^{sJ_\lambda}-1}{J_\lambda}{\nu}, & \lambda\neq 0\\
		s{\nu}, & \lambda=0
	\end{cases}, \qquad \qquad
	u(s)=
	\begin{cases}
		\left(\frac{|\lambda|s-\sin{(|\lambda|s)}}{2|\lambda|^2}\right)\frac{\lambda}{|\lambda|}|{\nu}|^2, & \lambda\neq 0 \\
		0, & \lambda=0
	\end{cases}.
\end{align}

The unique constant speed geodesic starting at the origin with initial momentum $(\nu, \lambda)\in \fv^*\oplus\fz^*$ is the projection to the $(x, u)$ component
\begin{align}
	s\mapsto (x(s), u(s)) \in G.\label{geodesic}
\end{align}

Motivated by Appendix \ref{momentumFunctions}, we view the geodesic flow in \textit{guiding center coordinates}: for $\nu\in S(\fv^*)$ and $\lambda\neq0$, left translation by $(J_\lambda^{-1}\nu, 0)$ moves the guiding center of \eqref{geodesic} to the origin and yields the curve
\begin{align*}
	\gamma_{\nu, \lambda}(s)=\left(J_\lambda^{-1}\nu, 0\right)(x(s), u(s))
\end{align*}
of \eqref{param}.
Moreover since $H$ and $\lambda$ are integrals of motion, the flow restricts to the $(H, \lambda)$-level sets
\begin{align*}
	\varphi_{s}: G\times S(\fv^*)\times\{\lambda\}\to G\times S(\fv^*)\times\{\lambda\}.
\end{align*}
The map
\begin{align*}
	\Phi: G\times S(\fv^*)\times(\fz^*\mz)\to\cG\setminus\cG_0, \qquad (x, u, \nu, \lambda)\longmapsto(x, u)\gamma_{\nu, \lambda},
\end{align*}
is onto. It is not injective: reparameterizing by $s\mapsto s+s_0$ in \eqref{param} gives
\begin{align}\label{e:reparam}
	\gamma_{\nu, \lambda}(s+s_0)=\big(0, \tfrac{s_0}{2}R\hat\lambda\big)\,\gamma_{e^{s_0J_\lambda}\nu, \lambda}(s),
\end{align}
so that $\Phi(x, u, \nu, \lambda)=\Phi\big((x, u)(0, \tfrac{s_0}{2}R\hat\lambda), e^{s_0J_\lambda}\nu, \lambda\big)$ for every $s_0\in \R$, and these are the only coincidences (Lemma \ref{l:fibers}). The fibers of $\Phi$ are lines, and $\dim\cG=2\dim G-2$, as for the space of lines in a Euclidean space of the same dimension.
\subsection{The Space of Geodesics}\label{spaceGeodesics}
While not strictly necessary for the main results, we describe the space of $\lambda=0$ geodesics.

Let $\cG_{\nu, 0}$ be the subset of geodesics of the form $(x, u)\gamma_{\nu, 0}$, with $\gamma_{\nu, 0}(s)=(s\nu, 0)$, and for $\nu\in S(\fv^*)$ set $I_{\nu, 0}f(x, u):=\int_\R f((x, u)\gamma_{\nu, 0}(s))\,ds$, which is well defined on $G/\R\nu$. Since left translation by an element of $G$ does not change the direction $\nu$, or the charge $\lambda=0$, the group $G$ acts on $\cG_{\nu, 0}$. Thus we may parameterize that part of $\cG$ with $\lambda=0$ by
\begin{align*}
	G\times S(\fv^*)\to \cG_{0}; \quad (x, u, \nu)\mapsto (x, u)\gamma_{\nu, 0}
\end{align*}
for $(x, u)\in G$ and $\nu\in S(\fv^*)$.

Since the action of $G$ on $\cG_{\nu, 0}$ is transitive, we obtain, for any fixed $\nu\in S(\fv^*)$, a parameterization of $\cG_{\nu, 0}$ by $G$ unique up to the isotropy subgroup $\R \nu:=\{r\nu: r\in \R\}$ fixing $\gamma_{\nu, 0}$. Therefore
\begin{align}\label{degGeo}
	\cG_{\nu, 0}\cong &G/\R\nu\cong T_\nu S^{2n-1}\times\fz\\
	\nonumber (x, u)\gamma_{\nu, 0}\mapsto &(x, u)\R\nu \mapsto \left({{\rm Pr}}_\nu^\perp x, u-\tfrac12\langle x,\nu\rangle\omega(x,\nu)\right)
\end{align}
where ${{\rm Pr}}_\nu^\perp(x)=x-\langle x, \nu\rangle \nu$ is the orthogonal projection.
What's more, if $\cG_0$ is the subset of all geodesics with momentum $\lambda=0$, we have the identification
\begin{align*}
	\cG_0\to TS^{2n-1}\times \R^m.
\end{align*}
For any $\nu\in S(\fv^*)$, let $d\nu^\perp du$ be the measure on $T_\nu S(\fv^*)\times\fz=\{(\eta, u)\in \fv^*\times\fz: \nu\cdot\eta=0\}$. Equip $\cG_{\nu, 0}$ with the same measure through the identification \eqref{degGeo}.

For $\lambda\neq0$, the orbit--stabilizer identification \eqref{e:orbitStabMain} of $\cG_{\nu, \lambda}$ with $G_\lambda=G/\Gamma_\lambda$ rests on the following description of the fibers of $\Phi$.
\begin{lemma}[Fibers of the parameterization]\label{l:fibers}
	Let $\lambda, \lambda'\in \fz^*\mz$, $\nu, \nu'\in S(\fv^*)$ and $(x, u), (x', u')\in G$. Then $(x, u)\gamma_{\nu, \lambda}$ and $(x', u')\gamma_{\nu', \lambda'}$ are the same oriented unparameterized geodesic if and only if $\lambda=\lambda'$ and there is $s_0\in \R$ with $\nu'=e^{s_0J_\lambda}\nu$ and $(x', u')=(x, u)(0, \tfrac{s_0}{2}|\lambda|^{-1}\hat\lambda)$. In particular the stabilizer of $\gamma_{\nu, \lambda}$ in $G$ is $\Gamma_\lambda$.
\end{lemma}
\begin{proof}
	Sufficiency is \eqref{e:reparam}. Conversely, two parameterizations by arclength of the same oriented geodesic differ by a shift, so $(x', u')\gamma_{\nu', \lambda'}(s)=(x, u)\gamma_{\nu, \lambda}(s+s_0)$ for all $s$ and some $s_0$. The charge is the central momentum, a constant of motion invariant under left translation, so $\lambda=\lambda'$. By \eqref{e:reparam} we may replace $(x, u)$ by $(x, u)(0, \tfrac{s_0}{2}|\lambda|^{-1}\hat\lambda)$ and $\nu$ by $e^{s_0J_\lambda}\nu$, and assume $s_0=0$. Comparing horizontal components in $(x', u')\gamma_{\nu', \lambda}(s)=(x, u)\gamma_{\nu, \lambda}(s)$ gives $x'+e^{sJ_\lambda}J_\lambda^{-1}\nu'=x+e^{sJ_\lambda}J_\lambda^{-1}\nu$ for all $s$, hence $x=x'$ and $\nu=\nu'$; then the central components give $u=u'$. For the stabilizer, take $(x, u)=(0, 0)$ and $\nu'=\nu$: then $e^{s_0J_\lambda}\nu=\nu$ forces $s_0\in 2\pi|\lambda|^{-1}\Z$, and $(x', u')=(0, j\pi|\lambda|^{-2}\hat\lambda)\in \Gamma_\lambda$.
\end{proof}
Let $d(x, u)$ be the Haar measure on $G$ determined by the orthonormal frame $X_1, \ldots, X_{2n}$, $U_1, \ldots, U_m$. Also denote by $d(x, u)$ the Haar measure on $G/\Gamma_\lambda$ normalized so that
\begin{align}\label{Haar2}
	\int_Gf(x, u)d(x, u)=\int_{G/\Gamma_\lambda}\sum_{(0, u')\in \Gamma_\lambda}f(x, u+u')d(x, u), \quad f\in C_c(G).
\end{align}
Equip $\cG_{\nu, \lambda}$ with the measure $d\cG_{\nu, \lambda}(x, u)=d(x, u)$ through its identification with $G/\Gamma_\lambda$. The orbits $\cG_{\nu, \lambda}$ partition $\cG_\lambda$ and depend on $\nu$ only through its orbit under $\{e^{sJ_\lambda}: s\in \R\}$ (Lemma \ref{l:fibers}), so $\cG_\lambda=\bigsqcup_{[\nu]\in \CP(\lambda, \fv^*)}\cG_{\nu, \lambda}$, where $\CP(\lambda, \fv^*):=S(\fv^*)/\{e^{sJ_\lambda}: s\in \R\}$. We give $\CP(\lambda, \fv^*)$ the measure $d\nu$ induced from the normalized surface measure on $S(\fv^*)$, and $\cG_\lambda$ the measure $d\cG_\lambda:=d\nu\, d\cG_{\nu, \lambda}$. This is well defined, since changing the representative $\nu$ changes the identification \eqref{e:orbitStabMain} by a translation, which preserves Haar measure.

\section{Basic Properties of the X-ray Transform}\label{s:properties}
\subsection{Homogeneity}

For $\eps>0$ the dilation map $\delta_\eps(x, u)=(\eps x, \eps^2 u)$ is an isomorphism of the group $G$. Furthermore,
\begin{align*}
	\delta_\eps: \Gamma_{\eps\lambda}\ni (0, k\pi|\eps\lambda|^{-2}\hat{\lambda})\longmapsto (0, k\pi |\lambda|^{-2}\hat{\lambda})\in \Gamma_{\lambda},
\end{align*}
so $\delta_\eps: G_{\eps\lambda}\to G_{\lambda}$ is well-defined. Denote by $\delta_\eps^*$ the pullback operator on functions
\begin{align*}
	&\delta_\eps^*: L^1(G)\to L^1(G), & \delta_\eps^*f=f\circ \delta_\eps\\
	&\delta_\eps^*: L^1(G_{\lambda})\to L^1(G_{\eps\lambda}), & \delta_\eps^*g=g\circ\delta_\eps.
\end{align*}
\begin{proposition}[Homogeneity]\label{homogeneity} For $f\in L^1(G)$, $\nu\in S(\fv^*)$ and $\lambda\in \fz^*\setminus\{0\}$
	\begin{align*}
		I_{\nu, \eps\lambda}(\delta_\eps^* f)(x, u)=(1/\eps)\,\delta_\eps^*\left(I_{\nu, \lambda}f\right)(x, u), \quad \forall \eps>0.
	\end{align*}
\end{proposition}
\begin{proof}
	The main observation is that
\begin{align*}
	\delta_\eps \gamma_{\nu, \eps\lambda}(s)=\gamma_{\nu, \lambda}(\eps s)
\end{align*}
Thus,
\begin{align*}
	I_{\nu, \eps\lambda}(\delta_\eps^* f)(x, u)
	&=\int_\R \delta_\eps^*f((x, u)\gamma_{\nu, \eps\lambda}(s))ds\\
	&=\int_\R f(\delta_\eps(x, u)\delta_\eps\gamma_{\nu, \eps\lambda}(s))ds\\
	&=\int_\R f(\delta_\eps (x, u)\gamma_{\nu, \lambda}(\eps s))ds\\
	&=(1/\eps)\int_\R f(\delta_\eps(x, u)\gamma_{\nu ,\lambda}(s))ds\\
	&=(1/\eps)(I_{\nu, \lambda} f)(\delta_\eps(x, u))
\end{align*}
as desired.
\end{proof}
\subsection{Boundedness and Factorization}
\begin{proposition}
	For any $\nu\in S(\fv^*)$, the X-ray transform $I_{\nu, 0}: L^1(G)\to L^1(G/\R\nu)$ is well-defined and bounded.
\end{proposition}
\begin{proof}
	For $f\in C_c(G)$,
	\begin{align*}
		\|I_{\nu, 0}f\|_{L^1(G/\R\nu)}
		&\leq\int_{\R^m} \int_{x\perp \nu}|I_{\nu, 0}f(x, u)|d\nu^\perp(x) du\\
		&\leq\int_{\R^m} \int_{x\perp \nu}\int_\R |f((x, u)(s\nu, 0))|dsd\nu^\perp(x)du\\
		&=\int_{\R^m}\int_{x\perp \nu} \int_\R| f(x+s\nu, u+s [x, \nu]/2)|dsd\nu^\perp(x)du\\
		&=\int_{\R^m}\int_{x\perp \nu} \int_\R| f(x+s\nu, u)|dsd\nu^\perp(x)du\\
		&=\int_{\fz}\int_{\fv}|f(x, u)|dxdu=\|f\|_{L^1(G)}.
	\end{align*}
Therefore, $I_{\nu, 0}$ extends to a bounded map between $L^1$ spaces.
\end{proof}
Here we decompose the X-ray transform into a periodization map and a special case of the Pompeiu Transform. (See \cite{agranovskyInjectivityPompeiuTransform1994}).
\begin{proposition}\label{factorization}
	For any $\lambda\in \mathfrak{z}^*\setminus \{0\}$ and $\nu\in S(\fv^*)$, the X-ray transform $I_{\nu, \lambda}: L^1(G)\to L^1(G_\lambda)$ is well-defined, bounded and factors in the following way:
	\[
	\begin{tikzcd}
		L^1(G) \arrow[d, "P_{\lambda}" left=1] \arrow[r, "{I}_{\nu, \lambda}"] &L^1(G_\lambda)\\
		L^1(G_{\lambda}) \arrow[ur, "\widetilde{I}_{\nu, \lambda}" below=10, pos=0.75]
	\end{tikzcd}
	\]
	where the maps which we call \textit{central periodization} and the \textit{holonomy transform} are given by
	\begin{align*}
		P_{\lambda}f(x, u)=\sum_{k\in\mathbb{Z}}f(x, u+\pi k |\lambda|^{-2}\hat{\lambda}),
		\quad \quad
		\widetilde{I}_{\nu, \lambda}g(x, u)=\int_0^{2\pi |\lambda|^{-1}}g\left( (x, u)\gamma_{\nu, \lambda}(s)\right)ds
	\end{align*}
	for $f\in L^1(G)$ and $g\in L^1(G_\lambda)$.
\end{proposition}
\begin{proof}
	By the dominated convergence theorem we have
	\begin{align*}
		\int_{G_{\lambda}}P_{\lambda}f(x, u)d(x, u)=\int_G f(x, u)d(x, u), \quad f\in C_c^\infty(G).
	\end{align*}
	In particular $$\|P_{\lambda}f\|_{L^1(G_{\lambda})}\leq \|f\|_{L^1(G)}.$$

	Now $\widetilde{I}_{\nu, \lambda}f$ is just a group convolution operator on $G_{\lambda}$ by a compactly supported convolution kernel. Hence it is bounded. In particular,
	\begin{align*}
		|\widetilde{I}_{\nu, \lambda}f(x, u)|\leq \int_0^{2\pi |\lambda|^{-1}}|f\left((x, u)\gamma_{\nu, \lambda}(s)\right)|ds
	\end{align*}
	and
	\begin{align*}
		\int_{G_\lambda}|\widetilde{I}_{\nu, \lambda} f(x, u)|d(x, u)
		&=\int_0^{2\pi |\lambda|^{-1}} \int_{G_\lambda}|f\left(x, u\right)|d(x, u)ds.
	\end{align*}
	Therefore
	\begin{align*}
		\|\widetilde{I}_{\nu, \lambda} f\|_{L^1(G_\lambda)}\leq 2\pi |\lambda|^{-1}\|f\|_{L^1(G_{\lambda})}.
	\end{align*}
	The decomposition follows from the helical symmetry \eqref{e:helicalSym} of the geodesics. For $f\in C_c^\infty(G)$
	\begin{align*}
		\widetilde{I}_{\nu, \lambda}\circ P_{\lambda}f(x, u)
		&=\int_0^{2\pi |\lambda|^{-1}}\sum_{k\in\mathbb{Z}}f\left((x, u+\pi k|\lambda|^{-2}\hat{\lambda})\gamma_{\nu, \lambda}(s)\right)ds\\
		&=\int_0^{2\pi |\lambda|^{-1}}\sum_{k\in\mathbb{Z}}f\left((x, u)\gamma_{\nu, \lambda}(s+2\pi k|\lambda|^{-1})\right)ds	\\
		&=\sum_{k\in\mathbb{Z}}\int_{2\pi k |\lambda|^{-1}}^{2\pi (k+1) |\lambda|^{-1}}f\left((x, u)\gamma_{\nu, \lambda}(s)\right)ds	\\
		&=\int_\R f\left((x, u)\gamma_{\nu, \lambda}(s)\right)ds=I_{\nu, \lambda}f(x, u)
	\end{align*}
	where the second to last equality follows from the dominated convergence theorem.

	By the previous two statements, we have
	\begin{align*}
		\| I_{\nu, \lambda} f\|_{L^1(G_\lambda)}\leq 2\pi |\lambda|^{-1}\|f\|_{L^1(G)}.
	\end{align*}
	In particular $I_{\nu, \lambda}$ extends to a bounded map from $L^1(G)$ to $L^1(G_\lambda)$.
\end{proof}

\begin{proposition}
		Let $\chi\in M(\fz^*)$ be any positive compactly supported Radon measure for which
		\begin{align}\label{centralWeight}
			\int_{\fz^*}|\lambda|^{-1}d\chi(\lambda)<\infty,
		\end{align}
	and let $\chi d\cG$ be the measure on the space of geodesics defined by
	\begin{align*}
		\int_{\cG}F \chi d\cG=\int_{\lambda\in \fz^*}\int_{\cG_\lambda}F d\cG_\lambda d\chi(\lambda), \quad F\in C_c^\infty(\cG).
	\end{align*}
		The full X-ray transform $I: L^1(G)\to L^1(\cG, \chi d\cG)$ is well-defined and bounded.
\end{proposition}
\begin{proof}
	\begin{align*}
		\|If\|_{L^1(\cG, \chi d\cG)}
		&=\int_{\cG}|If| \chi d\cG\\
		&=\int_{\lambda\in \fz^*}\int_{\cG_\lambda}|If|\, d\cG_\lambda\, d\chi(\lambda)\\
		&=\int_{\lambda\in \fz^*}\int_{\CP(\lambda, \fv^*)}\int_{G_{\lambda}}|I_{\nu, \lambda}f(x, u)|\, d(x, u)\, d\nu\, d\chi(\lambda)\\
		&\leq2\pi \int_{\fz^*}|\lambda|^{-1}d\chi(\lambda)\,\|f\|_{L^1(G)}
	\end{align*}
	which proves boundedness in view of the assumption \eqref{centralWeight}.
\end{proof}

\section{Representation Theory of $G$ and $G_\lambda$}\label{s:representations}
The Fourier transform on $G$ is defined by means of irreducible unitary representations.
There are many choices of conventions for representations of the Heisenberg group. We choose to use the Bargmann-Fock representation because certain intertwining operators are more explicit than for the Schr\"{o}dinger representation. %

\subsection{Bargmann-Fock Representation on the Heisenberg group}\label{BF}
Bargmann-Fock space is the following Hilbert space (equipped with the evident inner product)
\begin{align*}
	\cF(\C^n):=\{F: \C^n\to \C, \text{ holomorphic: }\frac{1}{\pi^n}\int_{\C^n}|F(\zeta)|^2e^{-|\zeta|^2}<\infty\}.
\end{align*}
An orthonormal basis is given by
\begin{align*}
	\omega_\alpha=\zeta^\alpha/\sqrt{\alpha!}, \quad \zeta\in \C^n, \; \alpha\in \N^n
\end{align*}
For each $h\in \R^*$ the Bargmann-Fock representation is the map
\begin{align*}
	{\pi}_h: \heis_n \to \cU(\cF(\C^n))
\end{align*}
given by
\begin{align*}
	{\pi}_{h}(z, t)F(\zeta):=e^{i ht-\sqrt{h/2}\zeta\cdot\overline{z}-\tfrac{ h}{4}|z|^2}F(\zeta+\sqrt{h/2}z), \quad F\in \cF(\C^n), \quad h>0,
\end{align*}
and ${\pi}_h(z, t)={\pi}_{|h|}(\overline{z}, -t)$ for $h<0$. Despite the dependence on a parameter $h\in \R^*$, we simply refer to this as \textit{the} Bargmann-Fock representation. Here $\zeta\cdot \overline{z}=\zeta_1\overline{z}_1+\cdots + \zeta_n\overline{z}_n$. The Bargmann-Fock representation is a strongly continuous unitary irreducible representation of the Heisenberg group on $\cF(\C^n)$.

\begin{remark}
	There are many other normalizations out there, but we have chosen the one for which
	\begin{itemize}
		\item ${\pi}_h(z, t)=e^{iht}{\pi}_h(z, 0)$
		\item ${\pi}_*(X_j-iY_j)=\sqrt{2h}\,\partial_{\zeta_j}$ and ${\pi}_*(X_j+iY_j)=-\sqrt{2h}\,\zeta_j$ for $h>0$, so that $\pi_*(-\Delta_{\rm{sr}})=h\sum_{j=1}^n(2\zeta_j\partial_{\zeta_j}+1)$ for the sub-Laplacian $\Delta_{\rm{sr}}=\sum_j X_j^2+Y_j^2$, with eigenvalue $h(2l+n)$ on the homogeneous polynomials of degree $l$.
	\end{itemize}
\end{remark}

\subsection{Bargmann-Fock Representation on $H$-type groups}\label{BargmannFock}
To construct representations on a general $H$-type group, we follow \cite{mullerSharpLpEstimates2008}. Given $\mu\in \mathfrak{z}^*\setminus\{0\}$ consider the operator $J_\mu:\mathfrak{v}\to\mathfrak{v}$. If $|\mu|=1$, then $J_\mu^2=-I$. In particular, $J_\mu$ has eigenvalues $\pm i$, and since it is orthogonal, there exists an orthonormal basis
\begin{align}
	P_1^\mu, \ldots, P^\mu_n, Q_1^\mu, \ldots, Q_n^\mu\label{symplecticBasis}
\end{align}
of $\mathfrak{v}$ which is symplectic with respect to the form $\omega_\mu$.

This means that, for every $\mu\in \mathfrak{z}^*\setminus \{0\}$, there exists an orthogonal transformation $R_\mu=R_{\mu/|\mu|}: \R^{2n}\to \mathfrak{v}$ such that
\begin{align}\label{normalForm}
	J_\mu=|\mu|R_\mu J R_\mu^t
\end{align}
where $J=\begin{bmatrix}
	0& -I \\
	I& 0
\end{bmatrix}$
with respect to the basis \eqref{symplecticBasis}.

\begin{remark}
    In general there is no continuous global assignment $\mu\mapsto R_\mu$ for all $\mu\in \fz^*\mz$. Indeed, the set of $J_\mu$ for $\mu\in \fz^*\mz$ forms an $S^{m-1}$ sphere over which the set of possible $R_\mu$ for $\mu\in \fz^*\mz$, forms a principal $U(n)$ bundle. One may still choose $R_\mu$ smoothly locally in $\mu$.
\end{remark}
\begin{lemma}
	The mapping $\alpha_\mu: G\to \mathbb{H}_n$ defined by
	\begin{align*}
		\alpha_\mu(x, u)=(R_\mu^t x, \tfrac{\mu\cdot u}{|\mu|}), \quad (x, u)\in \mathbb{R}^{2n}\times\mathbb{R}^m
	\end{align*}
	is a homomorphism.
\end{lemma}

Then define the irreducible representations of $G$ on $\cF(\C^n)$ by
\begin{align}\label{e:repH}
	\pi_\mu:={\pi}_{|\mu|}\circ\alpha_\mu.
\end{align}
\begin{remark}
	When $\mu\in \fz^*\mz$ we assume it is clear that $\pi_\mu$ is a representation of an $H$-type group $G$ with center $\fz$; when $h=|\mu|$ is a scalar, $\pi_h$ is a representation of the Heisenberg group $\heis_n$.
\end{remark}
Since $R_\mu$ is orthogonal, $d\alpha_\mu$ carries the sub-Laplacian $\Delta_{\rm{sr}}=\sum_{j=1}^{2n}X_j^2$ of $G$ to that of $\heis_n$, so by the remark in Section \ref{BF}, $d\pi_\mu(-\Delta_{\rm{sr}})=|\mu|\sum_{j=1}^n(2\zeta_j\partial_{\zeta_j}+1)$, with eigenvalue $|\mu|(2l+n)$ on the homogeneous polynomials of degree $l$.
\subsection{Intertwiners}
Let $S\in Sp(2n, \R)$ and $h>0$. Then $(z, t)\mapsto\pi_h(Sz, t)$ is an irreducible unitary representation of $\heis_n$ with the same central character as $\pi_h$, so by the Stone--von Neumann theorem there is a unitary operator $\tau(S)$ on $\cF(\C^n)$, unique up to a phase, such that
\begin{align}\label{e:intertwine}
    \pi_h(Sz, t)=\tau(S)\circ\pi_h(z, t)\circ\tau(S)^*, \quad (z, t)\in \heis_n.
\end{align}
The projective representation $S\mapsto\tau(S)$ is the metaplectic representation \cite[Ch.~4]{follandHarmonicAnalysisPhase1989}. We use $\tau(S)$ only through the conjugation, which does not depend on the phase. On the unitary group we fix a choice:
\begin{definition}\label{d:tau}
	For $U\in U(n)=Sp(2n, \R)\cap O(2n, \R)$, identified with the unitary group of $\C^n\cong(\R^{2n}, J)$, and $F\in \cF(\C^n)$, define
	\begin{align*}
		{\tau}(U)F(\zeta):=F(U^*\zeta), \quad \zeta \in \C^n.
	\end{align*}
\end{definition}
This is a unitary representation of $U(n)$ which preserves the homogeneous polynomials of each degree and satisfies \eqref{e:intertwine}:
\begin{lemma}\label{l:rotateAxis}
    Let $u\in \R^{2n}$ with $|u|=1$. For any $U\in U(n)$ the operator $\tau(U)$ is unitary on $\cF(\C^n)$, preserves $\cE_l:={\rm Span}_\C\{\omega_\alpha: |\alpha|=l\}$ for every $l\in \N$, and satisfies
    \begin{align*}
        \pi_h(Uz, t)=\tau(U)\circ\pi_h(z, t)\circ\tau(U)^*, \quad (z, t)\in \heis_n, \; h>0.
    \end{align*}
\end{lemma}

\begin{proof}
    The group $U(n)$ acts transitively on the unit sphere of $\C^n$. The substitution $\zeta\mapsto U^*\zeta$ preserves the Gaussian measure $\pi^{-n}e^{-|\zeta|^2}$, so $\tau(U)$ is unitary, and it preserves the homogeneous polynomials of each degree. Since $(U^*\zeta)\cdot\overline{z}=\zeta\cdot\overline{Uz}$ and $|Uz|=|z|$,
    \begin{align*}
        \tau(U)\circ{\pi}_{h}(z, t)\circ\tau(U)^*F(\zeta)
        &=e^{i ht-\sqrt{h/2}\,(U^*\zeta)\cdot\overline{z}-\tfrac{ h}{4}|z|^2}F(\zeta+\sqrt{h/2}\, Uz)
        ={\pi}_h(Uz, t)F(\zeta).\qedhere
    \end{align*}
\end{proof}

\begin{lemma}\label{l:tauGauss}
	Let $h>0$, $S\in Sp(2n, \R)$, and let $\tau(S)$ be a unitary operator satisfying \eqref{e:intertwine}. Write $Sz=Mz+N\bar z$ for $z\in \C^n\cong(\R^{2n}, J)$. Then $M$ is invertible, $\Theta:=\bar NM^{-1}$ is symmetric, and for every $q\in \C^n$
	\begin{align}\label{e:tauGauss}
		\tau(S)e^{q\cdot\zeta}=C_q\, e^{\frac12\zeta^t\Theta\zeta+((M^t)^{-1}q)\cdot\zeta}, \qquad C_q\in \C\mz.
	\end{align}
\end{lemma}
\begin{proof}
	Since $\omega(u, v)=v^tJu={\rm Im}(\bar u\cdot v)$, the identity $\omega(Su, Sv)=\omega(u, v)$, applied to $(u, v)$, $(iu, v)$ and $(u, iv)$, is equivalent to
	\begin{align*}
		M^*M-N^t\bar N=I, \qquad M^t\bar N=\bar N^tM.
	\end{align*}
	Hence $M^*M\geq I$, so $M$ is invertible; $\Theta=(M^t)^{-1}(M^t\bar N)M^{-1}$ is symmetric; $\Theta M=\bar N$; and $M^t(\bar M-\Theta N)=M^t\bar M-\bar N^tN=\overline{M^*M-N^t\bar N}=I$, that is,
	\begin{align}\label{e:MTheta}
		\bar M-\Theta N=(M^t)^{-1}.
	\end{align}
	Let $s:=\sqrt{h/2}$, $g_S:=\tau(S)\omega_0$, and $u\cdot\partial_\zeta:=\sum_j u_j\partial_{\zeta_j}$ for $u\in \C^n$. For $v\in \C^n$ the map $\R\ni\eps\mapsto\pi_h(\eps v, 0)\omega_0=e^{-s^2\eps^2|v|^2/2-s\eps\,\bar v\cdot\zeta}$ is differentiable into $\cF(\C^n)$, hence so is $\eps\mapsto\pi_h(\eps Sv, 0)g_S=\tau(S)\pi_h(\eps v, 0)\omega_0$. Point evaluations are continuous on $\cF(\C^n)$, so differentiating at $\eps=0$,
	\begin{align*}
		-\tau(S)(\bar v\cdot\zeta)=(Sv)\cdot\partial_\zeta g_S-(\overline{Sv}\cdot\zeta)\, g_S.
	\end{align*}
	Adding the same identity for $iv$, divided by $i$, and using $S(iv)=iMv-iN\bar v$,
	\begin{align*}
		(Mv)\cdot\partial_\zeta g_S=(\bar Nv\cdot\zeta)\, g_S, \quad v\in \C^n, \qquad\text{that is,}\qquad \partial_\zeta g_S=\Theta\zeta\, g_S.
	\end{align*}
	So $\partial_\zeta\big(e^{-\frac12\zeta^t\Theta\zeta}g_S\big)=0$ and $g_S=C_0\,e^{\frac12\zeta^t\Theta\zeta}$, with $C_0\neq0$ since $\tau(S)$ is unitary.

	For $q\in \C^n$, $e^{q\cdot\zeta}=e^{|q|^2/2}\pi_h(-\bar q/s, 0)\omega_0$. By \eqref{e:intertwine}, with $w:=S\bar q=M\bar q+Nq$,
	\begin{align*}
		\tau(S)e^{q\cdot\zeta}=e^{|q|^2/2}\pi_h(-w/s, 0)g_S=C_0\, e^{\frac12(|q|^2-|w|^2)+\bar w\cdot\zeta+\frac12(\zeta-w)^t\Theta(\zeta-w)}.
	\end{align*}
	The coefficient of $\zeta$ in the exponent is $\bar w-\Theta w=(\bar M-\Theta N)q+(\bar N-\Theta M)\bar q=(M^t)^{-1}q$ by $\Theta M=\bar N$ and \eqref{e:MTheta}. This is \eqref{e:tauGauss} with $C_q=C_0\,e^{\frac12(|q|^2-|w|^2+w^t\Theta w)}$.
\end{proof}

\subsection{Finite Dimensional Representations on $H$-type groups}
For any $\eta\in \fv^*$ define
\begin{align}\label{finiteRep}
	\pi_{\eta, 0}(x, u)=e^{i\eta(X) }, \quad (x, u)=\exp(X+U)\in G, \quad \text{with } X\in \fv, \; \text{and }U\in \fz.
\end{align}
Then $\pi_{\eta, 0}: G\to U(1)$ is a one dimensional unitary representation of $G$ on $\C$.

\subsection{The Unitary Dual}
The \textit{unitary dual} $\widehat{G}$ of a group $G$ is the set of its irreducible unitary representations up to unitary equivalence. That is, the set of pairs $(\cH_\pi, \pi)$ where $\cH_\pi$ is a Hilbert space, and $\pi: G\to \cU(\cH_\pi)$ is a strongly continuous group homomorphism from $G$ to the set of unitary operators on $\cH_\pi$. When $G$ is a nilpotent Lie group, Kirillov theory establishes a bijection between $\widehat{G}$ and orbits of $\fg^*$ by the coadjoint action.

When $G$ is an $H$-type group, the Kirillov correspondence associates to each $\eta\in \fv^*$ the representation $\pi_{(\eta, 0)}$ in \eqref{finiteRep} acting on $\cH_{\pi_{\eta, 0}}=\mathbb{C}$, and to each $\mu\in \fz^*\mz$ the infinite dimensional representation $\pi_{\mu}$ acting on $\cH_{\pi_\mu}=\cF(\C^n)$. Since $\fv^*\cup \fz^*\mz$ exhausts all representatives of the coadjoint orbits of $\fg^*$, all elements of $\widehat{G}$ are obtained in this way \cite{kirillovUnitaryRepresentationsNilpotent1962}. Therefore $\widehat{G}$ is parameterized by $\fv^*\cup \fz^*\setminus\{0\}$:
\begin{align*}
	\widehat{G}=\{\text{class of } \pi_{\eta, 0}: \eta\in \fv^*\}\cup \{\text{class of } \pi_\mu: \mu \in \fz^*\setminus\{0\}\}.
\end{align*}
\subsection{Representations on $G_\lambda$}

\begin{proposition}[Compatibility of momenta]\label{propCompatibility}
	For $\mu\in \fz^*\mz$ and the associated representation $\pi_\mu$ of $G$, there exists a representation of $G_\lambda$ also denoted by $\pi_\mu$ making the following diagram
	\[
	\begin{tikzcd}
		G \arrow[d] \arrow[r, "\pi_\mu"] & \cU(\cF(\C^n))\\
		G_{\lambda} \arrow[ur, "\pi_\mu" below=7, pos=0.5]
	\end{tikzcd}
	\]
	commute if and only if
	\begin{align*}
		\mu\cdot \lambda=2|\lambda|^3k, \quad \text{ for some } k\in \Z.
	\end{align*}
	That is, $\pi_\mu$ descends from an irreducible unitary representation of $G$ to one of $G_\lambda$ if and only if $\mu$ is an element of the dual of the subgroup $\Gamma_\lambda$ excluding zero.
\end{proposition}
\begin{proof}
	Fix $\lambda\in \fz^*\mz$. Then $\pi_\mu$ is a representation of $G_\lambda$ if and only if $\Gamma_\lambda\subset \ker{\pi_\mu}$. Since $\pi_\mu(x, u)=e^{i\mu\cdot u}\pi_\mu(x, 0)$, we have
	\begin{align*}
		\ker\pi_\mu
		&=\{(x, u)\in G: \pi_\mu(x, u)={{\rm Id}}_{\cF(\C^n)}\}\\
		&=\{(0, u)\in G: u\cdot \mu\in 2\pi\Z \}.
	\end{align*}
	Then $\Gamma_\lambda=\{(0, l\pi|\lambda|^{-2}\, \hat{\lambda}): l\in \Z\}\subset \ker{\pi_\mu}$ if and only if
	\begin{align*}
		\pi |\lambda|^{-2}\,\hat{\lambda}\cdot \mu =2\pi k, \quad k\in \Z.
	\end{align*}
	This is precisely the condition that $\mu$ lies in the annihilator $\Gamma_\lambda^*:=\{\mu\in \fz^*: \mu\cdot\omega\in 2\pi \Z, \text{for all }\omega\in \Gamma_\lambda\}$ of the central lattice. We have

	\begin{align}
		\mu\cdot \lambda=2|\lambda|^3k, \quad k\in \mathbb{Z} \label{compadability}
	\end{align}
as desired.
\end{proof}
\begin{remark}
	Equation \eqref{compadability} implies
		\begin{align*}
		\mu=\mu_\perp+2k|\lambda|\lambda
	\end{align*}
where $\mu_\perp\in \fz^*$ such that $\mu_\perp\cdot\lambda=0$.
\end{remark}
\begin{remark}
	Note $k=0$ in Equation \eqref{compadability} if and only if $\lambda, \mu\in \fz^*\mz$ are perpendicular. This is only possible when $m=\dim \fz >1$.
\end{remark}

By Proposition \ref{propCompatibility}, we may view the unitary dual $\widehat{G}_\lambda$ as a subset of $\widehat{G}$.
\begin{align*}
	\widehat{G}_\lambda=\{\text{class of }\pi_{\eta, 0}: \eta\in \fv^*\}\cup \{\text{class of } \pi_\mu: \mu\in \Gamma_\lambda^*\mz \}\subset \widehat{G}
\end{align*}
where
\begin{align*}
	\Gamma_\lambda^*:=\{\mu\in \fz^*: \mu\cdot \lambda=2k|\lambda|^3, \, k\in \Z\}\subset \fz^*
\end{align*}
is the dual to the discrete subgroup $\Gamma_\lambda$.

\begin{remark}\label{conjRmk}
	For any $\lambda\in \fz^*\mz$ the geodesic $\gamma_{\nu, \lambda}(s)$ starting at $s=0$ reaches its first conjugate point at $s=2\pi R$, $R=|\lambda|^{-1}$ (see \cite{montgomeryTourSubriemannianGeometries2006}), and the two points in conjugacy differ by a group translation of $(0, \pi R^2 \hat{\lambda})\in G$ with $\hat{\lambda}=|\lambda|^{-1}\lambda$. Indeed
	\begin{align*}
		\gamma_{\nu, \lambda}(s+2\pi R)=(0, \pi R^2\hat{\lambda})\gamma_{\nu, \lambda}(s)
	\end{align*}
	On the other hand, the representation $\pi_{\mu}$ which is a map from $G$ to unitary operators on Bargmann-Fock space $\cF(\C^n)$, while satisfying the left-invariance $\pi_{\mu}((x, u)\gamma_{\nu, \lambda}(s)) =\pi_{\mu}(x, u)\pi_{\mu}(\gamma_{\nu, \lambda}(s))$ also satisfies
	\begin{align*}
		\pi_\mu(0, \pi R^2\hat{\lambda})=\pi_{\mu}(0, 0)=I
	\end{align*}
	whenever $\mu\cdot\lambda=2k|\lambda|^3$, and therefore
	\begin{align*}
		\pi_\mu((x, u)\gamma_{\nu, \lambda}(s+2\pi R))=\pi_\mu((x, u)\gamma_{\nu, \lambda}(s)).
	\end{align*}
\end{remark}
\subsection{Schr\"{o}dinger Representation on the Heisenberg group}
Denote by $\mathcal{U}(L^2(\mathbb{R}^n))$ the set of unitary operators acting on $L^2(\mathbb{R}^n)$. For each $h\in \mathbb{R}\setminus \{0\}$ the map
\begin{align*}
	\rho_h: \mathbb{H}_n\longrightarrow \mathcal{U}(L^2(\mathbb{R}^n))
\end{align*}
\begin{align*}
	[\rho_h(x, y, t)\phi](\xi)=e^{ih(t+\frac{1}{2}xy)}e^{i\sqrt{h}y\xi}\phi(\xi+\sqrt{|h|}x), \quad \phi\in L^2(\mathbb{R}^n)
\end{align*}
where $\sqrt{h}=\text{sgn}(h)\sqrt{|h|}$. Here we are more or less using the conventions of Strichartz \cite[p. 358]{strichartzLpHarmonicAnalysis1991} and \cite{fischerQuantizationNilpotentLie2016}.

This is a strongly continuous unitary representation of $\mathbb{H}_n$ on $L^2(\mathbb{R}^n)$. These \textit{Schr\"{o}dinger representations} are irreducible and pairwise inequivalent. Moreover, by the Stone-von Neumann Theorem, they exhaust all of the irreducible unitary representations, up to unitary equivalence, that are nontrivial on the center of $\mathbb{H}_n$.

\begin{remark}
	There are many other realizations of the Schr\"{o}dinger representation, but we have chosen the one for which
	\begin{itemize}
		\item $\rho_h(z, t)=e^{iht}\rho_h(z, 0)$
		\item $d\rho_h(X_j)=\sqrt{|h|}\,\partial_{\xi_j}$
		\item $d\rho_h(Y_j)=i\sgn(h)\sqrt{|h|}\,\xi_j$, $1\leq j\leq n$.
	\end{itemize}
\end{remark}
The Bargmann-Fock Representation ${\pi}_h$ is unitarily equivalent to the Schr\"{o}dinger representation $\rho_{ h}$ with the intertwining map $B: L^2(\R^n)\to \cF(\C^n)$ given by
\begin{align}\label{Bargmann}
	Bf(\zeta)=\pi^{-n/4}e^{-\frac{1}{2}\zeta\cdot\zeta}\int_{\R^n}f(\xi)e^{\sqrt{2}\xi\cdot\zeta -\frac{1}{2}|\xi|^2}d\xi
\end{align}
That is, $\pi_h(z, t)\circ B=B\circ\rho_h(z, t).$ This can be seen by taking the generating series for the Hermite functions:
\begin{align*}
	\pi^{-n/4}e^{\sqrt{2}\xi\cdot\zeta-\zeta\cdot\zeta/2-|\xi|^2/2}=\sum_{\alpha\in \N^n}\Phi_\alpha(\xi)\omega_\alpha(\zeta).
\end{align*}
\subsection{Schr\"{o}dinger Representation on $H$-type groups}
Returning to $H$-type groups, let $\mu\in \fz^*\setminus\{0\}$.
The map
\begin{align*}
	\rho_\mu:=\rho_{|\mu|}{\circ}\,\alpha_\mu
\end{align*}
is an irreducible unitary representation of $G$ acting on $L^2(\mathbb{R}^n)$. In fact, any irreducible representation of $G$ with central character $e^{i\mu\cdot u}$ is unitarily equivalent to $\pi_\mu$. All irreducible unitary representations that are nontrivial on the center of $G$ are unitarily equivalent to one such $\pi_\mu$.
\begin{remark}
	We assume it is understood that when $\mu\in \fz^*\mz$, the map $\rho_\mu$ is the Schr\"{o}dinger representation on the $H$-type group $G$; when $h=|\mu|$ is a scalar, $\rho_h$ is a representation of the Heisenberg group $\heis_n$.
\end{remark}

\subsection{The Weyl transform}\label{s:Weyl}
For $\mu\in \fz^*\mz$ and $\phi\in L^1(\fv)$ define
\begin{align}\label{e:WeylDef}
	W_\mu(\phi):=\int_\fv\phi(x)\pi_\mu(x, 0)^*dx.
\end{align}
By \eqref{e:repH}, $W_\mu(\phi)=\int_{\R^{2n}}(\phi\circ R_\mu)(z)\,\pi_{|\mu|}(z, 0)^*dz$ is the Weyl transform on $\heis_n$ at the parameter $|\mu|$ \cite{thangaveluHarmonicAnalysisHeisenberg1998, follandHarmonicAnalysisPhase1989} of $\phi\circ R_\mu$. Since $R_\mu$ is orthogonal, $\phi$ and $\phi\circ R_\mu$ have the same $L^p$ norms, so properties of the Weyl transform on $\heis_n$ transfer to $W_\mu$ verbatim: $W_\mu$ is injective on $L^1(\fv)$ \cite{thangaveluHarmonicAnalysisHeisenberg1998}. When $m=1$ and $\mu>0$ one may take $R_\mu={\rm id}$, and $W_\mu$ is the classical Weyl transform at the parameter $\mu$.

\begin{proposition}[Plancherel identity for the Weyl transform]\label{p:PlancherelWeyl}
	Let $\mu\in \fz^*\mz$ and $\phi\in L^1\cap L^2(\fv)$. Then $W_\mu(\phi)$ is Hilbert--Schmidt and
	\begin{align}\label{e:PlancherelWeyl}
		\|W_\mu(\phi)\|^2_{HS}=(2\pi)^n|\mu|^{-n}\|\phi\|^2_{L^2(\fv)}.
	\end{align}
\end{proposition}
\begin{proof}
	Write $h=|\mu|$ and $\phi_\mu:=\phi\circ R_\mu$. By \eqref{e:repH} and the unitary equivalence $B$ of \eqref{Bargmann} between $\pi_{|\mu|}$ and $\rho_{|\mu|}$, the operator $W_\mu(\phi)$ has the same Hilbert--Schmidt norm as $W_h(\phi_\mu):=\int_{\R^{2n}}\phi_\mu(z)\rho_h(z, 0)^*dz$. From the definition of $\rho_h$, $W_h(\phi_\mu)$ has integral kernel
	\begin{align*}
		K(\xi, \xi')=h^{-n/2}\int_{\R^n}\phi_\mu\Big(\frac{\xi-\xi'}{\sqrt h}, y\Big)e^{-i\sqrt h\, y\cdot\frac{\xi+\xi'}{2}}dy,
	\end{align*}
	and the change of variables $(\xi, \xi')\mapsto\big(\frac{\xi-\xi'}{\sqrt h}, \sqrt h\frac{\xi+\xi'}{2}\big)$ has unit Jacobian, so by the Euclidean Plancherel theorem in $y$,
	\begin{align*}
		\|W_h(\phi_\mu)\|_{HS}^2=h^{-n}\int_{\R^{2n}}\Big|\int_{\R^n}\phi_\mu(a, y)e^{-iy\cdot b}dy\Big|^2da\, db=(2\pi)^nh^{-n}\|\phi_\mu\|^2_{L^2(\R^{2n})}.
	\end{align*}
	Since $R_\mu$ is orthogonal, $\|\phi_\mu\|_{L^2(\R^{2n})}=\|\phi\|_{L^2(\fv)}$, which is \eqref{e:PlancherelWeyl}.
\end{proof}

\section{Group Fourier Analysis}\label{GFT}
Motivated by the Peter--Weyl Theorem for compact Lie groups, we consider the unitary dual as a generalized frequency space for the group Fourier transform. See \cite{fischerQuantizationNilpotentLie2016} for a detailed exposition on the topic, and \cite{follandHarmonicAnalysisPhase1989, gellerFourierAnalysisHeisenberg1977, thangaveluHarmonicAnalysisHeisenberg1998} for the harmonic analysis of the Heisenberg group.
For any $\pi\in \widehat{G}$, we define the group Fourier transform of an integrable function $f\in L^1(G)$ to be the operator $\cF(f)(\pi)\in \mathcal{L}(\cH_\pi)$ defined by
\begin{align*}
	\mathcal{F}(f)(\pi)=\int_Gf(g)\pi(g)^*dg.
\end{align*}
For $\mu\in \fz^*\mz$ or $\eta\in \fv^*$, we write $\cF(f)(\mu)=\cF(f)(\pi_\mu)$ and $\cF(f)(\eta, 0)=\cF(f)(\pi_{(\eta, 0)})$. For $f\in L^1(G)$, its Fourier transform $\pi\mapsto \cF(f)(\pi)$ is a bounded family of bounded operators on $\cH_\pi$, indexed by $\widehat{G}=\fv^*\cup \fz^*\mz$. The Fourier transform extends to an isometry from $L^2(G)$ onto the Hilbert space $L^2(\widehat{G})$ of all measurable families $\{A(\mu)\}_{\mu\in \fz^*\mz}$ of operators on $\cH_{\pi_\mu}$ which are Hilbert-Schmidt for almost every $\mu\in \fz^*\mz$, with norm
\begin{align*}
	\|A\|:=\left(\int_{\fz^*\mz}\|A(\mu)\|^2_{HS(\cH_\mu)}d\pi(\mu)\right)^{\frac{1}{2}}<\infty
\end{align*}
where $d\pi(\mu)=c_0|\mu|^nd\mu$, for a constant $c_0>0$ depending on $n$ and $m$, is the \textit{Plancherel measure}.
This is the unique measure on $\widehat{G}$ for which the Plancherel formula holds:
\begin{align*}
	\int_G |f(g)|^2dg=c_0\int_{\fz^*\mz}\|\cF(f)(\mu)\|^2_{HS(\cH_\mu)}|\mu|^nd\mu.
\end{align*}
Note that only the $\fz^*\mz$ part of the unitary dual appears in the Plancherel formula. This is because the subset of $\widehat{G}$ corresponding to $\fv^*$ of finite dimensional representations has measure zero.
The Plancherel Formula yields an inversion formula for $f\in \mathcal{S}(G)$:
\begin{align*}
	f(g)=\int_{\widehat{G}}{\rm tr}\left(\pi(g)\mathcal{F}(f)(\pi)\right)d\pi.
\end{align*}
\subsection{Fourier Transform on $G_\lambda$}
Fix $\lambda\in \fz^*\mz$. For $f\in C_c^\infty(G_\lambda)$ and $\mu\in \fz^*\setminus\{0\}$ satisfying \eqref{compadability} let
\begin{align*}
	\mathcal{F}_\lambda(f)(\pi_\mu)=\int_{G_\lambda}f(g)\pi_\mu^*(g)dg
\end{align*}
Note that the Fourier transform $\cF_\lambda(f)(\pi_\mu)$ on $G_\lambda$ is only defined for $\mu\in \Gamma_\lambda^*$.

For $\eta\in \fv^*$ we have
\begin{align*}
	\cF_\lambda(f)(\pi_{(\eta, 0)})=\int_{\fv}f_0(x)e^{-i\langle \eta, x\rangle}dx
\end{align*}
where $f_0(x)=\int_{\fz/\Gamma_\lambda}f(x, u)du$.

\begin{theorem}[Plancherel theorem on $G_\lambda$]\label{t:PlancherelQuotient}
	Let $F\in L^1(G_\lambda)\cap L^2(G_\lambda)$ and, for $\mu\in \Gamma_\lambda^*$, let $F^\mu(x)=\int_{\fz/\Gamma_\lambda}F(x, u)e^{-i\mu\cdot u}du$ as in \eqref{e:centralFT}. Write $\mu=\mu_\perp+2k|\lambda|^2\hat\lambda$ with $\mu_\perp\in \lambda^\perp$ and $k\in \Z$. Then
	\begin{align}\label{e:PlancherelCentral}
		\int_{G_\lambda}|F(x, u)|^2d(x, u)=\frac{|\lambda|^2}{\pi(2\pi)^{m-1}}\sum_{k\in \Z}\int_{\lambda^\perp}\big\|F^{\mu_\perp+2k|\lambda|^2\hat\lambda}\big\|^2_{L^2(\fv)}d\mu_\perp,
	\end{align}
	where for $m=1$ the integral over $\lambda^\perp=\{0\}$ is evaluation at $\mu_\perp=0$. Moreover, for almost every $\mu\in \Gamma_\lambda^*\mz$ (every $\mu$ when $m=1$), both sides finite,
	\begin{align}\label{e:PlancherelHS}
		\|F^\mu\|^2_{L^2(\fv)}=(2\pi)^{-n}|\mu|^n\,\|\cF_\lambda(F)(\pi_\mu)\|^2_{HS},
	\end{align}
	and $\|F^0\|^2_{L^2(\fv)}=(2\pi)^{-2n}\int_{\fv^*}|\cF_\lambda(F)(\pi_{(\eta, 0)})|^2d\eta$.
\end{theorem}
\begin{proof}
	For a.e. $x\in \fv$ the function $F(x, \cdot)$ lies in $L^1\cap L^2(\fz/\Gamma_\lambda)$. Since $\fz/\Gamma_\lambda\cong\lambda^\perp\times\R/\pi|\lambda|^{-2}\Z$ has dual group $\Gamma_\lambda^*$, Parseval on $\lambda^\perp$ and on the circle of length $\pi|\lambda|^{-2}$ give
	\begin{align*}
		\int_{\fz/\Gamma_\lambda}|F(x, u)|^2du=\frac{1}{(2\pi)^{m-1}}\cdot\frac{|\lambda|^2}{\pi}\sum_{k\in \Z}\int_{\lambda^\perp}\big|F^{\mu_\perp+2k|\lambda|^2\hat\lambda}(x)\big|^2d\mu_\perp,
	\end{align*}
	and integrating in $x$ gives \eqref{e:PlancherelCentral}. For $\mu\in \Gamma_\lambda^*\mz$, Fubini and $\pi_\mu(x, u)=e^{i\mu\cdot u}\pi_\mu(x, 0)$ give $\cF_\lambda(F)(\pi_\mu)=\int_\fv F^\mu(x)\pi_\mu(x, 0)^*dx=W_\mu(F^\mu)$, and $F^\mu\in L^1(\fv)$ since $\|F^\mu\|_{L^1(\fv)}\leq\|F\|_{L^1(G_\lambda)}$, while $F^\mu\in L^2(\fv)$ for almost every $\mu$ by \eqref{e:PlancherelCentral} (and for every $\mu$ when $m=1$), so \eqref{e:PlancherelHS} is Proposition \ref{p:PlancherelWeyl} applied to $\phi=F^\mu$. The formula for $\|F^0\|$ is the Euclidean Plancherel theorem on $\fv$.
\end{proof}
\begin{corollary}\label{c:uniquenessQuotient}
	Let $F\in L^1(G_\lambda)$. If $F^\mu=0$ for every $\mu\in \Gamma_\lambda^*$, then $F=0$. In particular $F=0$ whenever $\cF_\lambda(F)(\pi)=0$ for all $\pi\in \widehat{G}_\lambda$.
\end{corollary}
\begin{proof}
	For a.e. $x$, $F(x, \cdot)\in L^1(\fz/\Gamma_\lambda)$ has vanishing Fourier transform on the dual group $\Gamma_\lambda^*$, hence vanishes a.e. If $\cF_\lambda(F)(\pi)=0$ for all $\pi\in \widehat{G}_\lambda$, then $F^0=0$ by Fourier uniqueness on $\fv$, and for $\mu\in \Gamma_\lambda^*\mz$, $W_\mu(F^\mu)=\cF_\lambda(F)(\pi_\mu)=0$ with $F^\mu\in L^1(\fv)$, so $F^\mu=0$ by the injectivity of $W_\mu$ on $L^1(\fv)$ (Section \ref{s:Weyl}).
\end{proof}

\subsection{Lemmas on the Fourier Transform}
The proof of the following lemma is immediate.
\begin{lemma}\label{dilationProperty}
	Fix $\eps>0$. For $f\in L^1(G)$,
\begin{align*}
	\cF(\delta_\eps^* f)(\mu)=\eps^{-(2n+2m)}\cF(f)(\mu/\eps^2), \quad \forall \mu\in \fz^*\mz,
\end{align*}
and for $f\in L^1(G/\Gamma_\lambda)$,
\begin{align*}
	\cF_{\eps\lambda}(\delta^*_\eps f)(\mu)=\eps^{-(2n+2m)}\cF_{\lambda}(f)(\mu/\eps^2), \quad \forall \mu\in \Gamma_{\eps\lambda}^*\mz.
\end{align*}
\end{lemma}
\begin{lemma}\label{PoissonSum}
	Fix $\lambda\in \mathfrak{z}^*\setminus\{0\}$, then for all $f\in L^1(G)$,
	\begin{align*}
		\mathcal{F}_\lambda\left(P_\lambda f\right)(\pi)
		=\mathcal{F}(f)(\pi),
		\quad \forall \pi\in \widehat{G}_\lambda.
	\end{align*}
	Here we use the abuse of notation $\pi=\pi\circ q_\lambda$, where $q_\lambda: G\to G_\lambda$ is the quotient map.
\end{lemma}
\begin{proof}
	For any $\phi_1, \phi_2\in \cH_\pi$ by Equation \eqref{Haar2}, we have
	\begin{align*}
		\int_{G_\lambda}P_\lambda f(x, u)\langle \pi(x, u)^*\phi_1, \phi_2\rangle dg_\lambda(x, u)
		=\int_G f(x, u)\langle \pi(x, u)^*\phi_1, \phi_2\rangle dg(x, u)
	\end{align*}
	The lemma follows.
\end{proof}
\begin{corollary}
	Fix $\lambda\in \mathfrak{z}^*\setminus\{0\}$, then for all $f\in L^1(G)$,
	\begin{align*}
		\mathcal{F}_\lambda\left(P_\lambda f\right)(\pi_{(\eta, 0)})
		=\mathcal{F}(f)(\pi_{(\eta, 0)}),
		\quad \forall \eta\in \fv^*.
	\end{align*}
\end{corollary}

\begin{corollary}
	Fix $\lambda\in \mathfrak{z}^*\setminus\{0\}$, then for all $f\in L^1(G)$,
	\begin{align*}
		\mathcal{F}_\lambda\left(P_\lambda f\right)(\pi_\mu)
		=\mathcal{F}(f)(\pi_\mu),
		\quad \forall \mu\in \fz^*\setminus\{0\} \text{ with }\mu\cdot \lambda=2k|\lambda|^3 \text{ and } k\in \mathbb{Z}.
	\end{align*}
\end{corollary}

\subsection{Convolution Operators on $G_\lambda$}
Define convolution	 on $G_\lambda$ by
\begin{align*}
	f*g(x)=\int_{G_\lambda}f(y)g(y^{-1}x)dy, \quad f, g\in C_c^\infty(G_\lambda)
\end{align*}
Let $\kappa\in \cE'(G_\lambda)$ be a compactly supported distribution on the Lie group $G_\lambda$. Convolution of a function with a distribution $f\in C_c^\infty(G_\lambda)$ is defined by
\begin{align*}
	f*\kappa (x)=\langle \kappa, \pi_R(x^{-1})\widetilde{f}\rangle
\end{align*}
where $\pi_R(x)f(y):=f(yx)$ and $\widetilde{f}(x):=f(x^{-1})$, hence $\pi_R(x^{-1})\widetilde{f}(y)=f(xy^{-1})$.

Suppose $A: C_c^\infty (G_\lambda)\to C_c^\infty(G_\lambda)$ is the right convolution operator given by
\begin{align*}
	Af=f*\kappa.
\end{align*}
Then there is a family of operators $\hat\kappa(\pi): \cH_\pi\to \cH_\pi$ indexed by $\pi\in \widehat{G}_\lambda$ such that
\begin{align*}
	\cF_\lambda(Af)(\pi)=\hat{\kappa}(\pi)\circ\cF_\lambda(f)(\pi), \quad \pi \in \widehat{G}_\lambda.
\end{align*}
In many contexts, the \textit{field of operators} $\{\hat{\kappa}(\pi): \pi\in \widehat{G}\}$ is called the \textit{operator-valued symbol} of $A$ \cite{fischerQuantizationNilpotentLie2016}.

If $\kappa=\int_{G_\lambda} d\sigma$ is just integration against a compactly supported Borel measure $\sigma$ on $G_\lambda$, then the multiplier of the convolution operator $A$ is given by
\begin{align*}
	\hat{\kappa}(\pi)=\int_{G_\lambda}\pi(x)^*d\sigma(x).
\end{align*}
The X-ray Transform $\widetilde{I}_{\nu, \lambda}$, in particular, is a right convolution operator on $G_\lambda$:
\begin{align}\label{convOp}
	\widetilde{I}_{\nu, \lambda}f=f*\kappa_{\nu, \lambda}, \quad f\in C_c^\infty(G_\lambda),
\end{align}
where $\kappa_{\nu, \lambda}$ is integration against the Dirac measure on the curve $\gamma_{\nu, \lambda}(s)^{-1}$:
\begin{align*}
	\langle \kappa_{\nu, \lambda}, f\rangle=\int_{G_\lambda}f(x)\delta_{\gamma_{\nu, \lambda}^{-1}}(x)dx=\int_0^{2\pi|\lambda|^{-1}}f(\gamma_{\nu, \lambda}(s)^{-1})ds.
\end{align*}
Therefore
\begin{align}
	\hat\kappa_{\nu, \lambda}(\pi)=\int_0^{2\pi|\lambda|^{-1}}\pi(\gamma_{\nu, \lambda}(s))ds, \quad \pi\in \widehat{G}_\lambda. \label{multiplier}
\end{align}
\subsection{Homogeneity of the symbol}
The analog of Proposition \ref{homogeneity} in frequency space is the following proposition, which describes the homogeneity of the symbol:
\begin{proposition}\label{homofourier}
	Let $\kappa_{\nu, \lambda}$ be the convolution kernel given in \eqref{convOp} of the holonomy transform $\widetilde{I}_{\nu, \lambda}$. Then for all $\eps>0$ and all $\mu\in \Gamma_\lambda^*\mz$,
	\begin{align*}
		\hat{\kappa}_{\nu, \lambda}(\mu)=\eps \hat{\kappa}_{\nu, \eps\lambda}(\eps^2\mu).
	\end{align*}
\end{proposition}
\begin{proof}
	By \eqref{multiplier}, $\pi_{\eps^2\mu}=\pi_\mu\circ\delta_\eps$ (Lemma \ref{dilationProperty}) and $\delta_\eps\gamma_{\nu, \eps\lambda}(s)=\gamma_{\nu, \lambda}(\eps s)$ (proof of Proposition \ref{homogeneity}),
	\begin{align*}
		\hat{\kappa}_{\nu, \eps\lambda}(\eps^2\mu)&=\int_0^{2\pi|\eps\lambda|^{-1}}\pi_\mu\big(\gamma_{\nu, \lambda}(\eps s)\big)\, ds=\eps^{-1}\int_0^{2\pi|\lambda|^{-1}}\pi_\mu\big(\gamma_{\nu, \lambda}(s)\big)\, ds=\eps^{-1}\hat{\kappa}_{\nu, \lambda}(\mu).\qedhere
	\end{align*}
\end{proof}

\subsection{Fourier Multipliers}
\begin{lemma}\label{l:scalarMultiplier} For $\lambda\in \fz^*\mz$, $\nu\in S(\fv^*)$ and $\eta\in \fv^*$, let $\eta_{\nu, \lambda}$ be the orthogonal projection of $\eta$ onto ${\rm Span}\{\nu, J_{\hat\lambda}\nu\}$. Then
	\begin{align*}
		\mathcal{F}_\lambda\left(\widetilde{I}_{\nu, \lambda} f\right)(\pi_{(\eta, 0)})=2\pi|\lambda|^{-1}J_0(|\eta_{\nu, \lambda}|\,|\lambda|^{-1})\mathcal{F}_\lambda(f)(\pi_{(\eta, 0)}), \quad f\in L^1(G_\lambda),
	\end{align*}
\end{lemma}
\begin{proof}
	By equation \eqref{convOp} and \eqref{multiplier}, we are interested in the quantity
	\begin{align*}
		\int_0^{2\pi|\lambda|^{-1}}\pi_{(\eta, 0)}(\gamma_{\nu, \lambda}(s))ds=
	\int_0^{2\pi|\lambda|^{-1}}e^{i\langle \eta, J_\lambda^{-1}e^{sJ_\lambda}\nu\rangle}ds.
	\end{align*}
With $R_\lambda\in O(\R^{2n}, \fv)$ as in \eqref{normalForm}, $\widetilde{\eta}=R_\lambda^t\eta$ and $\widetilde{\nu}=R_\lambda^t\nu$,
\begin{align*}
	\langle \eta, J_\lambda^{-1}e^{sJ_\lambda}\nu\rangle
	=|\lambda|^{-1}\langle \widetilde{\eta}, J^{-1}e^{s|\lambda|J}\widetilde{\nu}\rangle=-|\lambda|^{-1}\langle \widetilde{\eta}, Je^{s|\lambda|J}\widetilde{\nu}\rangle.
\end{align*}
Let $V={{\rm Span}}_\R\{\widetilde{\nu},\; J\widetilde{\nu}\}=R_\lambda^t\,{\rm Span}\{\nu, J_{\hat\lambda}\nu\}$ and $\widetilde{\eta}=\widetilde{\eta}_V+\widetilde{\eta}_\perp$ be the orthogonal decomposition with respect to $V$ and $V^\perp$, so that $|\widetilde{\eta}_V|=|\eta_{\nu, \lambda}|$. Since $Je^{s|\lambda|J}\widetilde{\nu}\in V$, writing $\widetilde{\eta}_V=|\widetilde{\eta}_V|e^{\theta J}\widetilde{\nu}$ gives
\begin{align*}
	\langle \widetilde{\eta}, Je^{s|\lambda| J}\widetilde{\nu}\rangle=|\widetilde{\eta}_V|\langle \widetilde{\nu}, Je^{(s|\lambda|-\theta)J}\widetilde{\nu}\rangle=-|\widetilde{\eta}_V|\sin{(s|\lambda|-\theta)}.
\end{align*}
Substituting $s\mapsto s/|\lambda|$,
\begin{align*}
	\int_0^{2\pi|\lambda|^{-1}}e^{i|\lambda|^{-1}|\eta_{\nu, \lambda}| \sin{(s|\lambda|-\theta)}}ds=|\lambda|^{-1}\int_{-\pi}^{\pi}e^{i|\lambda|^{-1}|\eta_{\nu, \lambda}| \sin{s}}ds=2\pi|\lambda|^{-1}J_0(|\lambda|^{-1}|\eta_{\nu, \lambda}|)
\end{align*}
as desired.
\end{proof}
\begin{definition}\label{d:J}
	For $\lambda, \mu\in \fz^*\setminus\{0\}$ with $\mu\cdot\lambda =2k|\lambda|^3$ for some $k\in \mathbb{Z}$ and $\nu\in S(\fv^*)$ define
	\begin{align*}
		\cJ_{\nu, \lambda}(\mu)
		&=\frac{1}{2\pi}\int_0^{2\pi }e^{iks}\pi_\mu\left(R_\lambda\frac{ e^{sJ}}{|\lambda|J}R_\lambda^{{\rm{t}}}\nu, 0\right)ds
	\end{align*}
\end{definition}

\begin{lemma}\label{FourierSliceLemma} For $f\in L^1(G_\lambda)$, and $\lambda, \mu\in \mathfrak{z}^*\setminus\{0\}$ with $\lambda\cdot\mu=2k|\lambda|^3$, $k\in \mathbb{Z}$,
	\begin{align*}
	\mathcal{F}_\lambda\left(\widetilde{I}_{\nu, \lambda} f\right)(\pi_\mu)=2\pi|\lambda|^{-1}\mathcal{J}_{\nu, \lambda}( \mu)\circ\mathcal{F}_\lambda(f)(\pi_\mu)
	\end{align*}
\end{lemma}
\begin{proof}

	By \eqref{multiplier}, we just need to compute
	\begin{align*}
	\hat{\kappa}_{\nu, \lambda}(\mu)
	&=\int_0^{2\pi|\lambda|^{-1}}\pi_\mu(\gamma_{\nu, \lambda}(s))ds\\
	&=\int_0^{2\pi |\lambda|^{-1}}\pi_\mu\left(\frac{e^{sJ_\lambda}}{J_\lambda}\nu, \frac{s}{2|\lambda|^2}\lambda^\sharp\right)ds\\
	&=\int_0^{2\pi |\lambda|^{-1}}e^{ik|\lambda|s}\pi_\mu\left(\frac{e^{sJ_\lambda}}{J_\lambda}\nu, 0\right)ds\\
	&=\frac{1}{|\lambda|}\int_0^{2\pi }e^{iks}\pi_\mu\left(\frac{e^{sJ_\lambda/|\lambda|}}{J_\lambda}\nu, 0\right)ds\\
	&:=2\pi|\lambda|^{-1}\mathcal{J}_{\nu, \lambda}(\mu)
	\end{align*}
with $\gamma_{\nu, \lambda}$ given by \eqref{param}.
\end{proof}

\begin{remark}\label{r:scaleJ}

	It is clear from a direct computation, or from Proposition \ref{homofourier} that for $\eps>0$
	\begin{align*}
		\cJ_{\nu, \eps\lambda}(\eps^2\mu)=\cJ_{\nu, \lambda}(\mu), \quad \mu\in \Gamma_\lambda^*\mz.
	\end{align*}
	Therefore, we may assume without loss of generality that $|\lambda|=1$, and later extend using homogeneity.
\end{remark}

\begin{proof}[Proof of Theorem \ref{FourierSliceTheorem}]
	The Fourier slice theorem is now a formal consequence of Lemmas \ref{FourierSliceLemma} and \ref{PoissonSum}, and Proposition \ref{factorization}. Let $f\in L^1(G)$, $\nu\in S(\fv^*)$ and $\lambda\in \fz^*\mz$ be given.

	For all $\eta\in \fv^*$,
	\begin{align*}
		\cF_\lambda(I_{\nu, \lambda}f)(\eta, 0)
		&=\cF_\lambda(\widetilde{I}_{\nu, \lambda}\circ P_\lambda f)(\eta, 0)\\
		&=2\pi|\lambda|^{-1}J_0(|\eta_{\nu, \lambda}|\,|\lambda|^{-1})\cF_\lambda(P_\lambda f)(\eta, 0)\\
		&=2\pi|\lambda|^{-1}J_0(|\eta_{\nu, \lambda}|\,|\lambda|^{-1})\cF( f)(\eta, 0).
	\end{align*}
	 Similarly, for all $\mu\in \Gamma_\lambda^*\mz$,
	\begin{align*}
		\cF_{\lambda}(I_{\nu, \lambda}f)(\mu)
		&=\cF_\lambda(\widetilde{I}_{\nu, \lambda}\circ P_\lambda f)(\mu)\\
		&=2\pi|\lambda|^{-1}\cJ_{\nu, \lambda}(\mu)\circ \cF_\lambda(P_\lambda f)(\mu)\\
		&=2\pi|\lambda|^{-1}\cJ_{\nu, \lambda}(\mu)\circ \cF( f)(\mu)
	\end{align*}
as desired.
\end{proof}

\begin{remark}\label{r:parallelTransport}
	We may view the weight in \eqref{e:twistedSlice} in terms of parallel transport. For $\mu\in\Gamma_\lambda^*\mz$ with $\mu\cdot\lambda=2k|\lambda|^3$, set $x(s):=J_\lambda^{-1}e^{sJ_\lambda}\nu$, as in Remark \ref{r:twistedSlice}. Let $A_\mu:=\tfrac12\omega_\mu(y, dy)$, a one-form on $\fv$ with $dA_\mu=\omega_\mu$. By \eqref{e:htypePolar}, $\tfrac12\omega_\mu(x(s), \dot x(s))=\tfrac12|\lambda|^{-2}\mu\cdot\lambda=k|\lambda|$, so for $y(s):=x+x(s)$,
	\begin{align*}
		e^{ik|\lambda|s}e^{\frac{i}{2}\omega_\mu(x, x(s))}=e^{\frac{i}{2}\omega_\mu(x, x(0))}\exp\Big(i\int_0^sA_\mu(\dot y(t))\,dt\Big).
	\end{align*}
	Thus, the central Fourier mode of $I_{\nu, \lambda}f$ at frequency $\mu$ is the ray transform of $f^\mu$ along the magnetic geodesics $y(s)$ of charge $\lambda$ in $\fv$ \cite{dairbekovBoundaryRigidityMagnetic2007}, attenuated by the parallel transport of the unitary connection $d-iA_\mu$, whose curvature is $-i\omega_\mu$ \cite{holmanAttenuatedGeodesicXray2018, paternainGeometricInverseProblems2023}. %
\end{remark}

\begin{remark}\label{commutative}
	Suppose that $G=\heis_1$ is the first Heisenberg group. Then the Fourier slice theorem is essentially the statement that the following diagram commutes:
	\[
	\begin{tikzcd}
		L^1(\heis_1) \arrow[dd, bend right=50, "I_{\nu, \lambda}"']\arrow[d, "P_\lambda"] \arrow[r, "\mathcal{F}"] &C\left(\mathbb{R}^*, S_\infty\right)\arrow[d, "\text{res}"]\\
		L^1(\heis_1/\Gamma_\lambda) \arrow[r, "\mathcal{F}_{\lambda}"] \arrow[d, "\widetilde{I}_{\nu, \lambda}"]& \bigoplus_{\mu\in\Gamma_\lambda^*\mz}S_\infty\arrow[d, "2\pi|\lambda|^{-1}\mathcal{J}_{\nu, \lambda}(\mu)"]\\
		L^1(\heis_1/\Gamma_\lambda) \arrow[r, "\mathcal{F}_{\lambda}"] & \bigoplus_{\mu\in\Gamma_\lambda^*\mz}S_\infty
	\end{tikzcd}
	\]
	Here $S_\infty$ is the space of compact operators on $\cF(\C)$, $C(\R^*, S_\infty)$ is the set of continuous maps from $\R^*\to S_\infty$, and the sums run over $\Gamma_\lambda^*\mz=2\lambda^2\Z^*$.
\end{remark}

\section{Proofs of the Main Theorems}\label{s:proofs}
The symbol is computed in the Bargmann--Fock representation, and the computation splits into two cases. It is the average of the representation over one cycle of the helix, a loop in the plane $P_\nu\subset\fv$ spanned by the direction $\nu$ and its rotation $J_{\hat\lambda}\nu$. The calculation of the average is split into cases depending on the K\"{a}hler angle of $P_\nu$ with respect to the complex structure $J_{\hat\mu}$ that the frequency induces on $\fv$. The cosine of this angle is $\hat\mu\cdot\hat\lambda$, which on $\Gamma_\lambda^*$ takes the values $2k|\lambda|^2/|\mu|$, $k\in \Z$ (Remark \ref{r:kahlerAngle}). When the frequency is not orthogonal to the charge the plane is symplectic for $\omega_\mu$, a metaplectic intertwiner moves the loop into a complex line. Then the symbol becomes a weighted shift on the Bargmann--Fock basis with Laguerre-function weights, the computation of \cite{flynnInjectivityHeisenbergXray2021} on the Heisenberg group (Proposition \ref{specDecompShift}). When the frequency is orthogonal to the charge, which occurs only for $\dim\fz\geq2$, the plane is isotropic for $\omega_\mu$, the loop can be moved into a Lagrangian subspace. Then in the Schr\"{o}dinger model the symbol is multiplication by the Bessel function $J_0$, exactly as for the Euclidean mean value transform of Section \ref{MainResults} (Proposition \ref{p:specDecompPerp}). This is why we work with both models. A rotation of a complex line acts on Bargmann--Fock space by precomposition by a rotation on $\C^n$, whereas translation along a Lagrangian subspace acts in the Schr\"{o}dinger model by multiplication by characters, and the Bargmann transform passes between the two.

\subsection{Entry Functions}
Here we follow \cite[page 18]{thangaveluHarmonicAnalysisHeisenberg1998}.
For $k=0, 1, 2 \ldots $, let $$\phi_k(x):=(2^k\sqrt{\pi}k!)^{-\frac{1}{2}}(-1)^ke^{x^2/2}\frac{d^k}{dx^k}e^{-x^2}$$
be the $k$th Hermite function. For a multi-index $\alpha\in \N^n$ define the multi-Hermite function
\begin{align*}
	\Phi_\alpha(x)=\prod_{j=1}^n\phi_{\alpha_j}(x_j).
\end{align*}
They are eigenfunctions for the Hermite operator $H=-\Delta+|x|^2$: $H\Phi_\alpha=(2|\alpha|+n)\Phi_\alpha$ and for the Fourier transform: $\mathcal{F}(\Phi_\alpha)=(-i)^{|\alpha|}\Phi_\alpha$.

For $\alpha, \beta\in \N^n$ the special Hermite functions are
\begin{align}\label{specialHermite}
	\Phi_{\alpha, \beta}(z)=(2\pi)^{-\frac{n}{2}}\langle \rho(z)\Phi_\alpha, \Phi_\beta\rangle , \quad z\in \C^n.
\end{align}
The special Hermite functions for $n=1$ are related to the general $n\geq 1$ case via
\begin{align*}
	\Phi_{\alpha, \beta}(z)=\prod_{j=1}^n \Phi_{\alpha_j, \beta_j}(z_j), \quad z\in \mathbb{C}^n.
\end{align*}
Explicitly, a computation for $a, b\in \N$, and $z\in \C$ gives:
\begin{align*}
	\Phi_{a, b}(z)=(2\pi)^{-\frac{1}{2}}
	\begin{cases}
		\sqrt{\frac{b!}{a!}}\left(\frac{+1}{\sqrt{2}}z\right)^{a-b}L_b^{(a-b)}(\frac{1}{2}|z|^2)e^{-|z|^2/4} & a\geq b\\
		\sqrt{\frac{a!}{b!}}\left(\frac{-1}{\sqrt{2}}\overline{z}\right)^{b-a}L_a^{(b-a)}(\frac{1}{2}|z|^2)e^{-|z|^2/4} & a\leq b.
	\end{cases}
\end{align*}
Here $L_a^b(x)$ is the \textit{generalized Laguerre polynomial.}
The \textit{entry functions}, or \textit{matrix coefficients} on the $n$th Heisenberg group $\heis_n$ are given for $h>0$ by
\begin{align}\label{entryFunc}
	E_{\alpha\beta}^h(z, t)=\langle \rho_h(z, t)\Phi_\alpha, \Phi_\beta\rangle .
\end{align}
They are related to special Hermite functions by
\begin{align*}
	E^h_{\alpha\beta}(z, t)=(2\pi)^{\frac{n}{2}}\Phi_{\alpha, \beta}(\sqrt{h}z)e^{iht}.
\end{align*}
\begin{remark}
	When $h=1$, we will omit the superscript in the entry functions.
\end{remark} What follows are useful symmetries of the entry functions.
\begin{lemma}\label{cylinSym} For $\alpha, \beta\in \N^n$ and $h>0$ we have
	\begin{align*}
		E_{\alpha\beta}^h(e^{sJ}z, t)=e^{is|\alpha-\beta|}E_{\alpha\beta}^h(z, t)
	\end{align*}
	for all $(z, t)\in \heis_n$. Here $|\alpha-\beta|=\sum_{i}\alpha_i-\beta_i$.
	Furthermore
	\begin{align*}
		E_{\alpha\beta}^h(J^{-1}z, t)=(-i)^{|\alpha-\beta|}E_{\alpha\beta}^h(z, t).
	\end{align*}
\end{lemma}

\subsection{Spectral Decomposition of $\cJ_{\nu, \lambda}(\mu)$}
We begin by computing the spectral decomposition of $\cJ_{\nu, \lambda}(\mu)$ in the special case when $\mu\parallel \lambda$.
\begin{lemma}\label{l:homoH}
    Let $\mu, \lambda\in \fz^*\mz$ with $\mu\parallel \lambda$, and let $c=\hat{\mu}\cdot\hat{\lambda}=\pm 1$. Then $(z, t)\mapsto (R_\mu^tR_\lambda z, ct)$ is an automorphism of $\heis_n$.
\end{lemma}
\begin{proof}
    Write $A:=R_\mu^tR_\lambda$. By the group law \eqref{product}, the map $(z, t)\mapsto (Az, ct)$ is an automorphism of $\heis_n$ if and only if $A^tJA=cJ$. Since $\mu\mapsto J_\mu$ is linear and $\hat{\mu}=c\hat{\lambda}$, we have $J_{\hat\mu}=cJ_{\hat\lambda}$, so \eqref{normalForm} gives $R_\mu J R_\mu^t=cR_\lambda J R_\lambda^t$. Conjugating by $R_\mu^t$ and using $c^2=1$ gives $AJA^t=cJ$, and conjugating by $A^t\in O(2n)$ gives $A^tJA=cJ$.
\end{proof}
Recall from Definition \ref{d:tau} that for $U\in U(n)$ we take $\tau(U)F(\zeta)=F(U^*\zeta)$, $F\in \cF(\C^n)$.

When $\hat{\mu}=\hat{\lambda}$, Lemma \ref{l:homoH} gives $A^tJA=J$ for the orthogonal map $A=R_\mu^tR_\lambda$, so $R_\mu^tR_\lambda\in U(n)$. Where $\{\omega_\alpha\}_{\alpha\in \N^n}$ is the standard orthonormal basis of $\cF(\C^n)$, for $\mu, \lambda\in \fz^*\mz$ such that $\hat{\mu}=\hat{\lambda}$, define
\begin{align*}
    \omega^{\mu\lambda}_\alpha:=\tau(R_\mu^tR_\lambda)\omega_\alpha, \quad \alpha\in \N^n.
\end{align*}
When $\hat{\mu}=-\hat{\lambda}$, Lemma \ref{l:homoH} gives $A^tJA=-J$, so $A\notin U(n)$ and no such basis is defined; that case is covered by Proposition \ref{specDecompShift}, and for the normal operator it reduces to $\mu\cdot\lambda>0$ by the symmetry $c_l(-k, r)=c_l(k, r)$ used in the proof of Corollary \ref{invCriteria}.
\begin{proposition}[Spectral Decomposition of $\mathcal{J}_{\nu, \lambda}(\mu)$]\label{specDecomp}
Fix $\lambda\in S(\fz^*)$ and $\nu\in S(\fv^*)$ and $\mu\in \Gamma_\lambda^*\mz$ such that $\mu\parallel\lambda$ and $\mu\cdot\lambda>0$. With $w=\sqrt{|\mu|}R_\lambda^t\nu$, and $\cJ_{\nu, \lambda}(\mu)$ given in Definition \ref{d:J},
\begin{align*}
\cJ_{\nu, \lambda}(\mu)\omega_\alpha^{\mu\lambda}=(2\pi)^{\frac{n}{2}}i^k\sum_{|\beta|=|\alpha|+k}\Phi_{\alpha\beta}(w)\omega_\beta^{\mu\lambda}
\end{align*}
where $\Phi_{\alpha\beta}$ are the special Hermite functions given by \eqref{specialHermite}.
\end{proposition}
\begin{proof}
	Write $A:=R_\mu^tR_\lambda$ and let $k>0$ be given by $\mu\cdot\lambda=2k$. Since $\hat{\mu}=\hat{\lambda}$, $A\in U(n)$ and $\{\omega^{\mu\lambda}_\alpha\}_{\alpha\in \N^n}$ is an orthonormal basis of $\cF(\C^n)$.
	By Definition \ref{d:J} and $\pi_\mu(x, 0)=\pi_{|\mu|}(R_\mu^tx, 0)$,
	\begin{align*}
		\cJ_{\nu, \lambda}(\mu)=\frac{1}{2\pi}\int_0^{2\pi}e^{iks}\pi_{|\mu|}\left(Ae^{sJ}J^{-1}R_\lambda^t\nu, 0\right)ds.
	\end{align*}
	By Lemma \ref{l:rotateAxis}, $\tau(A)^*\circ \pi_{|\mu|}(Az, 0)\circ \tau(A)=\pi_{|\mu|}(z, 0)$, and since $B\Phi_\alpha=\omega_\alpha$ for the Bargmann transform \eqref{Bargmann}, $\langle \pi_{|\mu|}(z, 0)\omega_\alpha, \omega_\beta\rangle=E^{|\mu|}_{\alpha\beta}(z, 0)$. Therefore
	\begin{align*}
		\langle \cJ_{\nu, \lambda}(\mu)\omega^{\mu\lambda}_\alpha, \omega^{\mu\lambda}_\beta\rangle
		=\frac{1}{2\pi}\int_0^{2\pi}e^{iks}E^{|\mu|}_{\alpha\beta}\left(e^{sJ}J^{-1}R_\lambda^t\nu, 0\right)ds.
	\end{align*}
	By Lemma \ref{cylinSym},
	\begin{align*}
		E^{|\mu|}_{\alpha\beta}\left(e^{sJ}J^{-1}R_\lambda^t\nu, 0\right)=(-i)^{|\alpha-\beta|}e^{is|\alpha-\beta|}E^{|\mu|}_{\alpha\beta}\left(R_\lambda^t\nu, 0\right),
	\end{align*}
	and $\int_0^{2\pi}e^{i(k+|\alpha-\beta|)s}ds=2\pi\delta(k+|\alpha-\beta|)$, so the entry vanishes unless $|\beta|=|\alpha|+k$, in which case $(-i)^{|\alpha-\beta|}=i^k$. Since $E^{|\mu|}_{\alpha\beta}(R_\lambda^t\nu, 0)=(2\pi)^{\frac{n}{2}}\Phi_{\alpha\beta}(w)$,
	\begin{align*}
		\langle \cJ_{\nu, \lambda}(\mu)\omega^{\mu\lambda}_\alpha, \omega^{\mu\lambda}_\beta\rangle&=(2\pi)^{\frac{n}{2}}i^k\,\Phi_{\alpha\beta}(w)\,\delta(k+|\alpha-\beta|).\qedhere
	\end{align*}
\end{proof}
Let $\cE_l:={ \rm Span}_{\C}\{\omega_\alpha: |\alpha|=l\}$ and ${{\rm Pr}}_l: \cF(\C^n)\to \cE_l$ be the orthogonal projection.
\begin{corollary}
	For all $l\in \N$, we have $\cJ_{\nu, \lambda}(\mu): \cE_l\to \cE_{l+k}$. That is $\cJ_{\nu, \lambda}(\mu)$ is lower block-diagonal with blocks given by
	\begin{align*}
		\cJ_{\nu, \lambda}(\mu)=i^k\sum_{l=0}^\infty {{\rm Pr}}_{l+k}\circ \pi_\mu(\nu, 0)\circ {{\rm Pr}}_{l}.
	\end{align*}
\end{corollary}

\begin{proof}
	Write $A:=R_\mu^tR_\lambda\in U(n)$. Since $\tau(A)$ preserves each $\cE_l$ (Lemma \ref{l:rotateAxis}), $\{\omega^{\mu\lambda}_\alpha: |\alpha|=l\}$ is an orthonormal basis of $\cE_l$. By Lemma \ref{l:rotateAxis} and $R_\mu^t\nu=AR_\lambda^t\nu$,
	\begin{align*}
		\langle \pi_\mu(\nu, 0)\omega^{\mu\lambda}_\alpha, \omega^{\mu\lambda}_\beta\rangle=E^{|\mu|}_{\alpha\beta}(R_\lambda^t\nu, 0)=(2\pi)^{\frac{n}{2}}\Phi_{\alpha\beta}(w).
	\end{align*}
	Comparing with Proposition \ref{specDecomp} gives $\cJ_{\nu, \lambda}(\mu)\circ{{\rm Pr}}_{l}=i^k\,{{\rm Pr}}_{l+k}\circ\pi_\mu(\nu, 0)\circ{{\rm Pr}}_{l}$ for all $l\in \N$.
\end{proof}
\begin{corollary}
	With the hypotheses of Proposition \ref{specDecomp},
	\begin{align*}
		\cJ_{\nu, \lambda}(\mu)^*\omega_{\gamma}^{\mu\lambda}=(2\pi)^{\frac{n}{2}}i^k\sum_{|\beta|=|\gamma|-k}\Phi_{\gamma\beta}(w)\omega_\beta^{\mu\lambda}.
	\end{align*}
	In particular $\cJ_{\nu, \lambda}(\mu)^*: \cE_{l+k}\to \cE_{l}$.
\end{corollary}
\begin{proof}
	By the explicit formula for the special Hermite functions, $\overline{\Phi_{\beta\gamma}(w)}=(-1)^{|\gamma-\beta|}\Phi_{\gamma\beta}(w)$, so
	\begin{align*}
		\langle \cJ_{\nu, \lambda}(\mu)^*\omega_\gamma^{\mu\lambda}, \omega_\beta^{\mu\lambda}\rangle
		&=\overline{\langle\cJ_{\nu, \lambda}(\mu)\omega_\beta^{\mu\lambda}, \omega_\gamma^{\mu\lambda}\rangle}\\
		&=(2\pi)^{\frac{n}{2}}(-i)^k\,\overline{\Phi_{\beta\gamma}(w)}\,\delta(k+|\beta-\gamma|)
		=(2\pi)^{\frac{n}{2}}i^k\,\Phi_{\gamma\beta}(w)\,\delta(k+|\beta-\gamma|).\qedhere
	\end{align*}
\end{proof}

When $\mu\not\parallel\lambda$, the operator $R_\mu^tR_\lambda$ is no longer unitary, so it doesn't correspond to a change of basis in $\cE_l$. Instead we will make a symplectic transformation to put the curve $R_\mu^te^{sJ_\lambda}\nu$ into a normal form.

\begin{lemma}\label{l:normalForm} Let $\mu, \lambda\in \fz^*\mz$, $c=\hat\mu\cdot\widehat \lambda$, $\nu\in S(\fv^*)$, and $w(s):=R_\mu^tR_\lambda e^{sJ}R_\lambda^t\nu \in \R^{2n}$. Let
    \begin{align*}
        u_1:=R_\mu^t\nu, \quad u_2:=R_\mu^tJ_{\hat{\lambda}}\nu.
    \end{align*}

    If $c\neq 0$ then there exists $S\in Sp(2n, \R)$ such that $S(\sqrt{|c|}\, e_1)=u_1$ and $S(\sqrt{|c|}\, Je_1)=\sgn(c)u_2$. Any such $S$ satisfies
    \begin{align*}
        S^{-1}w(s)=\sqrt{|c|}\, e^{\sgn(c)sJ}e_1.
    \end{align*}
    If $c=0$ then there exists $S\in Sp(2n, \R)$ such that $S(Je_1)=u_1$ and $S(Je_2)=u_2$. Any such $S$ satisfies
    \begin{align*}
        S^{-1}w(s)=(\cos s\, Je_1+\sin s\, Je_2).
    \end{align*}
\end{lemma}
\begin{remark}\label{r:kahlerAngle}
	The number $c=\hat\mu\cdot\hat\lambda$ is the cosine of the angle between $\mu$ and $\lambda$ in $\fz^*$. Equivalently, it is the cosine of the K\"{a}hler angle, with respect to the complex structure $J_{\hat\mu}$, of the plane $P_\nu:={\rm Span}\{\nu, J_{\hat\lambda}\nu\}\subset\fv$ containing the horizontal projection of $\gamma_{\nu, \lambda}$. For the oriented orthonormal basis $(\nu, J_{\hat\lambda}\nu)$,
	\begin{align*}
		\omega_{\hat\mu}(\nu, J_{\hat\lambda}\nu)=\langle J_{\hat\mu}\nu, J_{\hat\lambda}\nu\rangle=\hat\mu\cdot\hat\lambda.
	\end{align*}
	By the $H$-type identity \ref{e:htypePolar} this angle does not depend on $\nu$. The three cases of Lemma \ref{l:normalForm} are the cases in which $P_\nu$ is $J_{\hat\mu}$-invariant ($c=\pm1$), symplectic with $\omega_{\hat\mu}|_{P_\nu}$ equal to $|c|$ times the area form ($0<|c|<1$), and totally real ($c=0$). On $\Gamma_\lambda^*$ with $|\lambda|=1$ one has $c=2k/|\mu|$, so the K\"{a}hler angle is quantized on each slice $\mu\cdot\lambda=2k$.
\end{remark}
\begin{proof}
    Let $A:=R_\mu^tR_\lambda$. Since $R_\lambda JR_\lambda^t=J_{\hat{\lambda}}$ by \eqref{normalForm}, we have $u_1=AR_\lambda^t\nu$ and $u_2=AJR_\lambda^t\nu$. In particular $u_1, u_2$ are unit vectors and
    \begin{align*}
        w(s)=(\cos{s})\, u_1+(\sin s)\, u_2.
    \end{align*}
    Let $\omega_{\R^{2n}}(u, v):=v^tJu$. Then
    \begin{align*}
        \omega_{\R^{2n}}(u_1, u_2)
        &=-\nu^tR_\lambda J A^tJ(R_\mu^t\nu)\\
        &=-\nu^t J_{\hat\mu} J_{\hat\lambda} \nu\\
        &=\langle J_{\hat\mu} \nu, J_{\hat\lambda} \nu\rangle.
    \end{align*}
    Using $J_{\hat\mu} J_{\hat\lambda}+J_{\hat\lambda} J_{\hat\mu}=-2(\hat{\mu} \cdot \hat{\lambda}) I$ and $|\nu|=1$, we obtain
    \begin{align*}
        \omega(u_1, u_2)=\hat{\mu}\cdot \hat{\lambda}=c.
    \end{align*}
    Thus $P_0:={\rm Span}\{ u_1, u_2\}$ is a symplectic subspace of $\R^{2n}$ when $c\neq 0$. Furthermore $\omega(e_1, Je_1)=(Je_1)^tJe_1=1$.

    For $c>0$,
    \begin{align*}
        \omega(\sqrt{c}\, e_1, \sqrt{c}\, Je_1)=c=\omega(u_1, u_2),
    \end{align*}
    so the map $S(\sqrt{c}\, e_1)=u_1$, $S(\sqrt{c}\, Je_1)=u_2$ extends to an element $S\in Sp(2n, \R)$. For any such $S$,
    \begin{align*}
        S\left(\sqrt{c}\, e^{sJ}e_1\right)=(\cos s)\, u_1+(\sin s)\, u_2=w(s).
    \end{align*}
    For $c<0$,
    \begin{align*}
        \omega(\sqrt{|c|}\, e_1, \sqrt{|c|}\, Je_1)=|c|=\omega(u_1, -u_2),
    \end{align*}
    so the map $S(\sqrt{|c|}\, e_1)=u_1$, $S(\sqrt{|c|}\, Je_1)=-u_2$ extends to an element $S\in Sp(2n, \R)$. For any such $S$,
    \begin{align*}
        S\left(\sqrt{|c|}\, e^{-sJ}e_1\right)=(\cos s)\, u_1+(\sin s)\, u_2=w(s).
    \end{align*}
    When $c=0$, the subspace $P_0\subset \R^{2n}$ is isotropic. Since $\langle u_1, u_2\rangle=0$, they are linearly independent. Since $\{Je_1, Je_2\}$ and $\{u_1, u_2\}$ each span a two dimensional isotropic subspace of $\R^{2n}$, the map $S(Je_1)=u_1$, $S(Je_2)=u_2$ extends to an element $S\in Sp(2n, \R)$. For any such $S$, $S((\cos s\, Je_1+\sin s\, Je_2))=w(s)$.
\end{proof}

\begin{remark}
    Lemma \ref{l:normalForm} may be read at the level of complex structures. Let $K:=R_\mu^tJ_{\hat{\lambda}}R_\mu$, so that $K\in O(2n)$, $K^2=-I$, and
    \begin{align*}
        w(s)=e^{sK}R_\mu^t\nu, \quad Ku_1=u_2.
    \end{align*}
    Since $J=R_\mu^tJ_{\hat{\mu}}R_\mu$, the relation $J_{\hat{\mu}}J_{\hat{\lambda}}+J_{\hat{\lambda}}J_{\hat{\mu}}=-2c\, I$ gives
    \begin{align*}
        \omega(Kx, Ky)+\omega(x, y)=2c\,\langle Kx, y\rangle, \quad x, y\in \R^{2n},
    \end{align*}
    so $K$ preserves $\omega$ if and only if $J=cK$, that is $|c|=1$, in which case $K=cJ$. Since the symplectic conjugates of $J$ are exactly the $\omega$-compatible complex structures, for $0<|c|<1$ there is no $S\in Sp(2n, \R)$ with $S^{-1}KS=J$. Lemma \ref{l:normalForm} instead conjugates the single $K$-orbit through $R_\mu^t\nu$ to a $J$-orbit, at the cost of the factor $\sqrt{|c|}$. For $c=0$, $K$ anti-preserves $\omega$ and the orbit spans the isotropic plane $P_0$.
\end{remark}

\begin{lemma}\label{l:Jconjugate} Let $\lambda, \mu\in \fz^*\setminus\{0\}$ and $\nu\in S(\fv^*)$ with $\mu\cdot\lambda =2k|\lambda|^3$ for some $k\in \mathbb{Z}$.
    Suppose $\lambda\in S(\fz^*)$. The operator $\cJ_{\nu, \lambda}(\mu)$ is unitarily equivalent to
    \begin{align*}
        \frac{i^k}{2\pi}\int_0^{2\pi}e^{iks}\pi_{|\mu|}(w(s), 0)ds
    \end{align*}
    such that $w(s)$ is the curve given for $c=\hat\mu\cdot\hat\lambda>0$ by
    \begin{align*}
        w(s)=\sqrt{c}e^{sJ}e_1
    \end{align*}
    and $c<0$ by
    \begin{align*}
        w(s)=\sqrt{|c|}e^{-sJ}e_1
    \end{align*}
    and for $c=0$ by
    \begin{align*}
        w(s)=\cos s\, Je_1+\sin s\, Je_2.
    \end{align*}
\end{lemma}
\begin{proof}
Since $J^{-1}e^{sJ}=e^{(s-\pi/2)J},$ the substitution $s\mapsto s+\pi/2$ in Definition \ref{d:J} gives
\begin{align*}
    \cJ_{\nu, \lambda}(\mu)=\frac{i^k}{2\pi}\int_0^{2\pi}e^{iks}\pi_{\mu}\left(R_\lambda e^{sJ}R_\lambda^t\nu, 0\right)ds.
\end{align*}
Let $w(s):=R_\mu^tR_\lambda e^{sJ}R_\lambda^t\nu$ be as in Lemma \ref{l:normalForm}. Then by equation \eqref{e:repH}
	\begin{align*}
	\pi_\mu\left(R_\lambda e^{sJ}R_\lambda^{{\rm{t}}}\nu, 0\right)
		=\pi_{|\mu|}\left(w(s), 0\right).
	\end{align*}
    For any $S\in Sp(2n, \R)$, let $\tau(S)$ be a unitary operator satisfying \eqref{e:intertwine}, which exists and is unique up to a phase. Then
     \begin{align*}
        \tau(S)^*\circ \cJ_{\nu, \lambda}(\mu)\circ\tau(S)=\frac{i^k}{2\pi}\int_0^{2\pi}e^{iks}\pi_{|\mu|}(S^{-1}w(s), 0)ds.
    \end{align*}
    By Lemma \ref{l:normalForm}, we are then free to choose $S$ to put $w(s)$ into the normal form above.
\end{proof}

\begin{proposition}[Spectral Decomposition of $\cJ_{\nu, \lambda}(\mu)$, $\mu\cdot\lambda\neq 0$]\label{specDecompShift}
Let $\lambda\in S(\fz^*)$, $\nu\in S(\fv^*)$, and $\mu\in \Gamma_\lambda^*\mz$ with $\mu\cdot\lambda=2k$, $k\in \Z^*$. Let $S\in Sp(2n, \R)$ be given by Lemma \ref{l:normalForm}, and $\tau(S)$ as in Lemma \ref{l:Jconjugate}. Then for all $\alpha\in \N^n$,
\begin{align}\label{e:shift}
	\tau(S)^*\circ\cJ_{\nu, \lambda}(\mu)\circ\tau(S)\,\omega_\alpha=\sqrt{2\pi}\, i^k\, \Phi_{\alpha_1, \alpha_1+|k|}\big(\sqrt{2|k|}\big)\,\omega_{\alpha+|k|e_1}
\end{align}
where $e_1:=(1, 0, \ldots, 0)\in \N^n$ and $\Phi_{ab}$ are the special Hermite functions given by \eqref{specialHermite}. In particular, the singular values of $\cJ_{\nu, \lambda}(\mu)$ are $\sigma_{\alpha_1}$, $\alpha\in \N^n$, where for $a\in \N$
\begin{align*}
	\sigma_a:=\sqrt{\tfrac{a!}{(a+|k|)!}}\,\big(\sqrt{|k|}\big)^{|k|}\,\big|L_a^{(|k|)}\big(|k|\big)\big|\, e^{-|k|/2},
\end{align*}
and $\cJ_{\nu, \lambda}(\mu)$ is injective if and only if $L_a^{(|k|)}(|k|)\neq 0$ for all $a\in \N$.
\end{proposition}
\begin{proof}
Write $\cJ':=\tau(S)^*\circ\cJ_{\nu, \lambda}(\mu)\circ\tau(S)$, $c:=\hat{\mu}\cdot\hat{\lambda}$, and $\varrho:=\sqrt{2|k|}$. Since $|\lambda|=1$ we have $c=2k/|\mu|$, so ${\rm sgn}(c)={\rm sgn}(k)$ and
\begin{align*}
	\sqrt{|\mu|}\sqrt{|c|}=\varrho.
\end{align*}
Suppose $k>0$. By Lemma \ref{l:Jconjugate},
\begin{align*}
	\cJ'=\frac{i^k}{2\pi}\int_0^{2\pi}e^{iks}\,\pi_{|\mu|}\big(\sqrt{c}\, e^{sJ}e_1, 0\big)\, ds.
\end{align*}
Since $B\Phi_\alpha=\omega_\alpha$ for the Bargmann transform \eqref{Bargmann}, we have $\langle \pi_{|\mu|}(z, 0)\omega_\alpha, \omega_\beta\rangle=E^{|\mu|}_{\alpha\beta}(z, 0)$ for $z\in \R^{2n}$. Then by Lemma \ref{cylinSym},
\begin{align*}
	\langle \cJ'\omega_\alpha, \omega_\beta\rangle
	&=\frac{i^k}{2\pi}\int_0^{2\pi}e^{i(k+|\alpha-\beta|)s}\, ds\; E^{|\mu|}_{\alpha\beta}\big(\sqrt{c}\, e_1, 0\big)\\
	&=i^k\,\delta(k+|\alpha-\beta|)\,(2\pi)^{\frac{n}{2}}\,\Phi_{\alpha\beta}(\varrho e_1).
\end{align*}
When $k<0$, Lemma \ref{l:Jconjugate} gives the orbit $e^{-sJ}$ and the same computation yields the factor $\delta(k-|\alpha-\beta|)$. In both cases $\langle \cJ'\omega_\alpha, \omega_\beta\rangle=0$ unless $|\beta|=|\alpha|+|k|$.

Since $\Phi_{\alpha\beta}(z)=\prod_{j=1}^n\Phi_{\alpha_j\beta_j}(z_j)$ and $\Phi_{ab}(0)=(2\pi)^{-\frac{1}{2}}\delta_{ab}$,
\begin{align*}
	\Phi_{\alpha\beta}(\varrho e_1)=(2\pi)^{-\frac{n-1}{2}}\,\Phi_{\alpha_1\beta_1}(\varrho)\prod_{j=2}^n\delta_{\alpha_j\beta_j}.
\end{align*}
Together with $|\beta|=|\alpha|+|k|$ this forces $\beta=\alpha+|k|e_1$, and \eqref{e:shift} follows.

Since the vectors $\omega_{\alpha+|k|e_1}$ are orthonormal, \eqref{e:shift} gives
\begin{align*}
	(\cJ')^*\cJ'\,\omega_\alpha=2\pi\,\big|\Phi_{\alpha_1, \alpha_1+|k|}(\varrho)\big|^2\,\omega_\alpha
\end{align*}
and by the explicit formula for the special Hermite functions, with $\tfrac{\varrho}{\sqrt{2}}=\sqrt{|k|}$ and $\tfrac{\varrho^2}{2}=|k|$,
\begin{align*}
	\sqrt{2\pi}\,\big|\Phi_{a, a+|k|}(\varrho)\big|=\sqrt{\tfrac{a!}{(a+|k|)!}}\,\left(\tfrac{\varrho}{\sqrt{2}}\right)^{|k|}\big|L_a^{(|k|)}\big(\tfrac{\varrho^2}{2}\big)\big|\, e^{-\varrho^2/4}=\sigma_a.
\end{align*}
Since $\tau(S)$ is unitary, $\cJ_{\nu, \lambda}(\mu)^*\cJ_{\nu, \lambda}(\mu)$ is unitarily equivalent to $(\cJ')^*\cJ'$, and $\cJ_{\nu, \lambda}(\mu)$ is injective if and only if $\sigma_a>0$ for all $a\in \N$.
\end{proof}
\begin{proposition}[Spectral Decomposition of $\cJ_{\nu, \lambda}(\mu)$, $\mu\cdot\lambda= 0$]\label{p:specDecompPerp}
    Let $\lambda\in S(\fz^*)$, $\nu\in S(\fv^*)$, and $\mu\in \Gamma_\lambda^*\mz$ with $\mu\cdot\lambda=0$. Let $S\in Sp(2n, \R)$ be given by Lemma \ref{l:normalForm}, $\tau(S)$ as in Lemma \ref{l:Jconjugate}, and let $B$ be the Bargmann transform \eqref{Bargmann}. Then
    \begin{align*}
        B^*\circ \tau(S)^*\circ\cJ_{\nu, \lambda}(\mu)\circ\tau(S)\circ B
    \end{align*}
    is multiplication on $L^2(\R^n)$ by
	\begin{align*}
		m(\xi)=\frac{1}{2\pi}\int_0^{2\pi}e^{i\sqrt{|\mu|}(\xi_1\cos s+\xi_2\sin s)}\, ds=J_0\big(\sqrt{|\mu|}\sqrt{\xi_1^2+\xi_2^2}\big)
	\end{align*}
	Therefore $\cJ_{\nu, \lambda}(\mu)$ is self adjoint and injective for every $\nu\in S(\fv^*)$ and every $\mu\in \fz^*\mz$ with $\mu\cdot \lambda=0$.
\end{proposition}
\begin{proof}
	Since $\mu\cdot\lambda=0$ we have $k=0$ and $c=\hat{\mu}\cdot\hat{\lambda}=0$, so Lemma \ref{l:Jconjugate} gives
	\begin{align*}
		\tau(S)^*\circ\cJ_{\nu, \lambda}(\mu)\circ\tau(S)=\frac{1}{2\pi}\int_0^{2\pi}\pi_{|\mu|}(w(s), 0)\, ds, \quad w(s)=\cos s\, Je_1+\sin s\, Je_2.
	\end{align*}
	Since $B$ intertwines $\rho_{|\mu|}$ and $\pi_{|\mu|}$,
	\begin{align*}
		B^*\circ\tau(S)^*\circ\cJ_{\nu, \lambda}(\mu)\circ\tau(S)\circ B=\frac{1}{2\pi}\int_0^{2\pi}\rho_{|\mu|}(w(s), 0)\, ds.
	\end{align*}
	In the coordinates $(x,y)\in\R^{2n}$, the curve $w(s)=(x(s),y(s))$ has
\begin{align*}
    x(s)=0, \qquad y(s)=(\cos s,\sin s,0,\ldots,0),
\end{align*}
so by the definition of the Schr\"{o}dinger representation
	\begin{align*}
		\rho_{|\mu|}(w(s), 0)\phi(\xi)=e^{i\sqrt{|\mu|}(\xi_1\cos s+\xi_2\sin s)}\phi(\xi).
	\end{align*}
	Averaging in $s$ gives multiplication by $m(\xi)$. Writing $\xi_1\cos s+\xi_2\sin s=\sqrt{\xi_1^2+\xi_2^2}\,\cos(s-\theta)$ and using the integral representation $J_0(t)=\frac{1}{2\pi}\int_0^{2\pi}e^{it\cos s}\, ds$ gives the Bessel form of $m$. Since $m$ is real, the multiplication operator is self adjoint. Its zero set $\{\xi\in \R^n: \sqrt{|\mu|}\sqrt{\xi_1^2+\xi_2^2}\in J_0^{-1}(0)\}$ is a countable union of cylinders over circles, hence has measure zero, so multiplication by $m$ is injective on $L^2(\R^n)$.
\end{proof}

\subsection{Normal Operator $\cN_{\nu, \lambda}$}
\begin{definition}
	For $\nu\in S(\fv^*)$ and $\lambda\in S(\fz^*)$ the normal operator is
	\begin{align*}
		\cN_{\nu, \lambda}(\mu)=\cJ_{\nu, \lambda}(\mu)^*\circ\cJ_{\nu, \lambda}(\mu), \quad \mu\in \Gamma_\lambda^*\mz.
	\end{align*}
\end{definition}
When $\mu$ and $\lambda$ are parallel, the spectrum of $\cN_{\nu, \lambda}(\mu)$ is easily computed.
\begin{proposition}[Spectral Decomposition of the Normal Operator]\label{p:normalParallel}
	Fix $\lambda\in S(\fz^*)$ and $\nu\in S(\fv^*)$ and $\mu\in \Gamma_\lambda^*\mz$ with $\mu\parallel\lambda$ and $\mu\cdot\lambda>0$. For all $\alpha\in \N^n$,
	\begin{align*}
	\cN_{\nu, \lambda}(\mu)\omega_\alpha^{\mu\lambda}
	&=(2\pi)^n(-1)^k\sum_{\substack{|\beta|=|\alpha|+k \\ |\gamma|=|\alpha|}}\Phi_{\alpha\beta}(w)\Phi_{\beta\gamma}(w)\omega_\gamma^{\mu\lambda}.
	\end{align*}
with $w=\sqrt{|\mu|}R_\lambda^t\nu$.
\end{proposition}

In particular, for such $\mu$, the normal operator is block-diagonal in the sense that $\cN_{\nu, \lambda}(\mu): \cE_l\to \cE_l$ for all $l\in \N$, with block diagonal terms given by
\begin{align*}
	\cN_{\nu, \lambda}(\mu)=\sum_{l=0}^\infty {{\rm Pr}}_{l}\circ \pi_\mu(-\nu, 0)\circ {{\rm Pr}}_{l+k}\circ\pi_\mu(\nu, 0)\circ {{\rm Pr}}_{l}.
\end{align*}

\begin{proof}
	With respect to the basis $\{\omega_\alpha\}_{\alpha\in \N^n}$,
	\begin{align*}
		\cN_{\nu, \lambda}(\mu)\omega_\alpha^{\mu\lambda}
		&=(2\pi)^{\frac{n}{2}}i^k\sum_{|\beta|=|\alpha|+k}\Phi_{\alpha\beta}(w)\cJ_{\nu, \lambda}(\mu)^*\omega_\beta^{\mu\lambda}\\
		&=(2\pi)^n(-1)^k\sum_{|\beta|=|\alpha|+k}\sum_{|\gamma|=|\beta|-k}\Phi_{\alpha\beta}(w)\Phi_{\beta\gamma}(w)\omega_\gamma^{\mu\lambda}\\
		&=(2\pi)^n(-1)^k\sum_{\substack{|\beta|=|\alpha|+k \\ |\gamma|=|\alpha|}}\Phi_{\alpha\beta}(w)\Phi_{\beta\gamma}(w)\omega_\gamma^{\mu\lambda}.\qedhere
	\end{align*}
\end{proof}
The following lemma will assist subsequent computations when $\mu \not\parallel\lambda$.
\begin{lemma}\label{l:sU}
	Let $\lambda\in S(\fz^*)$ and $\mu\in \fz^*\mz$ be not parallel with $\mu\cdot\lambda=2k$ for some $k\in \Z$. Write $\lambda=c\hat\mu+s{\hat\mu}_\perp$, where $c=2k/|\mu|$, $s=\sqrt{1-c^2}$, and ${\hat\mu}_\perp\in \fz^*$ is the unit vector along $\lambda-c\hat\mu$. Then $J_{\hat\mu}$ and $J_{{\hat\mu}_\perp}$ anticommute and generate a quaternionic structure on $\fv$. In particular $n=2d$ is even. The unitary group of this quaternionic structure,
	\begin{align*}
		Sp(d):=\{U\in O(\fv): UJ_{\hat\mu}=J_{\hat\mu}U \text{ and } UJ_{{\hat\mu}_\perp}=J_{{\hat\mu}_\perp}U\},
	\end{align*}
	acts transitively on $S(\fv^*)$, and for all $U\in Sp(d)$ and $\nu\in S(\fv^*)$,
	\begin{align}\label{e:sU}
		\cN_{U\nu, \lambda}(\mu)=\tau(\widetilde U)\circ\cN_{\nu, \lambda}(\mu)\circ\tau(\widetilde U)^*, \qquad \widetilde U:=R_\mu^tUR_\mu\in U(n),
	\end{align}
	with $\tau(\widetilde U)$ given by Lemma \ref{l:rotateAxis}.
\end{lemma}
\begin{proof}
	The polarized $H$-type condition gives $J_{\hat\mu}J_{{\hat\mu}_\perp}+J_{{\hat\mu}_\perp}J_{\hat\mu}=-2(\hat\mu\cdot{\hat\mu}_\perp)I=0$, so $J_{\hat\mu}$, $J_{{\hat\mu}_\perp}$, and $J_{\hat\mu}J_{{\hat\mu}_\perp}$ make $\fv$ a quaternionic vector space, and $2n=4d$. Every unit vector extends to a quaternionic orthonormal basis, and any two such bases are related by an element of $Sp(d)$, so $Sp(d)$ acts transitively on $S(\fv^*)$.

	Let $U\in Sp(d)$. Since $U$ is orthogonal and commutes with $J_{\hat\mu}$, equation \eqref{normalForm} gives $\widetilde UJ=J\widetilde U$, so $\widetilde U\in U(n)$. Since $\eta\mapsto J_\eta$ is linear, $U$ also commutes with $J_{\lambda}=cJ_{\hat\mu}+sJ_{{\hat\mu}_\perp}$.

	By \eqref{normalForm} and $|\lambda|=1$, the integrand in Definition \ref{d:J} is $e^{ik\theta}\pi_\mu(z(\theta, \nu), 0)$, where
	\begin{align*}
		z(\theta, \nu):=R_\lambda\frac{e^{\theta J}}{J}R_\lambda^{t}\nu=\frac{e^{\theta J_\lambda}}{J_\lambda}\nu, \quad \theta\in[0, 2\pi],
	\end{align*}
	so $z(\theta, U\nu)=Uz(\theta, \nu)$ because $U$ commutes with $J_\lambda$. By \eqref{e:repH}, $R_\mu^tU=\widetilde UR_\mu^t$, and Lemma \ref{l:rotateAxis}, for all $x\in \R^{2n}$,
	\begin{align*}
		\pi_\mu(Ux, 0)=\pi_{|\mu|}(\widetilde UR_\mu^tx, 0)=\tau(\widetilde U)\circ\pi_{|\mu|}(R_\mu^tx, 0)\circ\tau(\widetilde U)^*=\tau(\widetilde U)\circ\pi_\mu(x, 0)\circ\tau(\widetilde U)^*.
	\end{align*}
	Applying this with $x=z(\theta, \nu)$ and integrating in $\theta$,
	\begin{align*}
		\cJ_{U\nu, \lambda}(\mu)=\tau(\widetilde U)\circ\cJ_{\nu, \lambda}(\mu)\circ\tau(\widetilde U)^*.
	\end{align*}
	Taking adjoints and composing gives \eqref{e:sU}.
\end{proof}

\subsection{Averaged Normal Operator}

\leavevmode The full normal operator $I^*I$ is undefined when given Lebesgue measure on charges, since a geodesic of charge $\lambda$ spends length of order $|\lambda|$ in any fixed neighborhood of a point, and the integral over $\fz^*$ that defines its kernel diverges.

\begin{definition} For $\lambda\in \fz^*\mz$ and $\mu\in \Gamma_\lambda^*\mz$, define the averaged normal operator by
	\begin{align*}
		\cN_\lambda(\mu):=\int_{S(\fv^*)}\cN_{\nu, \lambda}(\mu)d\sigma(\nu)
	\end{align*}
Here $d\sigma(\nu)$ is the normalized measure on the unit sphere in $\fv^*$, and the integral converges in the strong operator topology since $\nu\mapsto\cN_{\nu, \lambda}(\mu)$ is strongly continuous with $\|\cN_{\nu, \lambda}(\mu)\|\leq 1$.
\end{definition}
\begin{remark}\label{r:scaleN}
	The assignment $(\lambda, \mu)\mapsto \cN_{\lambda}(\mu)$ is a field of operators indexed by $\mathfrak{R}=\{(\lambda, \mu)\in (\fz^*\mz) \times(\fz^*\mz): \mu\in \Gamma_\lambda^*\}$. Note that, $\cN_{\eps\lambda}(\eps^2\mu)=\cN_{\lambda}(\mu)$ for $\eps>0$ by Remark \ref{r:scaleJ}. Therefore $\cN$ is determined by its restriction to the set $\mathfrak{R}_1$ of all $(\lambda, \mu)\in \mathfrak{R}$ such that $|\lambda|=1$.
\end{remark}

\begin{proposition}[Diagonalization of $\cN_\lambda(\mu)$]\label{p:avgDiag}
	Let $\lambda\in S(\fz^*)$ and $\mu\in \Gamma_\lambda^*\mz$ with $\mu\cdot\lambda=2k$, $k\in \Z$. Then for every $l\in \N$, $\cE_l$ is an eigenspace of $\cN_\lambda(\mu)$, with eigenvalue
	\begin{align}\label{e:avgTrace}
		c_l=\binom{l+n-1}{l}^{-1}{\rm tr}\left({\rm Pr}_l\circ\cN_{\nu, \lambda}(\mu)\circ{\rm Pr}_l\right)
	\end{align}
	for any $\nu\in S(\fv^*)$.
\end{proposition}
\begin{proof}
	Fix $\nu\in S(\fv^*)$. Let
    $K:=Sp(d)$ when $\mu$ and $\lambda$ are not parallel. When $\mu\parallel\lambda$, instead let $K:=U(\lambda):=\{U\in O(\fv): UJ_\lambda=J_\lambda U\}$, the unitary group of the complex structure on $\fv$ induced by $J_\lambda$. For $U\in K$ we have $\widetilde U=R_\mu^tUR_\mu\in U(n)$ and $\cN_{U\nu, \lambda}(\mu)=\tau(\widetilde U)\circ\cN_{\nu, \lambda}(\mu)\circ\tau(\widetilde U)^*$. For $K=Sp(d)$ this is Lemma \ref{l:sU}. For $K=U(\lambda)$, every $U$ commutes with $J_{\hat\mu}=\pm J_{\hat\lambda}$, and the proof of Lemma \ref{l:sU} applies verbatim. In either case $K$ acts transitively on $S(\fv^*)$, so the pushforward of its Haar measure under $U\mapsto U\nu$ is the measure $d\sigma$ and
	\begin{align*}
		\cN_\lambda(\mu)=\int_{|\nu'|=1}\cN_{\nu', \lambda}(\mu)\,d\sigma(\nu')=\int_{K}\cN_{U\nu, \lambda}(\mu)\,dU.
	\end{align*}
	By left invariance of the Haar measure, for every $V\in K$,
	\begin{align*}
		\tau(\widetilde V)\circ\cN_\lambda(\mu)\circ\tau(\widetilde V)^*=\int_{K}\cN_{VU\nu, \lambda}(\mu)\,dU=\cN_\lambda(\mu),
	\end{align*}
	so $\cN_\lambda(\mu)$ commutes with $\tau(\widetilde U)$ for all $U\in K$. Each layer $\cE_l$ is preserved by these operators (Lemma \ref{l:rotateAxis}) and is irreducible under them. For $K=U(\lambda)$ the group $\{\widetilde U: U\in U(\lambda)\}$ is all of $U(n)$, whose action by $\tau$ on $\cE_l$ is the classical irreducible one. For $K=Sp(d)$ the layer carries the representation ${\rm Sym}^l$ of the standard representation of $Sp(d)$, which is irreducible \cite[Ch.~17]{fultonRepresentationTheory1991}. For $n\geq 2$ the layers have pairwise distinct dimensions $\dim\cE_l=\binom{l+n-1}{l}$, so they are pairwise inequivalent. For $n=1$, $K=U(\lambda)\cong U(1)$ and, by Definition \ref{d:tau}, $\tau(U)$ acts on the line $\cE_l=\C\zeta^l$ by the scalar $\bar{U}^{\,l}$, so the layers are the pairwise distinct characters $U\mapsto\bar{U}^{\,l}$ of $U(1)$. In either case Schur's lemma gives $\cN_\lambda(\mu)=\sum_{l}c_l\,{\rm Pr}_l$. Taking traces against ${\rm Pr}_l$,
	\begin{align*}
		c_l\dim\cE_l=\int_{K}{\rm tr}\left({\rm Pr}_l\circ\cN_{U\nu, \lambda}(\mu)\circ{\rm Pr}_l\right)dU={\rm tr}\left({\rm Pr}_l\circ\cN_{\nu, \lambda}(\mu)\circ{\rm Pr}_l\right),
	\end{align*}
	where the last equality uses that $\tau(\widetilde U)$ commutes with ${\rm Pr}_l$.
\end{proof}

\begin{proposition}[Eigenvalues of $\cN_\lambda(\mu)$]\label{p:avgEval}
	Under the hypotheses of Proposition \ref{p:avgDiag},
	\begin{align}\label{e:avgEval}
		c_l=\binom{l+n-1}{l}^{-1}\frac{1}{2\pi}\int_0^{2\pi}e^{ik(\theta-\sin\theta)}\,L_l^{(n-1)}\!\left(2|\mu|\sin^2\tfrac{\theta}{2}\right)e^{-|\mu|\sin^2\frac{\theta}{2}}\,d\theta.
	\end{align}
	In particular, $c_l$ depends only on $l$, $k$, and $|\mu|$.
\end{proposition}
\begin{proof}
	Let $\nu\in S(\fv^*)$ and set $\nu':=J_\lambda^{-1}\nu$. By \eqref{normalForm}, $R_\lambda \frac{e^{sJ}}{J}R_\lambda^t=e^{sJ_\lambda}J_\lambda^{-1}$, so Definition \ref{d:J} reads
	\begin{align*}
		\cJ_{\nu, \lambda}(\mu)=\frac{1}{2\pi}\int_0^{2\pi}e^{iks}\pi_\mu(x(s), 0)\,ds, \qquad x(s):=e^{sJ_\lambda}\nu',
	\end{align*}
	and therefore
	\begin{align*}
		\cN_{\nu, \lambda}(\mu)=\frac{1}{4\pi^2}\int_0^{2\pi}\!\!\int_0^{2\pi}e^{ik(s'-s)}\,\pi_\mu(x(s), 0)^*\pi_\mu(x(s'), 0)\,ds\,ds'.
	\end{align*}
	By the group law \eqref{product}, $\pi_\mu(x(s), 0)^*\pi_\mu(x(s'), 0)=e^{-\frac{i}{2}\omega_\mu(x(s), x(s'))}\pi_\mu(x(s')-x(s), 0)$. Write $\mu=2k\lambda+\mu_\perp$ with $\mu_\perp\cdot\lambda=0$, so $J_\mu=2kJ_\lambda+J_{\mu_\perp}$. Since $e^{s'J_\lambda}$ is orthogonal and commutes with $J_\lambda$,
	\begin{align*}
		\omega_\lambda(x(s), x(s'))=\langle J_\lambda e^{(s-s')J_\lambda}\nu', \nu'\rangle=\sin(s'-s),
	\end{align*}
	and by the anticommutation $J_{\mu_\perp}e^{sJ_\lambda}=e^{-sJ_\lambda}J_{\mu_\perp}$,
	\begin{align*}
		\omega_{\mu_\perp}(x(s), x(s'))=\langle e^{-(s+s')J_\lambda}J_{\mu_\perp}\nu', \nu'\rangle=0,
	\end{align*}
	since $J_{\mu_\perp}$ and $J_\lambda J_{\mu_\perp}$ are skew-symmetric. Thus, with $\theta:=s'-s$,
	\begin{align*}
		\pi_\mu(x(s), 0)^*\pi_\mu(x(s'), 0)=e^{-ik\sin\theta}\,\pi_\mu(x(s')-x(s), 0), \qquad |x(s')-x(s)|=2\left|\sin\tfrac{\theta}{2}\right|.
	\end{align*}
	By \eqref{e:repH} and \eqref{entryFunc}, ${\rm tr}({\rm Pr}_l\circ\pi_\mu(z, 0)\circ{\rm Pr}_l)=\sum_{|\alpha|=l}E^{|\mu|}_{\alpha\alpha}(R_\mu^tz, 0)$. By the explicit formula for the special Hermite functions and the identity $\sum_{|\alpha|=l}\prod_{j=1}^nL^{(0)}_{\alpha_j}(x_j)=L^{(n-1)}_l(x_1+\cdots+x_n)$ \cite{thangaveluHarmonicAnalysisHeisenberg1998},
	\begin{align*}
		{\rm tr}\left({\rm Pr}_l\circ\pi_\mu(z, 0)\circ{\rm Pr}_l\right)=L_l^{(n-1)}\!\left(\tfrac{|\mu|}{2}|z|^2\right)e^{-|\mu||z|^2/4}.
	\end{align*}
	The integrand depends only on $\theta$, so the double integral collapses to a single one, and
	\begin{align*}
		{\rm tr}\left({\rm Pr}_l\circ\cN_{\nu, \lambda}(\mu)\circ{\rm Pr}_l\right)=\frac{1}{2\pi}\int_0^{2\pi}e^{ik(\theta-\sin\theta)}\,L_l^{(n-1)}\!\left(2|\mu|\sin^2\tfrac{\theta}{2}\right)e^{-|\mu|\sin^2\frac{\theta}{2}}\,d\theta.
	\end{align*}
	Together with \eqref{e:avgTrace} this gives \eqref{e:avgEval}.
\end{proof}
Propositions \ref{p:avgDiag} and \ref{p:avgEval} are an instance of the formula for the Weyl transform of a radial function \cite[Theorem 1.4.3]{thangaveluHarmonicAnalysisHeisenberg1998}, since $\cN_\lambda(\mu)$ is the Weyl transform of a radial measure on $\fv$.

\begin{proposition}[Spectral Decomposition of $\cN_\lambda(\mu)$, $\mu\parallel\lambda$]\label{specDecompNormal} Let $\lambda\in S(\fz^*)$ and $\mu\in \Gamma_\lambda^*\mz$ with $\mu\parallel\lambda$ and $\mu\cdot\lambda=2k>0$. Then the eigenvalue of $\cN_\lambda(\mu)$ on $\cE_l$ given by Proposition \ref{p:avgDiag} is
	\begin{align}\label{SVDNormal}
		c_l=(2\pi)^{n}\tfrac{l!(n-1)!}{(l+n-1)!}\sum_{\substack{|\beta|=l \\ |\gamma|=l+k }}|\Phi_{\beta\gamma}(w)|^2
	\end{align}
where $w=\sqrt{|\mu|}R_\lambda^t \nu$, for any $\nu\in S(\fv^*)$. The right-hand side is independent of $\nu$ by \eqref{e:avgTrace}.
\end{proposition}
\begin{proof}
	Since $\hat{\mu}=\hat{\lambda}$, the operator $\tau(R_\mu^tR_\lambda)$ is unitary and preserves each layer (Lemma \ref{l:rotateAxis}), so $\{\omega^{\mu\lambda}_\beta: |\beta|=l\}$ is an orthonormal basis of $\cE_l$. By Proposition \ref{specDecomp},
	\begin{align*}
		{\rm tr}\left({\rm Pr}_l\circ\cN_{\nu, \lambda}(\mu)\circ{\rm Pr}_l\right)=\sum_{|\beta|=l}\left\|\cJ_{\nu, \lambda}(\mu)\,\omega^{\mu\lambda}_\beta\right\|^2=(2\pi)^n\sum_{\substack{|\beta|=l \\ |\gamma|=l+k}}|\Phi_{\beta\gamma}(w)|^2,
	\end{align*}
	and \eqref{e:avgTrace} gives \eqref{SVDNormal}.
\end{proof}
\begin{corollary}\label{invCriteria}
	Let $\lambda\in S(\fz^*)$ and $\mu\in \Gamma_\lambda^*\mz$ with $\mu\parallel\lambda$ and $\mu\cdot \lambda=2k$, and let $w=\sqrt{|\mu|}R_\lambda^t\nu$ for any $\nu\in S(\fv^*)$. The operator $\cN_\lambda(\mu): \cF(\C^n)\to \cF(\C^n)$ is injective if and only if for every $l\in \N$, there exist $\beta, \gamma\in \N^n$ with $|\beta|=l$ and $|\gamma|=l+|k|$ for which $\Phi_{\beta\gamma}(w)\neq 0$.
\end{corollary}
\begin{proof}
	By Proposition \ref{p:avgDiag}, $\cN_\lambda(\mu)=\sum_l c_l\,{\rm Pr}_l$ and the layers $\cE_l$ span $\cF(\C^n)$, so $\cN_\lambda(\mu)$ is injective if and only if $c_l\neq 0$ for every $l\in \N$. By Proposition \ref{p:avgEval}, $c_l$ depends only on $(l, k, |\mu|)$, and since $c_l$ is real, the substitution $\theta\mapsto-\theta$ in \eqref{e:avgEval} shows $c_l$ is unchanged under $k\mapsto-k$ at fixed $|\mu|$. So we may assume $k>0$ (parallelism and $\mu\neq0$ force $k\neq0$). Then the hypotheses of Proposition \ref{specDecompNormal} hold, and \eqref{SVDNormal} expresses $c_l$ as a positive multiple of $\sum_{|\beta|=l, |\gamma|=l+k}|\Phi_{\beta\gamma}(w)|^2$. The corollary follows.
\end{proof}

\begin{remark}[The first Heisenberg group]\label{r:nOne}
	 All statements but the diagonalization in Part 1 of Theorem \ref{t:stability} exclude $n=1$. For $n=1$ the average over $\nu$ is trivial, $\cN_\lambda(\mu)=\cJ_{\nu, \lambda}(\mu)^*\cJ_{\nu, \lambda}(\mu)$ by \eqref{e:reparam}, and each $\cE_l$ is a line, so by Proposition \ref{specDecompShift} $c_l=\sigma_l^2$ is the square of a singular value of $\cJ_{\nu, \lambda}(\mu)$,
	\begin{align*}
		c_l=\frac{l!}{(l+|k|)!}\,|k|^{|k|}e^{-|k|}\,L_l^{(|k|)}(|k|)^2,
	\end{align*}
	and Fej\'er's formula for Laguerre polynomials at a fixed point \cite[Theorem 8.22.1]{szegoOrthogonalPolynomials1975} gives
	\begin{align*}
		c_l=\frac{1+(-1)^k\sin\big(4\sqrt{|k|\,l}\big)}{2\pi\sqrt{|k|\,l}}+O(l^{-1}).
	\end{align*}
	Since $4\sqrt{|k|\,l}$ is dense modulo $2\pi$, $\liminf_l\sqrt{l}\,c_l=0$, so the loss exceeds half a derivative along a sparse sequence of layers and no estimate of the form \eqref{e:mainStab} holds. Whether an estimate with a larger loss holds depends on how closely $|k|$ approaches the zeros of $L_l^{(|k|)}$ as $l\to\infty$, a Diophantine question we leave open.
\end{remark}

 For $\mu\parallel\lambda$, $\cN_\lambda(\mu)$ is injective for every $k\neq0$ when $n>1$, and for odd $k$ when $n=1$.
\begin{proposition}\label{oddInteger}
	Fix $\lambda\in S(\fz^*)$, let $\mu=2k\lambda$ with $k\in \Z\mz$, and write $\cN_\lambda(\mu)=\sum_l c_l\,{\rm Pr}_l$ as in Proposition \ref{p:avgDiag}.
	\begin{itemize}
		\item[(a)] If $n =1$ the operator $\cN_\lambda(\mu):\cF(\C^n)\to \cF(\C^n)$ is injective for $k$ odd.
		\item[(b)] If $n>1$ then $c_l>0$ for every $l\in \N$.
	\end{itemize}
\end{proposition}
By Remark \ref{r:scaleN}, the normalization $|\lambda|=1$ entails no loss of generality: for $\lambda\in \fz^*\mz$ we have $\cN_\lambda(\mu)=\cN_{\hat\lambda}(|\lambda|^{-2}\mu)$ and $|\lambda|^{-2}\Gamma_\lambda^*=\Gamma_{\hat\lambda}^*$, so the proposition applies to all $\lambda\in \fz^*\mz$ after rescaling $\mu$.

Here we again see a contrast between the endpoint case $n=1$ which corresponds to first Heisenberg group $\heis_1$, and higher Heisenberg and $H$-type groups. When $n=1$ the layers $\cE_l$ are one dimensional so the normal operator is already radial and averaging gains nothing.

It is not surprising that injectivity is easier to achieve in the higher dimensional case. Since the space of geodesics grows like $\dim \cG= 2\dim G-2$ the problem of recoverability of an unknown function from the data of its X-ray transform becomes increasingly overdetermined as $\dim G$ increases.
\begin{proof}[Proof of Proposition \ref{oddInteger}]
	If $n=1$ then necessarily $m=1$ and hence $G=\heis_1$ is the standard Heisenberg group. Then the first part was proven in \cite{flynnInjectivityHeisenbergXray2021}.

	Suppose $n>1$ and write $\mu\cdot\lambda=2k$. Since $\cN_\lambda(\mu)$ is self-adjoint, $c_l$ is real, and the substitution $\theta\mapsto-\theta$ in \eqref{e:avgEval} shows that $c_l$ is unchanged under $k\mapsto-k$ at fixed $|\mu|$, so we may assume $k>0$.

	We bound the sum in \eqref{SVDNormal} below, where $w=\sqrt{|\mu|}R_\lambda^t\nu$ and $|w|^2=|\mu|$. Since $\langle\pi_1(z, 0)\omega_\beta, \omega_\gamma\rangle=E_{\beta\gamma}(z, 0)=(2\pi)^{\frac{n}{2}}\Phi_{\beta\gamma}(z)$ and, by Lemma \ref{l:rotateAxis}, $\pi_1(Uz, 0)=\tau(U)\circ\pi_1(z, 0)\circ\tau(U)^*$ with $\tau(U)$ unitary and preserving each layer, the sum in \eqref{SVDNormal} is invariant under $w\mapsto Uw$ for $U\in U(n)$. Choosing $U$ with $Uw=|w|e_1$ (Lemma \ref{l:rotateAxis}), we may assume $w=|w|e_1$. Since $\Phi_{\beta\gamma}(z)=\prod_{j=1}^n\Phi_{\beta_j\gamma_j}(z_j)$ and $\Phi_{ab}(0)=(2\pi)^{-\frac{1}{2}}\delta_{ab}$, only the terms with $\gamma=\beta+ke_1$ survive, and by the explicit formula following \eqref{specialHermite},
	\begin{align*}
		|\Phi_{\beta, \beta+ke_1}(|w|e_1)|^2=(2\pi)^{-n}\frac{\beta_1!}{(\beta_1+k)!}\left(\tfrac{|w|^2}{2}\right)^{k}L_{\beta_1}^{(k)}\!\left(\tfrac{|w|^2}{2}\right)^{2}e^{-|w|^2/2}.
	\end{align*}
	Counting the multi-indices $\beta$ with $|\beta|=l$ and $\beta_1=p$, Equation \eqref{SVDNormal} becomes, with $d_l:=\binom{l+n-1}{l}$,
	\begin{align*}
		c_l=\frac{1}{d_l}\sum_{p=0}^l\binom{l-p+n-2}{n-2}\frac{p!}{(p+k)!}\left(\tfrac{|w|^2}{2}\right)^{k}L_{p}^{(k)}\!\left(\tfrac{|w|^2}{2}\right)^{2}e^{-|w|^2/2}.
	\end{align*}
	Since $L_0^{(k)}\equiv1$, the $p=0$ term alone gives
	\begin{align*}
		c_l\geq\frac{1}{d_l}\binom{l+n-2}{n-2}\frac{1}{k!}\left(\tfrac{|\mu|}{2}\right)^{k}e^{-|\mu|/2}=\frac{n-1}{l+n-1}\,\frac{1}{k!}\left(\tfrac{|\mu|}{2}\right)^{k}e^{-|\mu|/2}>0.
	\end{align*}
	Every eigenvalue of $\cN_\lambda(\mu)$ is therefore positive; since the layers $\cE_l$ span $\cF(\C^n)$, $\cN_\lambda(\mu)$ is injective.
\end{proof}

\begin{proposition}[Injectivity of $\cN_\lambda(\mu)$]\label{p:Inj}
	Let $n>1$, $\lambda\in \fz^*\mz$ and $\mu\in \Gamma_\lambda^*\mz$. Then $c_l>0$ for every $l\in \N$, and $\cN_\lambda(\mu)$ is injective.
\end{proposition}
\begin{proof}
	By Remark \ref{r:scaleN} we may assume $|\lambda|=1$. Write $\mu\cdot\lambda=2k$ and fix $\nu\in S(\fv^*)$. By \eqref{e:avgTrace},
	\begin{align*}
		\tbinom{l+n-1}{l}\,c_l=\sum_{|\alpha|=l}\big\|\cJ_{\nu, \lambda}(\mu)\omega_\alpha\big\|^2,
	\end{align*}
	so $c_l=0$ if and only if $\cE_l\subset\ker\cJ_{\nu, \lambda}(\mu)$. If $k=0$, $\cJ_{\nu, \lambda}(\mu)$ is injective by Proposition \ref{p:specDecompPerp}. Let $k\neq0$, and let $S$ and $\tau(S)$ be as in Proposition \ref{specDecompShift}. By \eqref{e:shift}, as in the proof of Theorem \ref{t:singleDirection},
	\begin{align*}
		\ker\cJ_{\nu, \lambda}(\mu)=\tau(S)\cK, \qquad \cK:=\overline{\rm span}\{\omega_\alpha: \alpha_1\in \cA_{|k|}\}.
	\end{align*}
	Since $L_0^{(|k|)}=1$, $0\notin\cA_{|k|}$, so $\cK\perp\mathcal{W}_0:=\overline{\rm span}\{\omega_\alpha: \alpha_1=0\}$, and $\tau(S)\cK\perp\tau(S)\mathcal{W}_0$ since $\tau(S)$ is unitary. Hence
	\begin{align*}
		c_l=0\implies{\rm Pr}_l\,\tau(S)\mathcal{W}_0=0.
	\end{align*}
	Since $n\geq2$, $e^{t\zeta_2}=\sum_{p\in \N}t^p(p!)^{-1/2}\,\omega_{pe_2}\in \mathcal{W}_0$ for $t\in \C$. By Lemma \ref{l:tauGauss}, applied with $h=|\mu|$, $Sz=Mz+N\bar z$ and $\Theta=\bar NM^{-1}$, and with $b:=(M^t)^{-1}e_2\neq0$ and $\vartheta(\zeta):=\frac12\zeta^t\Theta\zeta$,
	\begin{align*}
		\tau(S)e^{t\zeta_2}=C_{te_2}\, e^{\vartheta+t\, b\cdot\zeta}.
	\end{align*}
	Since $\omega_\alpha=\zeta^\alpha/\sqrt{\alpha!}$, ${\rm Pr}_lF$ is the homogeneous part of degree $l$ of the Taylor series of $F$ at $0$, so
	\begin{align*}
		{\rm Pr}_l\,\tau(S)e^{t\zeta_2}=C_{te_2}\sum_{2i+j=l}\frac{t^j}{i!\, j!}\,\vartheta^i\,(b\cdot\zeta)^j.
	\end{align*}
	If ${\rm Pr}_l\,\tau(S)\mathcal{W}_0=0$ then, as $C_{te_2}\neq0$, the $\cE_l$-valued polynomial $\sum_{2i+j=l}\frac{t^j}{i!j!}\vartheta^i(b\cdot\zeta)^j$ vanishes for every $t\in \C$, and so does its coefficient $(b\cdot\zeta)^l/l!$ of $t^l$, contradicting $b\neq0$. Hence $c_l>0$ for every $l\in \N$. Since the $\cE_l$ span $\cF(\C^n)$, $\cN_\lambda(\mu)$ is injective.
\end{proof}

\subsection{Proofs of the main results}

\begin{lemma}\label{zeroSlice}
	Fix $\lambda\in \fz^*\mz$. If $I_{\nu, \lambda}f=0$ for all $\nu\in S(\fv^*)$, then
	\begin{align*}
		0=\cN_\lambda(\mu)\circ\cF(f)(\mu)
	\end{align*}
for all $\mu\in \Gamma_\lambda^*\mz$.
\end{lemma}

\begin{proof}
	By the Fourier slice theorem
	\begin{align*}
		0=\cJ_{\nu, \lambda}(\mu)\circ\cF(f)(\mu).
	\end{align*}
	Applying $\cJ_{\nu, \lambda}(\mu)^*$, we obtain
	\begin{align*}
		0=\cN_{\nu, \lambda}(\mu)\circ\cF(f)(\mu).
	\end{align*}
	Integrating over $\nu\in S(\fv^*)$, we have the result.
\end{proof}

We now combine Lemma \ref{zeroSlice} with the injectivity results of the previous subsection.
\begin{lemma}\label{PW}
	Suppose $f\in L^1(G)$ and $I_{\nu, \lambda}f=0$ for all $\nu\in S(\fv^*)$ and some $\lambda\in \fz^*\mz$.
	\begin{itemize}
		\item[1.] If $m=1$ and $n=1$, then $\cF(f)(\mu)=0$ for all $\mu=2k|\lambda|\lambda$, for $k\in 2\Z+1$
		\item[2.] If $m=1$ and $n>1$, then $\cF(f)(\mu)=0$ for all $\mu=2k|\lambda|\lambda$ for $k\in \Z\mz$
		\item[3.] If $m>1$, then $\cF(f)(\mu)=0$ for all $\mu\in \fz^*\mz$ such that $\mu\cdot\lambda\in 2|\lambda|^3\Z$.
	\end{itemize}
\end{lemma}
\begin{proof}
	By Lemma \ref{zeroSlice}, $\cN_\lambda(\mu)\circ\cF(f)(\mu)=0$ for all $\mu\in \Gamma_\lambda^*\mz$, so $\cF(f)(\mu)=0$ whenever $\cN_\lambda(\mu)$ is injective, and by Remark \ref{r:scaleN} injectivity may be checked at $|\lambda|=1$. Parts 1 and 2 follow from Proposition \ref{oddInteger}: when $\dim\fz=1$ every $\mu\in \Gamma_\lambda^*\mz$ has the form $2k|\lambda|\lambda$, and $\cN_\lambda(\mu)$ is injective for $k$ odd when $n=1$, and for all $k\in \Z\mz$ when $n>1$.

	For Part 3, $\dim\fz\geq2$ forces $n=2d>1$ (Lemma \ref{l:sU}), so $\cN_\lambda(\mu)$ is injective for every $\mu\in \Gamma_\lambda^*\mz$ by Proposition \ref{p:Inj}.
\end{proof}

The argument of Lemma \ref{PW} is quantitative. Averaging the slice theorem over $\nu\in S(\fv^*)$ expresses the $L^2$ norm of the data at a fixed frequency $\mu$ as a weighted norm of $\cF(f)(\mu)$, the weights being the eigenvalues $c_l$ of $\cN_\lambda(\mu)$. We first determine how these weights decay. Recall from Proposition \ref{p:avgEval} that $c_l=c_l(k, r)$ depends only on $l$, $k=\tfrac12\mu\cdot\lambda$ and $r=|\mu|$ when $|\lambda|=1$, and note that $r\geq|\mu\cdot\lambda|=2|k|$.
\begin{lemma}[Asymptotics of the eigenvalues of $\cN_\lambda(\mu)$]\label{l:eigenAsymp}
	Let $n>1$, $k\in \Z$, $r\geq 2|k|$, $r>0$, and $N:=2l+n$. Then
	\begin{align}\label{e:eigenAsymp}
		c_l(k, r)=\frac{\Gamma(n)}{\sqrt{\pi}\,\Gamma(n-\tfrac12)}\,(rN)^{-1/2}\big(1+o(1)\big), \qquad l\to\infty.
	\end{align}
	Consequently, if $c_l(k, r)>0$ for every $l\in \N$, there are constants $0<C_1\leq C_2$, depending on $n$, $k$ and $r$, such that
	\begin{align}\label{e:eigenBounds}
		C_1(rN)^{-1/2}\leq c_l(k, r)\leq C_2(rN)^{-1/2}, \qquad l\in \N.
	\end{align}
	Here $rN=r(2l+n)$ is the eigenvalue of $|\lambda|^{-2}\Lambda_\mu$ on $\cE_l$ when $r=|\lambda|^{-2}|\mu|$, so \eqref{e:eigenAsymp} says that $\cN_\lambda(\mu)=K_n|\lambda|\,\Lambda_\mu^{-1/2}(1+o(1))$, with $K_n$ as in the paragraph after Theorem \ref{t:stability}.
\end{lemma}
\begin{proof}
	Set $\alpha:=n-1\geq1$, $x(\theta):=2r\sin^2\tfrac{\theta}{2}$, and $\cL_l(x):=\binom{l+\alpha}{l}^{-1}e^{-x/2}L_l^{(\alpha)}(x)$, so that \eqref{e:avgEval} reads $c_l=\frac{1}{2\pi}\int_0^{2\pi}e^{ik(\theta-\sin\theta)}\cL_l(x(\theta))\,d\theta$. The substitution $\theta\mapsto2\pi-\theta$ fixes $x(\theta)$ and conjugates the phase, so
	\begin{align*}
		c_l=\frac{1}{\pi}\int_0^{\pi}\cos\big(k(\theta-\sin\theta)\big)\,\cL_l(x(\theta))\,d\theta.
	\end{align*}

	\textit{Step 1: Hilb's formula.} By \cite[Theorem 8.22.4]{szegoOrthogonalPolynomials1975}, whose $N$ is our $N/2=l+\tfrac{\alpha+1}{2}$, and $\Gamma(l+\alpha+1)/l!=\Gamma(\alpha+1)\binom{l+\alpha}{l}$, for $0<x\leq 2r$,
	\begin{align*}
		\cL_l(x)=\Gamma(\alpha+1)\,\big(\tfrac{Nx}{2}\big)^{-\alpha/2}J_\alpha\big(\sqrt{2Nx}\big)+E_l(x),
	\end{align*}
	where $E_l(x)=\binom{l+\alpha}{l}^{-1}x^{-\alpha/2}\epsilon_l(x)$ and Szeg\H{o}'s remainder satisfies, for a constant $c>0$, $\epsilon_l(x)=x^{5/4}O(l^{\alpha/2-3/4})$ if $c/l\leq x\leq 2r$ and $\epsilon_l(x)=x^{\alpha/2+2}O(l^{\alpha})$ if $0<x\leq c/l$, uniformly in $x$. Since $\frac{l^\alpha}{\alpha!}\leq\binom{l+\alpha}{l}=\prod_{j=1}^\alpha\frac{l+j}{j}\leq\frac{(l+\alpha)^\alpha}{\alpha!}$, in the first range $|E_l(x)|\leq CN^{-\alpha/2}x^{5/4-\alpha/2}l^{-3/4}$, and $x^{5/4-\alpha/2}$ is bounded by $(2r)^{5/4-\alpha/2}$ if $\alpha\leq2$ and by $(c/l)^{5/4-\alpha/2}$ if $\alpha\geq3$, giving $|E_l(x)|\leq Cl^{-3/4}$ in either case; in the second range $|E_l(x)|\leq Cx^2\leq C(c/l)^2$. Hence $\sup_{0<x\leq2r}|E_l(x)|=O(l^{-3/4})$ and
	\begin{align*}
		c_l=\frac{\Gamma(n)}{\pi}\int_0^{\pi}\cos\big(k(\theta-\sin\theta)\big)\,\big(\tfrac{Nx(\theta)}{2}\big)^{-\alpha/2}J_\alpha\big(\sqrt{2Nx(\theta)}\big)\,d\theta+O(l^{-3/4}).
	\end{align*}

	\textit{Step 2: away from $\theta=0$.} Fix $\delta\in(0, \pi)$. Since $|J_\alpha(t)|\leq Ct^{-1/2}$ for $t\geq1$ and $\tfrac12Nx(\theta)\geq rN\sin^2\tfrac{\delta}{2}$ on $[\delta, \pi]$, the integrand is bounded there by $C_\delta N^{-\alpha/2-1/4}\leq C_\delta N^{-3/4}$, using $\alpha\geq1$. So the integral over $[\delta, \pi]$ is $O(N^{-3/4})$.

	\textit{Step 3: near $\theta=0$.} On $[0, \delta]$ substitute $t=\sqrt{2Nx(\theta)}=2\sqrt{rN}\sin\tfrac{\theta}{2}$, so $dt=\sqrt{rN}\cos\tfrac{\theta}{2}\,d\theta$, $\theta=\theta_N(t):=2\arcsin\big(t/2\sqrt{rN}\big)$, and
	\begin{align*}
		&\frac{\Gamma(n)}{\pi}\int_0^{\delta}\cos\big(k(\theta-\sin\theta)\big)\big(\tfrac{Nx}{2}\big)^{-\alpha/2}J_\alpha\big(\sqrt{2Nx}\big)d\theta\\
        &\qquad =\frac{2^\alpha\Gamma(n)}{\pi\sqrt{rN}}\int_0^{T_N}\frac{\cos\big(k(\theta_N(t)-\sin\theta_N(t))\big)}{\cos(\theta_N(t)/2)}\,t^{-\alpha}J_\alpha(t)\,dt,
	\end{align*}
	with $T_N=2\sqrt{rN}\sin\tfrac{\delta}{2}\to\infty$. For fixed $t$, $\theta_N(t)\to0$, so the integrand converges pointwise to $t^{-\alpha}J_\alpha(t)$, and it is dominated by $\cos(\delta/2)^{-1}\,t^{-\alpha}|J_\alpha(t)|$, which is integrable on $(0, \infty)$: $t^{-\alpha}J_\alpha(t)$ is bounded near $0$ and $O(t^{-\alpha-1/2})$ at infinity, with $\alpha+\tfrac12>1$. By dominated convergence and Weber's integral \cite[Eq. 10.22.43]{DLMF},
	\begin{align*}
		\int_0^{T_N}\frac{\cos\big(k(\theta_N(t)-\sin\theta_N(t))\big)}{\cos(\theta_N(t)/2)}\,t^{-\alpha}J_\alpha(t)\,dt\longrightarrow\int_0^\infty t^{-\alpha}J_\alpha(t)\,dt=\frac{\Gamma(\tfrac12)}{2^\alpha\,\Gamma(\alpha+\tfrac12)}.
	\end{align*}
	Combining the three steps, $c_l=\frac{\Gamma(n)}{\sqrt\pi\,\Gamma(n-\frac12)}(rN)^{-1/2}(1+o(1))+O(N^{-3/4})$, which is \eqref{e:eigenAsymp}.

	For \eqref{e:eigenBounds}, \eqref{e:eigenAsymp} gives $l_0$ such that $\tfrac12K(rN)^{-1/2}\leq c_l\leq 2K(rN)^{-1/2}$ for $l\geq l_0$, where $K>0$ is the constant in \eqref{e:eigenAsymp}. The upper bound therefore holds for every $r$, since $c_l\leq\|\cN_\lambda(\mu)\|\leq1\leq(r(2l_0+n))^{1/2}(rN)^{-1/2}$ for $l<l_0$. For the lower bound, the finitely many values $c_l$, $l<l_0$, are positive by hypothesis.
\end{proof}

\begin{proof}[Proof of Theorem \ref{t:stability}]
	By Remark \ref{r:scaleN}, $\cN_\lambda(\mu)=\cN_{\hat\lambda}(|\lambda|^{-2}\mu)$, with $(|\lambda|^{-2}\mu)\cdot\hat\lambda=2k$ and $|\lambda|^{-2}|\mu|=r$, so for Part 1 we may assume $|\lambda|=1$. The diagonalization in Part 1 is Propositions \ref{p:avgDiag} and \ref{p:avgEval}; $c_l\geq0$ because $\cN_\lambda(\mu)=\int\cJ_{\nu, \lambda}(\mu)^*\cJ_{\nu, \lambda}(\mu)\, d\sigma(\nu)\geq0$, and the $\cE_l$ are the eigenspaces of $d\pi_\mu(-\Delta_{\rm{sr}})$, with eigenvalue $|\mu|(2l+n)$. The rest of Part 1 is Lemma \ref{l:eigenAsymp} and Proposition \ref{p:Inj}. For Part 2, write $A:=\cF(f)(\mu)$. By the Fourier slice theorem (Theorem \ref{FourierSliceTheorem}), $\cF_\lambda(I_{\nu, \lambda}f)(\pi_\mu)=2\pi|\lambda|^{-1}\cJ_{\nu, \lambda}(\mu)\circ A$, and since $A$ is Hilbert--Schmidt and $\cJ_{\nu, \lambda}(\mu)$ is bounded,
	\begin{align*}
		\big\|\cJ_{\nu, \lambda}(\mu)\circ A\big\|_{HS}^2={\rm tr}\big(A^*\circ\cN_{\nu, \lambda}(\mu)\circ A\big)=\sum_{\alpha\in \N^n}\big\langle\cN_{\nu, \lambda}(\mu)A\omega_\alpha, A\omega_\alpha\big\rangle.
	\end{align*}
	Each term is nonnegative and measurable in $\nu$, so by Tonelli's theorem, the definition of $\cN_\lambda(\mu)$ as a strong integral, and Part 1,
	\begin{align*}
		\int_{S(\fv^*)}\big\|\cJ_{\nu, \lambda}(\mu)\circ A\big\|_{HS}^2\,d\sigma(\nu)={\rm tr}\big(A^*\circ\cN_{\lambda}(\mu)\circ A\big)=\sum_l c_l\,\big\|{\rm Pr}_l\circ A\big\|_{HS}^2,
	\end{align*}
	and by the Fourier slice theorem the integral in \eqref{e:mainStab} is $4\pi^2|\lambda|^{-2}$ times the left side. Since $\Lambda_\mu=|\mu|(2l+n)$ on $\cE_l$ and $r=|\lambda|^{-2}|\mu|$, $(rN)^{-1/2}=|\lambda|\,\Lambda_\mu^{-1/2}$ on $\cE_l$ with $N=2l+n$, and inserting \eqref{e:eigenBounds}, whose hypothesis holds by Proposition \ref{p:Inj}, with its constants multiplied by $4\pi^2/|\lambda|$, gives \eqref{e:mainStab}.
\end{proof}

\begin{proof}[Proof of Theorem \ref{t:singleDirection}]
	By the Fourier slice theorem (Theorem \ref{FourierSliceTheorem}),
\begin{align*}
    0=\cF_\lambda(I_{\nu,\lambda}f)(\pi_\mu)=2\pi|\lambda|^{-1}\cJ_{\nu,\lambda}(\mu)\circ\cF(f)(\mu)
\end{align*}
for all $\mu\in\Gamma_\lambda^*\mz$. By Remark \ref{r:scaleJ}, $\cJ_{\nu, \lambda}(\mu)=\cJ_{\nu, \hat\lambda}(|\lambda|^{-2}\mu)$, so we may assume $|\lambda|=1$ and $\mu\cdot\lambda=2k$. In particular ${\rm ran}\,\cF(f)(\mu)\subset\ker\cJ_{\nu, \lambda}(\mu)$.

	If $k=0$, the operator $\cJ_{\nu, \lambda}(\mu)$ is injective by Proposition \ref{p:specDecompPerp}, so $\cF(f)(\mu)=0$.

	If $k\neq0$, Proposition \ref{specDecompShift} gives $\tau(S)^*\circ\cJ_{\nu, \lambda}(\mu)\circ\tau(S)\,\omega_\alpha=\sigma_{\alpha_1}'\,\omega_{\alpha+|k|e_1}$ with $|\sigma_a'|=\sigma_a$, and $\sigma_a=0$ if and only if $L_a^{(|k|)}(|k|)=0$, that is, $a\in \cA_{|k|}$. Since $\alpha\mapsto\alpha+|k|e_1$ is injective and the $\omega_{\alpha+|k|e_1}$ are orthonormal, $\ker(\tau(S)^*\circ\cJ_{\nu, \lambda}(\mu)\circ\tau(S))=\overline{\rm span}\{\omega_\alpha: \alpha_1\in \cA_{|k|}\}$, that is, $\ker\cJ_{\nu, \lambda}(\mu)$ is the closed span of the orthonormal family $\{\tau(S)\omega_\alpha: \alpha_1\in \cA_{|k|}\}$, with $S\in Sp(2n, \R)$ given by Lemma \ref{l:normalForm} and $\tau(S)$ as in Lemma \ref{l:Jconjugate}. Finally, $L_a^{(k)}(k)\neq0$ for all $a\in \N$ when $k$ is odd \cite{flynnInjectivityHeisenbergXray2021}, so $\cA_{|k|}=\emptyset$ for odd $k$.
\end{proof}

\begin{proof}[Proof of Corollary \ref{chargeFrequency}]
	This is Lemma \ref{PW}.
\end{proof}
\begin{proof}[Proof of Theorem \ref{MainResult}]
	Write $D_Z:=\bigcup_{\lambda\in Z}\left(\Gamma_{\lambda, \rm odd}^*\cup \lambda^\perp\right)\subset\bigcup_{\lambda\in Z}\Gamma_\lambda^*$. Suppose first that $D_Z$ is dense in $\fz^*$ and $Z$ is a zero set of $If$ for some $f\in L^1(G)$. That is, for each $\lambda\in Z$ there is $\nu_\lambda\in S(\fv^*)$ with $I_{\nu_\lambda, \lambda}f=0$. For every $\mu\in D_Z\mz$ there is $\lambda\in Z$, such that $\mu\cdot\lambda=2k|\lambda|^3$ for $k$ odd or zero. If $k=0$ then $\cF(f)(\mu)=0$ by Theorem \ref{t:singleDirection}(a). If $k$ is odd then again $\cF(f)(\mu)=0$ by Theorem \ref{t:singleDirection}(a). Thus $\cF(f)$ vanishes on the dense subset $D_Z\mz$ of $\fz^*\mz$. By \eqref{e:Weyl} and the injectivity of $W_\mu$ on $L^1(\fv)$ (Section \ref{s:Weyl}), $f^\mu=0$ for every $\mu\in D_Z\mz$. Since $f\in L^1(G)$, the map $\mu\mapsto f^\mu$ is continuous from $\fz^*$ to $L^1(\fv)$ by dominated convergence. Thus $f^\mu=0$ for every $\mu\in\fz^*$, and Fourier uniqueness in the central variable gives $f=0$. Hence $Z$ is an injectivity set.

	Conversely, suppose that $\bigcup_{\lambda\in Z}\Gamma_\lambda^*$ is not dense in $\fz^*$, and let $B\subset\fz^*$ be an open ball disjoint from every $\Gamma_\lambda^*$, $\lambda\in Z$; since $0\in \Gamma_\lambda^*$, $0\notin B$. Let $f(x, u)=\phi(x)\psi(u)$ with $0\neq\phi\in C_c^\infty(\fv)$ and $\psi\in \mathcal{S}(\fz)$ such that $\widehat{\psi}(\mu)=\int_\fz\psi(u)e^{-i\mu\cdot u}du$ is nonzero and supported in $B$. Then $0\neq f\in L^1(G)$ and $f^\mu=\widehat{\psi}(\mu)\phi=0$ for every $\mu\in \Gamma_\lambda^*$, $\lambda\in Z$. By \eqref{e:twistedSlice}, $(I_{\nu, \lambda}f)^\mu=f^\mu\times_\mu\sigma^\mu_{\nu, \lambda}=0$ for all $\mu\in \Gamma_\lambda^*$ and all $\nu\in S(\fv^*)$, so $I_{\nu, \lambda}f\in L^1(G_\lambda)$ vanishes by Corollary \ref{c:uniquenessQuotient}. Thus $Z$ is a zero set of $If$ with $f\neq0$, and $Z$ is not an injectivity set.

	Finally, we show that if $\bigcup_{\lambda\in Z}\Gamma_\lambda^*$ is dense in $\fz^*$ then so is $D_Z$. If there are $\lambda_i\in Z$ with $|\lambda_i|\to0$, then $\Gamma_{\lambda_i, \rm odd}^*$ is a family of parallel hyperplanes with spacing $4|\lambda_i|^2$, so every open ball meets $\Gamma_{\lambda_i, \rm odd}^*$ for $i$ large. Otherwise $\delta:=\inf_{\lambda\in Z}|\lambda|>0$, and every $\mu\in \Gamma_\lambda^*\setminus\lambda^\perp$ satisfies $|\mu|\geq|\mu\cdot\hat\lambda|\geq2\delta^2$; hence on the ball $\{|\mu|<2\delta^2\}$ the set $\bigcup_{\lambda\in Z}\Gamma_\lambda^*$ coincides with the cone $\bigcup_{\lambda\in Z}\lambda^\perp$, which is therefore dense in that ball and, being invariant under dilations, dense in $\fz^*$. Hence if $\bigcup_{\lambda\in Z}\Gamma_\lambda^*$ is dense then $D_Z$ is dense and $Z$ is an injectivity set.
\end{proof}
\begin{proof}[Proof of Corollary \ref{c:injSets}]
	By Theorem \ref{MainResult} it suffices, for each of the two conditions on $Z$, to show that $\bigcup_{\lambda\in Z}\Gamma_\lambda^*$ is dense in $\fz^*$.

	Suppose first that $Z$ contains a sequence $\lambda_i\to0$ with $\lambda_i\neq0$. Fix $\mu_0\in \fz^*$ and $r>0$. The set $\Gamma_{\lambda_i}^*$ consists of the hyperplanes $\{\mu: \mu\cdot\hat\lambda_i=2k|\lambda_i|^2\}$, $k\in \Z$ (points when $\dim\fz=1$), which have spacing $2|\lambda_i|^2$. Choosing $k$ with $|\mu_0\cdot\hat\lambda_i-2k|\lambda_i|^2|\leq|\lambda_i|^2$, the point $\mu_0+(2k|\lambda_i|^2-\mu_0\cdot\hat\lambda_i)\hat\lambda_i$ lies in $\Gamma_{\lambda_i}^*$ at distance at most $|\lambda_i|^2$ from $\mu_0$, which is less than $r$ for $i$ large. Hence $\bigcup_i\Gamma_{\lambda_i}^*$ is dense in $\fz^*$.

	Suppose next that $\dim\fz>1$ and that the set $D:=\{\hat\mu\in S(\fz^*): \hat\mu^\perp\cap Z\neq\emptyset\}$ of unit normals to hyperplanes through $0$ meeting $Z$ is dense in $S(\fz^*)$. If $\hat\mu\in D$ and $\lambda\in \hat\mu^\perp\cap Z$, then $\mu\cdot\lambda=0$ for every $\mu\in \R\hat\mu$, so $\R\hat\mu\subset\lambda^\perp\subset\Gamma_\lambda^*$. Thus $\bigcup_{\lambda\in Z}\Gamma_\lambda^*$ contains the cone $\{t\hat\mu: t\geq0, \ \hat\mu\in D\}$. Given $\mu\neq0$, choose $\hat\mu_j\in D$ with $\hat\mu_j\to\hat\mu$; then $|\mu|\hat\mu_j\to\mu$, so the cone is dense in $\fz^*$.

	For the last claim of (ii), let $C=P\cap S(\fz^*)$ be a great circle, $P\subset\fz^*$ a two-dimensional subspace, and suppose $\{\hat\lambda: \lambda\in Z\}$ is dense in $C$. Given $\hat\mu_0\in S(\fz^*)$, the subspace $\hat\mu_0^\perp\cap P$ has dimension at least one, so there is $c_0\in C$ with $c_0\cdot\hat\mu_0=0$. Choose $\lambda_j\in Z$ with $\hat\lambda_j\to c_0$, so that $\hat\mu_0\cdot\hat\lambda_j\to0$, and let $\hat\mu_j$ be the normalized orthogonal projection of $\hat\mu_0$ onto $\hat\lambda_j^\perp$, defined for $j$ large. Then $\lambda_j\in \hat\mu_j^\perp\cap Z$ and $\hat\mu_j\to\hat\mu_0$. Hence $D$ is dense in $S(\fz^*)$ and the first part of (ii) applies.

	Finally, suppose $f\in L^1(G)$ and $I_{\nu, \lambda}f=0$ in $L^1(G_\lambda)$ for almost every $(\nu, \lambda)\in S(\fv^*)\times\fz^*$. By Fubini's theorem, the set $Z:=\{\lambda\in \fz^*\mz: I_{\nu, \lambda}f=0 \text{ for some } \nu\in S(\fv^*)\}$ has full Lebesgue measure in $\fz^*$, so every ball about $0$ meets $Z$ and $Z$ contains a sequence $\lambda_i\to 0$. By definition $Z$ is a zero set of $If$, and by (i) it is an injectivity set, so $f=0$.
\end{proof}

\begin{proof}[Proof of Theorem \ref{PW2}]
	Since $\cA_{|k|}=\emptyset$ for odd $k$, Theorem \ref{t:singleDirection} gives $\cF(f)(\mu)=0$ for all $\mu\in \Gamma_{\lambda, \rm odd}^*$ and all $\lambda\in Z$. Under the additional hypothesis and for $n>1$, Lemma \ref{PW}, parts 2 and 3, gives $\cF(f)(\mu)=0$ for all $\mu\in \Gamma_\lambda^*\mz$ and all $\lambda\in Z$.

	Part 1. Let $\dim\fz=1$, so $\fz^*\mz\cong\R^*$, $Z=\{\lambda: |\lambda|\in[1, 1+\eps]\}$, and $\Gamma_\lambda^*\mz=\{2k|\lambda|\lambda: k\in \Z\mz\}$. As $\lambda$ ranges over $Z$, for fixed $k>0$ the points $2k|\lambda|\lambda$ fill $\{\mu: |\mu|\in[2k, 2k(1+\eps)^2]\}$. For odd $k$, consecutive intervals $[2k, 2k(1+\eps)^2]$ and $[2(k+2), 2(k+2)(1+\eps)^2]$ overlap once $2k(1+\eps)^2\geq 2(k+2)$, that is $k\geq 2/(2\eps+\eps^2)$. Hence $\cF(f)(\mu)=0$ for all $|\mu|\geq 2k_0$, and $f\in L^1_\Omega(G)$ with $c=2k_0$. Under the additional hypothesis and for $n>1$, every $k>0$ is available; the intervals for $k$ and $k+1$ overlap once $k\geq 1/(2\eps+\eps^2)$, and $c=2k_1$.

	Part 2. For $\hat\mu\in S(\fz^*)$ let
	\begin{align*}
		w(\hat\mu):=\max_{\lambda\in Z}\hat\mu\cdot\lambda-\min_{\lambda\in Z}\hat\mu\cdot\lambda,
	\end{align*}
	the width of $Z$ in the direction $\hat\mu$. If $w(\hat\mu)=0$ then $Z$ lies in an affine hyperplane $\{\hat\mu\cdot\lambda=c\}$, so $w>0$ on $S(\fz^*)$. Since $w$ is continuous and $S(\fz^*)$ is compact, $\delta=\min_{\hat\mu}w(\hat\mu)>0$. Fix $\mu\in \fz^*\mz$ with $|\mu|\geq 4/\delta$ and write $\mu=|\mu|\hat\mu$. Since $Z$ is compact and connected, $\mu\cdot Z:=\{\mu\cdot\lambda: \lambda\in Z\}$ is a closed interval, of length $|\mu|\,w(\hat\mu)\geq |\mu|\delta\geq 4$, so it contains a point of $2(2\Z+1)$. In particular, there is $\lambda\in Z$ with $\mu\cdot\lambda=2k$ for some odd $k$, that is, $\mu\in \Gamma_{\lambda, \rm odd}^*$, and $\cF f(\mu)=0$. Under the additional hypothesis it suffices that $\mu\cdot Z$ contain a point of $2\Z$, which holds once $|\mu|\geq 2/\delta$, giving $c=2/\delta$.

	For the cap $Z=C_h(\lambda_0)$, write $\cos\alpha=1-h$, so that $Z=\{\lambda\in S(\fz^*): \lambda\cdot\lambda_0\geq\cos\alpha\}$ and $\alpha\leq\tfrac{\pi}{2}$. The cap is compact, connected, and not contained in any affine hyperplane. Since $w(\hat\mu)=w(-\hat\mu)$, it suffices to consider $\hat\mu$ at angle $\varphi\in [0, \tfrac{\pi}{2}]$ from $\lambda_0$. Maximizing and minimizing $\hat\mu\cdot\lambda$ over the cap gives
	\begin{align*}
		w(\hat\mu)=
		\begin{cases}
			1-\cos(\alpha+\varphi), & 0\leq \varphi\leq \alpha,\\
			\cos(\varphi-\alpha)-\cos(\varphi+\alpha)=2\sin\alpha\sin\varphi, & \alpha\leq \varphi\leq \tfrac{\pi}{2}.
		\end{cases}
	\end{align*}
	Both expressions are increasing in $\varphi$, so the minimal width $\delta$ is attained in the radial direction. $\delta=w(\lambda_0)=1-\cos\alpha=h$, so $c=4/h$, or $2/h$ under the additional hypothesis.
\end{proof}

\begin{proof}[Proof of Corollary \ref{c:centralSupport}]
	By Theorem \ref{PW2} and \eqref{e:Weyl}, $W_\mu(f^\mu)=0$ for $|\mu|>c$, so $f^\mu=0$ for $|\mu|>c$ by the injectivity of $W_\mu$ on $L^1(\fv)$ (Section \ref{s:Weyl}). For $\phi\in C_c(\fv)$, the function $\mu\mapsto \int_\fv f^\mu(x)\phi(x)dx$ is the Fourier transform of $u\mapsto \int_\fv f(x, u)\phi(x)dx\in L^1(\fz)$, supported in $K$, and so extends to an entire function on $\C^m$. It vanishes on the open set $\{|\mu|>c\}$, hence identically, and $f=0$.
\end{proof}

\section{Concluding Remarks}\label{s:concluding}
The sub-Riemannian geodesics are also of interest because of their relation to the compatible Riemannian and pseudo-Riemannian geodesics. We remark that these techniques may be used to study the X-ray transforms $I^\eps$ and $I^{i\eps}$ associated to the family of left-invariant Riemannian and pseudo-Riemannian metrics
\begin{align*}
	g_\eps^\pm:=\langle\cdot, \cdot\rangle_\fv\pm \eps^{-2}\langle \cdot, \cdot \rangle_\fz.
\end{align*}
The metric $g_\eps^-$ has signature $(2n,m)$ and is Lorentzian only when $m=1$. In the momentum coordinates $(\nu,\lambda)$, the geodesic Hamiltonians are
\begin{align*}
    H_\eps^\pm(\nu,\lambda)=\tfrac12|\nu|^2\pm\tfrac12\eps^2|\lambda|^2.
\end{align*}
The central term is a Casimir for the Lie--Poisson bracket on $\fg^*$. Hence the horizontal projections of the $g_\eps^\pm$-geodesics with fixed initial point and covector are independent of $\eps$. For geodesics of nonzero charge with nonzero central advance per turn, the parameterization \eqref{e:orbitStabMain} and the Fourier slice argument adapt to $g_\eps^\pm$, as in \cite{flynnInjectivityHeisenbergXray2021} for the Riemannian metrics on the Heisenberg group. The sub-Riemannian question is therefore a natural starting point to study other left-invariant geometries as a perturbation of the sub-Riemannian one. We will explore this in a subsequent work.

It is our goal to generalize this work to study X-ray tomography on sub-Riemannian manifolds with less symmetry. Much progress has been made studying classical inverse problems on Riemannian manifolds using tools from microlocal analysis. On such a manifold the fixed-charge transform is no longer a convolution, and the role of the group Fourier transform passes to a pseudodifferential calculus whose symbols are fields of operators on the unitary duals of the osculating groups \cite{fermaniankammererQuantizationFilteredManifolds}.

Since $H$-type groups are the osculating models of certain step-two sub-Riemannian manifolds, by the lifting theorem of Rothschild and Stein \cite{rothschildHypoellipticDifferentialOperators1976}, the identity $\cN_\lambda(\mu)=K_n|\lambda|\Lambda_\mu^{-1/2}(1+o(1))$ for $n>1$ of Theorem \ref{t:stability} is the computation from which a parametrix construction for the normal operator on such a manifold would start, in a calculus on filtered manifolds \cite{kammererGeometricInvarianceSemiclassical2023, fermaniankammererQuantizationFilteredManifolds}. For $n>1$, Theorem \ref{t:stability} gives the principal order and the sharp loss at each nonzero visible central frequency. A principal-symbol statement in that calculus requires in addition asymptotics uniform in the central frequency, and this belongs to the variable-coefficient problem. The horizontal lines, $\lambda=0$, are the opposite extreme. On the Heisenberg group $\heis_n$ their normal operator needs no periodization, being convolution with $|x|^{1-2n}\otimes\delta_\fz$. No central element fixes a line, so every central frequency is visible. Strichartz \cite[Example 2]{strichartzLpHarmonicAnalysis1991} computed the symbol of this convolution, which at central frequency $\mu$ acts on the $l$-th eigenspace of $d\pi_\mu(-\Delta_{\rm{sr}})$ by a constant multiple of $|\mu|^{-1/2}\binom{l+n-1}{l}^{-1}$ times the coefficient of $w^l$ in $(1-w)^{\frac12-n}(1+w)^{-\frac12}$. For $n\geq2$ this is comparable to $d\pi_\mu\big((-\Delta_{\rm{sr}})^{-1/2}\big)$ uniformly in $\mu$ and $l$, the estimate which for $\lambda\neq0$ we obtain at each fixed $\mu$. For $n=1$ the generating function is $(1-w^2)^{-1/2}$, the odd eigenspaces are annihilated, and this is the non-injectivity on $\heis_1$. As $\lambda\to0$ the spacing $2|\lambda|^2$ of $\Gamma_\lambda^*$ tends to zero and the geodesics of small charge approach this picture (Corollary \ref{chargeFrequency}).

\section{Acknowledgments}
The author acknowledges the partial support of The Leverhulme Trust for this work via Research Project
Grant 2020-037 \textit{Quantum Limits for sub-elliptic Operators}.

The author acknowledges the partial support granted by the European Union –
NextGeneration EU Project ``NewSRG - New directions in SubRiemannian
Geometry” within the Program STARS@UNIPD 2021.

The author acknowledges the support of the EPSRC Early Career Fellowship EP/V001760/1 (J. Galkowski). For the purpose of open access, the author has applied a Creative Commons Attribution (CC BY) licence to any Author Accepted Manuscript version arising from this submission.

Claude (Anthropic, 2026) was used to help edit the exposition and the write-up of proofs, to check computations and references, and to draw Figure~\ref{fig:visible}.

The author thanks V\'{e}ronique Fischer for many informative discussions about the harmonic analysis of $H$-type groups in relation to this project.
\appendix

\section{Momentum functions}\label{momentumFunctions}
In this section let $Q$ be a smooth manifold and $G$ be a Lie group, acting smoothly on $Q$ on the left. For $v\in \fg$, let $v_Q$ be its infinitesimal generator on $Q$. That is
\begin{align*}
	v_Q(q)=\ddt \exp(tv)\cdot q.
\end{align*} The function
\begin{align*}
	\mathbf{J}: T^*Q \to \fg^*, \quad \langle\mathbf{J}(\alpha_q), v\rangle=\langle \alpha_q, v_Q\rangle=P_{v_Q}(\alpha_q)
\end{align*}
is called the momentum function for the action of $G$ on $T^*Q$.

If $Q=G$ so that $G$ acts on its cotangent bundle, the momentum map is given by
\begin{align}
	\mathbf{J}(\alpha_g)=TR_{g}^*\,\alpha_g, \quad \alpha_g\in T_g^*G \label{momentum}
\end{align}
If now $G$ is an $H$-type group, in view of \eqref{product} and \eqref{momentum} we have, for $g=(x, u)$, $\xi=X+U\in \fg$ with $X\in \fv$ and $U\in \fz$, and $\alpha_g\in T_g^*G$, $TL^*_g\alpha_g=(\nu, \lambda)\in \fg^*$,
\begin{align*}
	\mathbf{J}(\alpha_g)(\xi)
	&=\alpha_g(TR_g\, \xi)\\
	&=\langle (\nu, \lambda), Ad_{g^{-1}}\xi\rangle\\
	&=\langle (\nu, \lambda), (X, U-\omega(x, X))\rangle\\
	&=\langle \nu, X\rangle +\langle \lambda, U\rangle-\langle \lambda, \omega(x, X)\rangle\\
	&=\langle \nu, X\rangle+\langle \lambda, U\rangle-\langle J_\lambda x, X\rangle.
\end{align*}
So in the left trivialization, the momentum map is given by
	\begin{align*}
		\mathbf{J}(x, u, \nu, \lambda)=(\nu-J_\lambda(x), \lambda)\in\mathbb{R}^{2n}\times\mathbb{R}^m
	\end{align*}
Here we used the fact that
\begin{align*}
	Ad_{(x', u')}(x, u)=(x, u+\omega(x', x))
\end{align*}
then
\begin{align*}
	\langle Ad^*_{g^{-1}}(\nu, \lambda), (x, u)\rangle
	&=\langle (\nu, \lambda), Ad_{g^{-1}}(x, u)\rangle\\
	&=\langle (\nu, \lambda), (x, u-\omega(x', x))\rangle\\
	&=\langle \nu, x\rangle+\langle \lambda, u\rangle -\langle J_\lambda(x'), x\rangle\\
	&=\langle \lambda, u\rangle+\langle \nu-J_\lambda(x')^\flat, x\rangle\\
	&=\langle (\nu-J_\lambda(x')^\flat, \lambda), (x, u)\rangle.
\end{align*}
Here we use $\langle\cdot, \cdot\rangle$ for both the dual pairing on $\fg^*\times\fg$ and for the inner product on $\fg$. We insert the flat $\flat$ or sharp $\sharp$ where necessary to keep the notation consistent.

Let $s\mapsto \varphi_s(x, u; \nu, \lambda)$, given by \eqref{geodesic} be the unique complete unit speed geodesic starting at the origin with initial momentum $(\nu, \lambda)$. Then
\begin{align*}
	\mathbf{J}(\varphi_s(x, u; \nu, \lambda))
	&=Ad^*_{(x, u)^{-1}}\mathbf{J}(\varphi_s(0; \nu, \lambda))\\
	&=Ad^*_{(x, u)^{-1}}(e^{sJ_\lambda}\nu-J_\lambda J_\lambda^{-1}\left(e^{sJ_\lambda}-1\right)\nu, \lambda)\\
	&=Ad^*_{(x, u)^{-1}}(\nu, \lambda)\\
	&=(\nu-J_\lambda(x)^\flat, \lambda)
\end{align*}

In particular, $(\nu-J_\lambda(x)^\flat, \lambda)$ is an integral of motion. Note that by inspecting Equation \eqref{geodesics} for the geodesic flow, the central axis of the (helical) geodesic $s\mapsto\varphi_s(x, u; \nu, \lambda)$, with $\lambda\neq 0$, lies above the point $C(x, u, \nu, \lambda)=x-J_\lambda^{-1}\nu$ in $\fv$. The point $C=x-J_\lambda^{-1}\nu$ is called the \textit{guiding center}. Then $$\mathbf{J}(\varphi_s(x, u; \nu, \lambda))=\left(-J_\lambda\left(C(x, u, \nu, \lambda)\right)^\flat,\; \lambda\right).$$
Therefore, when $\lambda\neq 0$, the guiding center is a constant of motion. We choose the geodesics $\gamma_{\nu, \lambda}$ given by Equation \eqref{param} so that the guiding center is zero. Indeed, a simple computation gives
\begin{align*}
	\mathbf{J}(\gamma_{\nu, \lambda}(s))=(0, \lambda).
\end{align*} %

We also have the momentum map for right translations:
\begin{align*}
	\mathbf{J}_R: T^*G\to \fg^*, \quad \mathbf{J}_R(\alpha_g)=TL_g^*(\alpha_g).
\end{align*}
Then by a similar computation,
$\mathbf{J}_R(\varphi_s(x, u; \nu, \lambda))=(e^{sJ_\lambda}\nu, \lambda)$. Thus when $\lambda=0$, $\nu\in \fv^*$ is a constant of motion, and when $\lambda\neq 0$, the equivalence class of $\nu$ in $\C P(\fv^*, J_\lambda)$ is a constant of motion.

\printbibliography

\end{document}